\documentclass[12pt]{amsart}
\usepackage[pagebackref]{hyperref}
\usepackage{amsmath}
\usepackage{amssymb}
\usepackage{amscd}
\usepackage{amsfonts}
\usepackage{graphicx}
\usepackage{psfrag}
\usepackage{bbm}
\usepackage{latexsym}
\usepackage[all]{xy}
\usepackage[T1]{fontenc}
\usepackage{selinput}
\newtheorem{theorem}{Theorem}[section]
\newtheorem{lemma}[theorem]{Lemma}
\newtheorem{corollary}[theorem]{Corollary}
\newtheorem{proposition}[theorem]{Proposition}

\newtheorem{exa}[theorem]{Example}
\newenvironment{example}{\begin{exa} \em}{\end{exa}}
\newtheorem{exas}[theorem]{Examples}

\newtheorem{defini}[theorem]{Definition}
\newenvironment{definition}{\begin{defini} \em}{\end{defini}}

\newtheorem{rema}[theorem]{Remark}
\newenvironment{remark}{\begin{rema} \em}{\end{rema}}

\newtheorem{remas}[theorem]{Remarks}

\newtheorem{question}[theorem]{Question}

\newtheorem{problem}[theorem]{Problem}

\newtheorem{conjecture}[theorem]{Conjecture}

\newenvironment{equationth}{\stepcounter{theorem}\begin{equation}}{\end{equation}}

\def\cL{\mathcal L}
\def\a{\alpha}
\def\e{\varepsilon}
\def\t{\theta}
\def\0{\underline 0}

\DeclareMathOperator{\sing}{Sing}

\makeindex
                              
\begin{document}

\title {On Milnor's fibration theorem and its offspring after  50 years}
%\title {On Milnor's fibration theorem}

\author{Jos\'e Seade}
\address{Instituto de Matem\'aticas, Universidad Nacional Aut\'onoma de M\'exico.}
\email{jseade@im.unam.mx}
\thanks{Research partially supported by  CONACYT and PAPIIT-UNAM from M\'exico, and by the CNRS-UMI 2001, Laboratoire Solomon Lefschetz-LaSoL, Cuernavaca, Mexico.}
\dedicatory{To Jack, whose profoundness  and clarity of vision   \\  seep   into our appreciation  of the beauty and depth of mathematics.}

%\endtitle

%\date{}
%%\endauthor
%%\vskip.3cm
%\date {\it \qquad \quad  \qquad \  \qquad  \qquad \qquad  \qquad \qquad \qquad \qquad \qquad  Dedicated to Jack Milnor, \\ \qquad \qquad  \qquad \qquad \qquad  \qquad  \qquad \qquad  \qquad \qquad \qquad \qquad  a great mind and a great heart.}
%\affil
%\linebreak E-mail:  \endaffil

%\date{10 November 2004}

\def\Z{{\Bbb Z}}
\def\R{{\Bbb R}}
\def\S{{\Bbb S}}

\def\D{{\Bbb D}}
\def\C{{\Bbb C}}
\def\Q{{\Bbb Q}}

\def\B{{\mathbb B}}
 \def\RP{\mathbb{R P}}
\def\F{\mathcal F}
\def\a{\alpha}
\def\0{\underline 0}
\def\b{\beta }
\def\g{\gamma}
\def\ro{\varrho}
\def\na{\nabla}
\def\o{\omega}
\def\e{\varepsilon}
\def\d{\delta}
\def\SU{{\rm SU}(2)}
\def\O{\mathcal O}
\def\K{\mathcal K}
\def\t{\theta}
\def\vt{\vartheta}
\def\nn {\vskip 0.3cm \noindent }
\def\n {\noindent}
\def\nn {\vskip.2cm \noindent}
\def\l{\ell}
\def\la{\lambda}
\def\v {\vskip.1cm}
\def\vv {\vskip.2cm}
\def\grad{\overline\nabla_X f}
\def\s{\mathbb{S}}

\setcounter{section}{0}

%%\NoRunningHeads

\pagestyle{headings}

%\pagenumbering{roman}

\begin{abstract}

Milnor's fibration theorem is about the geometry and topology of real and complex analytic maps near their critical points, a ubiquitous theme in mathematics. As such, after 50 years, this has become a whole area of research on its own, with a vast literature, plenty of different viewpoints, a large progeny and connections with many other branches of mathematics. In this work we revisit the classical theory in both the real and complex settings, and we glance at some areas of current research and connections with other important topics. The purpose of this article is two-fold. On the one hand, it should serve as an introduction to the topic for non-experts, and on the other hand, it gives a wide perspective of some of the work on the subject that has been and is being done. It includes a vast literature for further reading.

\end{abstract}

\maketitle
\bibliographystyle{plain}

\section*{\bf Introduction}

Milnor's fibration theorem in \cite {Mi1} is a milestone in singularity theory that has opened the way to a myriad of insights and new understandings. This is a beautiful piece of mathematics, where many different branches, aspects and ideas, come together.
The theorem concerns the geometry and topology of analytic maps near their critical points.

Consider the simplest case,  a holomorphic map $(\C^{n+1},\0) \buildrel{f}\over {\to} (\C,0)$  taking the origin  into the origin, with an isolated critical point at $\0$. As an example one can have in mind the Pham-Brieskorn polynomials:
\begin{equationth}\label{Brieskorn singularity}
\qquad \qquad \qquad \qquad z \mapsto z_0^{a_0} + \cdots + z_n^{a_n} \quad , \; \quad a_i \ge 2 \; \,{\hbox {for all}} \;  \,i= 0,1,\cdots, n\;.
\end{equationth} \hskip-4pt
Since $f$ is analytic, there exists $r>0$ sufficiently small so that $0 \in \C$ is the only critical value of the restriction $f|_{\B_{r}}$, where $\B_r$  is the open ball of radius $r$ and center at $\0$. Set $$V:= f^{-1}(0) \;\; \,\; \hbox{and } \;\; V^* := V \setminus \{\0\}  \;.$$ 
So $V^*$ is an $n$-dimensional complex manifold. We know (see Section \ref{sec. conical}) that $V^*$ 
meets transversally every sufficiently small sphere $\s_\e$  in $\C^{n+1}$  centered at $\0$ and contained in $\B_r$. The manifold 
$L_V:= V \cap \s_\e$ is called the link of the singularity and its diffeomorphism type does not depend on the choice of the sphere.
Then Milnor's theorem in \cite{Mi1}  says that for every such sphere $\s_\e$ we have a smooth fiber bundle 
\begin{equationth}\label{fib-intro}
\varphi:= \frac{f}{|f|} \,: \, \s_\e \setminus L_V  \longrightarrow \s^1 \;.
\end{equationth} \hskip-4pt
In fact $f$ can have non-isolated critical points. 
The fibers $F_f$ are diffeomorphic to the complex manifolds obtained by considering a regular value $t $ sufficiently near $0 \in \C$ and looking at  the piece of $f^{-1}(t)$ contained within the open ball $\B_\e$ bounded by $\s_\e$. 

\begin{figure}
\centering
\includegraphics[height=6cm ]{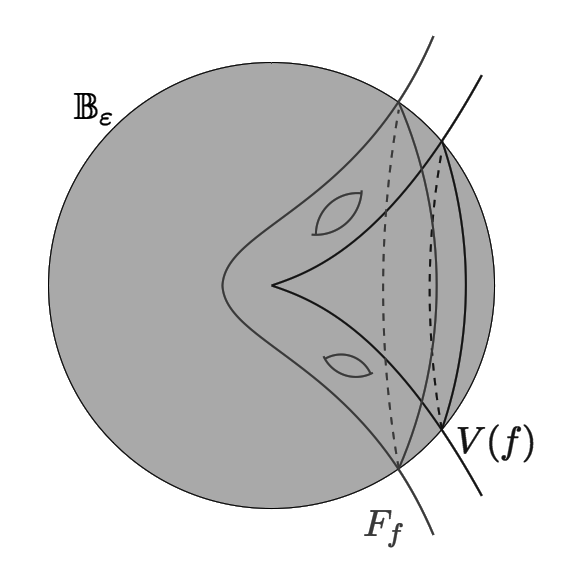}
\caption{The Milnor fiber $F_f$.}
\end{figure}

Early versions of this theorem arose when Brieskorn \cite{Brieskorn} discovered that  the fibration (\ref{fib-intro}) can be used to prove that 
exotic spheres show up naturally as links of isolated singularities of complex hypersurfaces. In fact Brieskorn proved that every exotic sphere that bounds a parallelizable manifold is diffeomorphic to the link of a singularity as in (\ref {Brieskorn singularity}), for appropriate $a_i$.
\

These kind of ideas led to studying  the local conical structure of all real and complex analytic sets, that we discuss in Section \ref {sec. conical}. 
 In the sequel we describe the fibration theorem for real and complex singularities, and various extensions of it. We discuss too the main ideas of its proof and we exemplify these with the case of the Pham-Brieskorn polynomials, where the main ingredients of the proof are shown in an elementary way.

 When the map in question is holomorphic and has an isolated critical point, the fiber is diffeomorphic to a $2n$-ball to which one attaches handles of middle index. The number of such handles is now called the Milnor number of the singularity.  If the critical point is non-isolated, then the Milnor fiber is diffeomorphic to a ball to which we must attach handles of various indices. The precise number of handles of each index is prescribed by the L\^e numbers of the singularity, a concept introduced by D. Massey in \cite{Massey0, Massey1}.

The study of the Milnor number and its generalizations, such as the L\^e numbers, has given rise to a vast literature, discussed in sections \ref{sec. topology} to  \ref {sec. Milnor classes}. In fact there are interesting relations with the theory of Chern classes for singular varieties: one has the Milnor classes, defined as the difference between the Schwartz-MacPherson and the Fulton-Johnson classes. For varieties which are complete intersections with isolated singularities, there is only one Milnor class, an integer  equal to the sum of the local Milnor numbers at the singular points. So Milnor classes generalize the Milnor number to varieties which may not be complete intersections and 
can have non-isolated singularities.

In section \ref{sec. equisingular} we glance at the important concept of 
 equisingularity. This   emerges from Zariski's seminal work  \cite {Zar1, Zar2, Zar3,Zar4,Zar5} and it has important ties with Milnor's fibration theorem.  The basic idea is that, if $V$ is a variety with non-isolated singularities,  we want to know whether  the singularity of $V$ at a given point  is ``worse than" or ``equivalent to"  its singularities at other nearby points; and what does ``equivalent" mean?
 
As mentioned above, singularity theory is a meeting point of various areas of mathematics. In Section \ref {sec. relations} we briefly discuss relations of complex singularities with open-book decompositions and fibered knots, contact structures and low dimensional manifolds. Later  we discuss relations with the so-called moment-angle manifolds, which spring from mathematical physics. For this we need to extend the discussion beyond the holomorphic realm.

Milnor proved in \cite {Milnor:ISH, Mi1} that the fibration theorem also holds for 
  analytic maps $(\R^{n},\0) \buildrel{f} \over {\to} (\R^p,0)$, $n \ge p \ge 1$, with  an isolated critical point, and  one has a fiber bundle:
\begin{equationth}\label{real}
\varphi \,: \, \s_\e \setminus L_V  \longrightarrow \s^1 \;.
\end{equationth} \hskip-4pt
Yet, in general the projection map $\varphi$ can only be taken as ${f}/{|f|}$ in a neighborhood of the link $L_V$. Also, as pointed out in \cite[Chapter11]{Mi1}, 
the  condition of having  an isolated critical point is very stringent. Generically the set of critical values has positive dimension, and even if the only critical value is $0$, it is fairly stringent to ask for having an isolated critical point. In sections \ref{foundations} to \ref{d-regularity}   
we discuss  Milnor fibrations for real analytic maps in general, starting with the isolated singularity case. We discuss a regularity condition that  is necessary and sufficient to assure that if we have the fibration (\ref{real}), we can take the projection $\varphi$  to be ${f}/{|f|}$ everywhere. This is called $d$-regularity.

We look also at  two  particularly interesting classes of  maps: the meromorphic germs and maps of the form $f \bar g$ with $f,g$ holomorphic. Notice that away from the set  $f g = 0$  one has $$\frac{f/g}{|f/g|} = \frac{f \bar g}{|f \bar g|}\,.$$
The study of Milnor fibrations for meromorphic germs began with a couple of papers \cite {GLM1, GLM2}  by Gusein-Zade {\it et al} (see also \cite{BP, BPS}). The study of Milnor fibrations and open-books defined by  functions  $f \bar g$ essentially began in \cite {Pi2, PS2}, though this type of maps has already appeared  in \cite{AC1, HR, Mi1}. 
We study these in a section  under the name  ``Mixed singularities", which by definition are critical points of analytic functions in the complex variables $z_1,\cdots, z_n$ and their conjugates. 

Polar weighted homogeneous singularities are a particularly interesting class of mixed functions which are reminiscent of the classical weighted homogeneous polynomials.  Now one has a weighted action of $ \s^1 \times \R^+ \cong \C^* $ on $\C^n$, and the function brings out scalars to some power. Unlike the classical case of weighted homogeneous functions,  now  $\s^1$ and $\R^+$ can act with  different weights. The paradigm examples of such singularities are the twisted Pham-Brieskorn polynomials, studied in \cite{RSV, Se5, Seade-librosing}:
$$
\qquad (z_1,...,z_n) \;   {\longrightarrow} \;  z_1^{a_1} \bar z_{\sigma(1)} + ... +
z_n^{a_n} \bar z_{\sigma(n)} \quad , \quad a_i \ge 2\,,
$$
where $\sigma$ is a permutation
of the  set $\{1,...,n\}$. If $\sigma$ is the identity, these are essentially equivalent to usual Pham-Brieskorn singularities (by \cite{RSV, Oka3}).

The notions of {\it weighted polar actions} and  {\it  polar singularities}  were introduced in \cite{Cis}, inspired by \cite{RSV}.
The  name {\it mixed singularity} was coined  by M. Oka \cite{Oka0} and his vast contributions to the subject have turned this into a whole new area of research.  These singularities also provide interesting open-book decompositions similar to Milnor's open-books but which do not appear in complex singularities (see for instance \cite{Hayde1, Hayde2, Hayde3, PS1}).

%%%%%%%%%%%%%%%%%%%

%%%%%%%%%%%%%%%%%%%

We finish this article with a specially interesting class of mixed singularities that spring  from  the study of holomorphic linear actions of $\C^m$ in $\C^n$, $0 < m \ll n$. 
Under appropriate  conditions, there is an open dense set in $\C^n$ of points that belong to  a type of orbits called ``Siegel leaves". These are parameterized by a $C^\infty$-manifold $V^* := V \setminus \{\0\}$, where $V$  
is a  complete intersection  defined by $2m$  real-valued quadratic polynomial equations; $V^*$ consists of the points where the orbits are tangent to the foliation in $\C^n \setminus \{0\}$ given by all spheres centered at the origin. The  manifold $V^*$  has a canonical complex structure that comes from being everywhere transversal to the orbits of the $\C^m$-action, and it has a canonical $\C^*$-action with compact quotient $V^* /\C^*$. These quotients are known as  {\it LVM-manifolds}; these have a  rich geometry and topology, and  they are a class of the so-called {\it moment-angle manifolds}, of relevance in mathematical-physics.

We begin this article by pointing out important relations between complex singularities and exotic spheres, and we close it with an  important relation between real singularities and moment-angle manifolds: this exemplifies  the  strong interactions singularity theory has with  other areas of mathematics. 

My aim with 
this article  is to  share with the readers some of the beauty and richness of singularity theory, a fascinating area of mathematics  which somehow begins with Isaac Newton and is built upon the foundational work of Hassler Whitney,  Ren\'e Thom,  John Milnor, Oscar Zariski, Heisuke Hironaka, Egbert Brieskorn, Vladimir Arnold and many others. Arnold used to say \cite {Arnold2} that the main goal in most problems of singularity theory is to understand the qualitative change in objects that depend on parameters, and which may come from analysis, geometry, physics or from some other science. Milnor's fibration theorem certainly fits within this framework.

I am happy to thank  Patrick Popescu-Pampu, Bernard Teissier, L\^e D\~ung Tr\'ang, Javier Fern\'andes de Bobadilla, David Massey, Jawad Snoussi, Jos\'e Luis Cisneros, Roberto Callejas-Bedregal, Michelle Morgado and Aurelio Menegon, for  wonderful conversations and useful comments that greatly enriched this presentation. I am also indebted to Martin Guest, for reading the manuscript and make me many helpful remarks, and to the referee, whose comments made this article far more readable.

\tableofcontents

%%%%%%%%%%
%%%%%%%%%%

%\newpage

\section*{}

\specialsection*{\bf PART I: COMPLEX SINGULARITIES}

\section{\bf The initial thrust: searching for homotopy spheres}\label{sec. initial}

In 1956 John Willard Milnor (b. 1931) surprised the world by finding \cite{Mi0} the first ``exotic spheres'':  7-dimensional smooth manifolds homeomorphic to ${\S}^7$ but with non-equivalent differentiable structures.

The  set of equivalence classes of smooth structures on the $n$-sphere $\s^n$ is a monoid where the operation is the connected sum. For $n \ne 4$ this monoid is   a group and it is isomorphic to the finite abelian group $\Theta_n$ of h-cobordism classes of oriented homotopy $n$-spheres,  with the connected sum as operation;
 the identity element is the standard $\s^n$.  This group was studied by Kervaire and Milnor in \cite{Ker-Mil}  for $n \ge 5$. They noticed that $\Theta_n$  contains a ``preferred subgroup'',  denoted 
$bP_{n+1} \subset \Theta_n$, of those homotopy spheres that bound a parallelizable manifold, {\it i.e.}, a manifold with trivial tangent bundle.  
For $n \ne 3$ odd this is a finite cyclic group which has finite index in $\Theta_n$. This cyclic group has order $1$ or $2$ for 
$n \equiv 1 \, (\hbox{mod} \, 4)$, but for $n \equiv 3 \, (\hbox{mod} \, 4)$ its order $ |bP_{4m}|$ grows more than exponentially:
$$ |bP_{4m}| \,=\, \big[2^{2m-2}(2^{2m-1}-1)\big] \, \cdot \, \big[ \hbox{numerator of} \,(\frac{4 \mathcal B_m}{m}) \big]\,\,,$$
where the $ \mathcal B_m$ are the Bernoulli numbers.  Thus for instance (see~\cite {Ker-Mil, Hir}), for $n = 7, 11, 15$ or $ 19\,$ there are, 
respectively, $ |bP_{n+1}| \,=\, 28,\, 992,\, 8128\,$ and $\,130816 \,$ non-equivalent differentiable structures on the
$n$-sphere that bound a parallelizable manifold. 

A question was how to construct those exotic spheres? An approach to this problem is by looking at the links of isolated complex singularities. This  led to Milnor's fibration theorem as we now explain. 

Let $(V,P)$ be a complex analytic variety of complex dimension $n$  in some affine space $\C^N$,  with a unique singular point at $P$.
Then
$ V^*= V \setminus \{P\}$ is a complex $n$-manifold.   One has:

\begin{proposition}\label{Milnor-cone}
There exists $\e>0$ sufficiently small, so that every sphere $\S_r$ in $\C^N$ of radius $r \le \e$ and center at $P$ meets $V^*$ transversally.
\end{proposition}

This  is proved in \cite[Chapter 2]{Mi1}  when $X$ is algebraic and it is a particular case of a  general theorem about the local conical structure of analytic sets, see Section \ref{sec. conical}.
It follows that if $(V,P)$ is as above, then {\it its link} $L_V := \S_\e \cap V$ is a smooth real analytic  manifold of dimension $2n-1$.

\begin{question}\label{question-homotopy spheres}
Can we know when $L_V $ is a homotopy sphere?, and if so, can we determine which element in $\Theta_n$ it represents?
\end{question}

For $n=1$ the question is trivial since $L_V $ is a union of circles, one for each branch of $V $. For $n=2$, if $V$ has a normal  singularity at $P$ then its link   is never simply connected, by \cite{Mum}.
When $V$ is a complex hypersurface, {\it i.e.}, defined by one single equation, Question \ref{question-homotopy spheres} was  answered  by the work of various people in the 1960s, most notably by E. Brieskorn, F. Hirzebruch and J. Milnor, see  \cite{Brieskorn, Hir2, Mi1}.

For instance, in dimension 7 we know from \cite{Ker-Mil} that there exist  28 classes of homotopy spheres including the standard one, and all of them bound a parallelizable manifold. Brieskorn proved in \cite{Brieskorn} that all these 28 classes can be represented  by the link $L_V $ of some hypersurface in $\C^5$ of the form 
$$ z_0^{a_0} + z_1^{a_1} + \cdots +  z_{4}^{a_4}  \, =\, 0 \,. 
$$
The fibration theorem is the culmination of a series of works by various authors, starting with F. Pham \cite{Pham}. 

Let $V$ be the zero-locus of an analytic map $(\C^{n+1},\0) \buildrel{f} \over {\to} (\C,0)$ with an isolated critical point at $\0$. Equip its link  
$L_V ^{2n-1}$   with its  natural differentiable structure as the transverse  intersection $L_V  = \S_\e \cap V$ of two smooth submanifolds of $\C^{n+1}$. One has a map:
\begin{equationth}\label{fibration-intro}
\varphi:= \frac{f}{|f|} : \S_\e \setminus L_V  \longrightarrow \S^1 \,,
\end{equationth} \hskip-4pt
and Milnor's fibration theorem says that this is a smooth  fiber bundle.

Milnor  also proves that the fiber
$F_t$ is diffeomorphic to the portion of a non-critical level $ f^{-1}(t)$ contained within the ball $\B_\e$ bounded by $\S_\e$.  
This implies that the normal bundle of $F_t$ is  trivial, being the inverse image of a regular value. Hence the tangent bundle $TF_t$ is stably trivial, {\it i.e.}, $F_t$ is stably parallelizable, and we know from \cite{Ker-Mil} that for compact connected manifolds with non-empty boundary, stably-parallelizable implies parallelizable. Thus we get:

\begin{proposition}
The link $L_V$ of every complex hypersurface isolated singularity bounds the fibers $F_\t$, which are parallelizable manifolds. 
\end{proposition}

The point is to know when $L_V $ is  a homotopy sphere, and when this happens, which element it represents in  $bP_{2n}$.  For this, Milnor proved in \cite {Mi1}:

\begin{proposition}
 The link of every isolated hypersurface singularity  in $\C^{n+1}$ is $(n-2)$-connected and  the fiber $F_t$ has the homotopy of a bouquet ${\bigvee} \S^n$ of spheres of middle dimension.
\end{proposition}

The number of spheres in the bouquet ${\bigvee} \S^n$ is strictly positive, unless $V$ has no singularity. This is now called the Milnor number $\mu$ of $f$ (see Section \ref{sec. topology}).

For $n > 2$ the link is simply connected and therefore the Hurewicz isomorphism implies that 
   the homology of $L_V $ also vanishes in dimensions $ i = 1, \cdots, n-2$. 
  Since $L_V $ is always  orientable,
 by the Poincar\'e duality isomorphism  its homology  vanishes in dimensions $n+i$, $ i = 1, \cdots, n-2$ as well. Thus 
the only possibly 
non-zero groups are in dimensions $i = n, n-1$ and of course $i = 0, 2n-1$ where they are isomorphic to the group of the integers (or the corresponding ring of 
coefficients).

If $H_{n-1}(L_V )$ vanishes then $H_{n}(L_V )$ also vanishes, by duality,  and $L_V $ is a homology sphere. If $n \ge 2$, then $L_V $ is 
simply connected by~\cite {Mi1}. Hence, if $n > 2$ 
then Smale's theorem in cite{Smale}  implies (Poincar\'e's Conjecture in dimensions $\ge 5$) that $L_V $ is actually homeomorphic
to $\s^{2n-1}$.  So the question is to decide when $H_{n-1}(L_V )$ vanishes, and it is here that the fibration theorem enters the scene:  
fix a fiber $F_0$ and notice that  $F_0$ and its complement in $\S_\e \setminus L_V $ have the same homotopy type. Consider the monodromy $h$ (a first return map) of the bundle (\ref{fibration-intro}) and  
 the induced representation in the middle homology of the fiber:
  $$h_* : H_n(F_0) \to H_n(F_0)\,.$$ One has the corresponding 
 Wang sequence  (see p. 68 in~\cite {Mi1}):
$$H_n(F_0) \buildrel
{h_* - I_*}\over\longrightarrow  H_n (F_0) \longrightarrow H_n(\s_\e - L_V ) \longrightarrow 0 \,\,.$$
Using this one arrives at Milnor's theorem 
  \cite[Theorem 8.5]{Mi1}, that for $n > 2$ the link  is a topological sphere if and only if the determinant of ${h_* - I_*}$ is $\pm 1$,  which was already  proved  by F. Pham  in \cite{Pham} for the Pham-Brieskorn polynomials.
  Using this, for instance Hirzebruch  proved in \cite [p. 20-21]{Hir2} that  the links of the following singularities are all homotopy spheres:
$$\qquad z_0^3 +  z_1^{6k-1} + z_2^2 +\cdots +  z_{2m}^2 = 0 \quad ; \;   k \ge1 \;, \; m \ge 2\;.$$
One has the following remarkable theorm of Brieskorn \cite[Korollar 2]{Brieskorn}:

\begin{theorem}
Every exotic sphere of dimension $m =2n-1 > 6$  that bounds a parallelizable manifold 
is the link 
of some hypersurface singularity of the form $$ z_0^{a_0} + z_1^{a_1} + \cdots +  z_{n}^{a_n}  \, =\, 0 \,, $$
for some appropriate integers $a_i \ge 2$, $i = 0, 1, \cdots, n$.
\end{theorem}

One may consider
singularities which are not hypersurfaces, and try to produce other elements in the homotopy of spheres. To my knowledge, little is known about this problem. 
If we consider   
complex isolated complete intersection singularities, one always has a Milnor fibration and the fibers are compact parallelizable manifolds with boundary the link, by \cite{Hamm}. So in these cases, if the link if an exotic sphere,  this is in $bP_{2n} \subset \Theta_{2n-1}$, which is the simplest and best understood part of $\Theta_{2n-1}$.

\vv

 I thank Patrick Popescu-Pampu for bringing to my attention the following interesting question posed by A. Durfee in \cite[Problem H, p.  252]{Web}:
 
\begin{question}\label{q. Durfee}
Does every exotic sphere occur as the link of an isolated complex singularity? 
\end{question}

A step for answering (\ref {q. Durfee}) is the question that Popescu-Pampu originally asked me:

\begin{question}
Does there exist a complex isolated singularity whose  link is a homotopy sphere that does not bound a parallelizable manifold? 
\end{question}

Such examples, if they exist, would produce elements in the most mysterious part of the groups $\Theta_n$.

\section{\bf An example: the Pham-Brieskorn singularities.}\label{sec. Brieskorn}

 Consider a Pham-Brieskorn polynomial $f: \C^n \to \C$:
\begin{equation*}%\label{d:twisted Pham-Brieskorn}
(z_1,...,z_n) \buildrel{f} \over {\mapsto} z_1^{a_1}  + ... +
z_n^{a_n},  \quad  \; a_i \ge 2\,.
\end{equation*}
The origin $\0 \in \C^n$ is its only critical point, so $V :=
f^{-1}(0)$ is a complex hypersurface with an isolated singularity at $\0$.
Let
$d$ be the lowest common multiple of the $a_i$ and for each
$i=1,\cdots, n$ set $d_i = d/a_i$. Then for every non-zero complex
number $\lambda \in \C^*$ one has  a $\C^*$-action on $\C^n$
determined by $$\lambda \cdot (z_1,\cdots,z_n) \mapsto
(\lambda^{d_1} z_1,\cdots,\lambda^{d_n} z_n)\,.$$ Observe that one has:
$$f(\lambda^{d_1} z_1,\cdots, \lambda^{d_n} z_n) = \lambda^d f(z_1, \cdots,z_n) \,,$$ so $f$  is {\it weighted homogeneous}.
This $\C^*$-action  has 
$\0$ as its only fixed point  and $V$ is an invariant set, union of
$\C^*$-orbits.  This has the following
 properties:

\medskip
$\bullet$  {\bf Property 1.} Restricting the action to the
positive real numbers $t \in \R^+$, we get a flow such that:

     - each orbit  converges to $\0$ as $t$ tends to $0$, and  it goes to
infinity as $t$ tends to $\infty$;

     - each  orbit is transversal to all spheres centered at $\0$. 
Hence  $V$ intersects transversally every $(2n - 1)$-sphere
$\s_r$ centered at $\0$, so   $K_r := V \cap
\s_r$ is a real codimension 2 smooth submanifold  of 
$\s_r$;

 -  Given arbitrary spheres $\s_r, \s_{r'}$
centered at $\0$, the flow gives a diffeomorphism from $\s_r$ into $
\s_{r'}$ taking $K_r$ into $K_{r'}$. Moreover, the flow  determines a 1-parameter group of
diffeomorphisms that exhibits the pair $(\C^n,V)$ as the cone over  $(\s_r, K_r)$.
We denote the manifold  $K_r$ by $L_f$ and call it {\it the link} (see \cite{Mi1, Du3}).

\medskip

$\bullet$ {\bf Property 2.} The argument
of the complex number $f(z)$  is  constant  on each orbit of the above flow, {\it
i.e.}, $f(z)/|f(z)| = f(tz)/|f(tz)|$ for all $ t \in \R^+$. 

\medskip

$\bullet$ {\bf Property 3.}  The restriction of the $\C^*$-action to $\s^1$ 
leaves invariant every sphere around $\0$.  Multiplication by $e^{i\t}$
in $\C^n$  transports each fiber $f^{-1}(\zeta)$  into the fiber
$f^{-1}(e^{i \t d} \cdot  \zeta)$.  Hence $\s^1$ acts on each tube $N(\delta):= f^{-1}(\partial \D_\delta) \,,$
where $\partial \D_\delta \cong
\s^1$ is the boundary of the disc in $\C$ of radius $\delta > 0$ and  
centered at $0$. A direct
computation shows that the orbits of this action are transverse to
the fibers of $f$. So we have a smooth  fiber bundle:
$
 N(\delta)  \buildrel {f} \over {\to} \partial \D_\delta \,.
$

\begin{figure}
\centering
\includegraphics[height=6cm ]{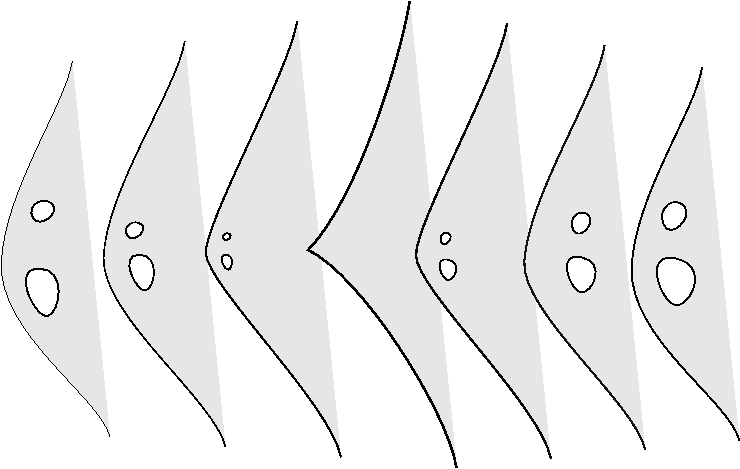}
\caption{The non-critical levels are all diffeomorphic.}
\end{figure}

\medskip
Now observe that
for each   line  $\cL_\t$ through the origin in $\C$,  we may consider the set
\begin{equation*}
 X_\t := \{z \in \C^n \, | \,  f(z) \in \cL_\t \} \,.
\end{equation*}
Each $X_\t $    is a real analytic
hypersurface with an isolated singularity at $\0$, their union
fills the entire  $\C^n$ and their intersection is $V$. By
Property 1, each $X_\t $ is transversal to all the
spheres, and by Property 3, the $\s^1$-action permutes these
hypersurfaces. Thus one has:

\medskip
$\bullet$ {\bf Property 4.} These varieties define a pencil in
$\C^n$, a sort of open-book where the binding is now the singular
variety $V$, and each of these varieties  is transverse to every
sphere around $\0$. If we remove $V$ from $\C^n$, for every ball
around $\0$ we get a fiber bundle
\begin{equationth}\label{global fibration}
\varphi = \frac{f}{|f|}: \B^{2n} \setminus V  \longrightarrow
\s^1\,.\end{equationth} 
The fiber
over a point $e^{i \t}$ is  a  connected component
of $X_\t \setminus V$. The other component is $f^{-1}(e^{-i \t})$. 
\medskip

\medskip
$\bullet$  {\bf Property 5.} We now focus our attention near the origin, say restricted to
the unit ball $\B^{2n}$  in $\C^{n}$. Since each $X_\t$ meets
transversally the  sphere $\s^{2n-1} = \partial \B^{2n}$,  the
intersection is a smooth codimension 1 submanifold of the sphere,
containing the link $L_f = V \cap \s^{2n-1}$. And since the orbits
of the $\s^1$-action preserve the sphere  $\s^{2n-1} $,  the restriction of $\varphi$ to $\s^{2n-1}$ defines
the classical {\it Milnor fibration}:
\begin{equationth}\label{fibration-on-spheres}
\varphi = \frac{f}{|f|}:  \s^{2n-1} \setminus L_f  \longrightarrow
\s^1\,\end{equationth} \hskip-4pt
 
\medskip
$\bullet$  {\bf Property 6.}
Since $V$ is transverse to
the unit sphere $\s^{2n-1} $ and each point in $V \setminus \{\0\}$ is regular, each fiber $f^{-1}(t)$ with $|t|$ sufficiently small 
 is also transverse
to $\s^{2n-1} $. Hence, if  $N(\delta)$ is as in Property 3 and we set $N(1,\d) := N(\delta) \cap \B^{2n}$, where the
$1$ means that the ball $\B^{2n}$ has radius 1, we have that the fiber
bundle described by Property 3 determines a fiber bundle:
\begin{equationth}\label{fibration-on-tubes}
f\colon N(1,\d) \to \partial \D_\delta \cong \s^1 \,.
\end{equationth}
This is the second classical version of Milnor's fibration for the
map $f$. 

\medskip
$\bullet$  {\bf Property 7.}
Notice that by Property 2, each $\R^+$-orbit is everywhere
transverse to  the tube $N(\d)$ and  transverse to the sphere
$\s^{2n-1} $,  and the complex numbers $f(z)$ have constant argument
along each orbit. Thence the integral lines of this action
determine a diffeomorphism between $N(1,d)$ and $\s^{2n-1} $ minus
the part of the sphere contained inside the open solid tube
$f^{-1}(\buildrel {\circ} \over {\D_\d})$. This determines the
classical equivalence between the  Milnor fibration
 in the sphere (\ref{Fibration Thm., version1}) and the Milnor-L\^e
fibration  in the tube (\ref{Fibration Thm., version2}).

%%%%%%%%%%%%%%%%
\medskip
We now remark that everything we said above works in exactly the
same way for all weighted homogeneous complex singularities, {\it
i.e.},  for all complex polynomials $f$ for which there is a
$\C^*$ action on $\C^n$  as above, for some positive integers
$\{d; d_1,\cdots,d_n\}$,
$$\lambda \cdot (z_1,\cdots,z_n) =
(\lambda^{d_1} z_1,\cdots,\lambda^{d_n} z_n) \,,$$ satisfying
that for all $\lambda \in \C^*$ and for all $z \in \C^n$ one has:
$$f(\lambda \cdot z) = \lambda^d \cdot f(z) \,.$$
These all have the same Properties 1 to 7. 
As we will see in the sequel, all real analytic isolated singularities can be equipped with flows that satisfy
properties analogous to  1, and 3 to 6, but not always 2. This implies that we have a fibration as in (\ref
{fibration-on-tubes}) and it can be carried to a fibration on the sphere as in (\ref{fibration-on-spheres}) but the projection map $\varphi$ may not always be taken to be $f/|f|$ away from a neighborhood of the link. Having also property  2 grants that $\varphi$ can be taken as $f/|f|$ everywhere, and this
is equivalent to the map-germ being $d$-regular, a concept that we discuss  in Section \ref{d-regularity}.

%%%%%%
%%%%%%
%%%%%%
%%%%%%

\section{\bf Local conical structure of analytic sets}\label{sec. conical}
Consider a reduced, equidimensional 
real analytic space $V$ of dimension $n$, defined 
in an open ball $\B_r(\0) \subset \R^N$ around the origin. Assume  $V$ contains the 
origin $\0$ and $V^* := V -\{\0\}$  is a real analytic manifold of dimension 
$n > 0$.  The following  is proved in \cite{Mi1} and it can be deduced from   \cite{Loj}:

\begin{theorem}[Milnor 1968]\label{thm. cone-isolated case}
There exists $\e>0$ sufficiently small, so that every sphere in $\R^N$ centered at $\0$ and with radius $\le \e$ intersects $V^*$ transversally. Moreover, there is a smooth 1-parameter family of diffeomorphisms $\{\gamma_t\}$, $t \in [0,\e)$, such that $\gamma_0$  is the identity and if $\S_{\e -t}$ denotes the sphere of radius ${\e -t}$, then each $\gamma_t$ carries the pair $(\S_{\e}, \S_{\e} \cap V)$ into $(\S_{\e -t}, \S_{\e -t} \cap V)$. The pair $(\B_\e, \B_\e \cap V)$ is homeomorphic to the cone over  $(\S_{\e}, \S_{\e} \cap V)$.
\end{theorem}

\begin{figure}
\centering
\includegraphics[height=6cm ]{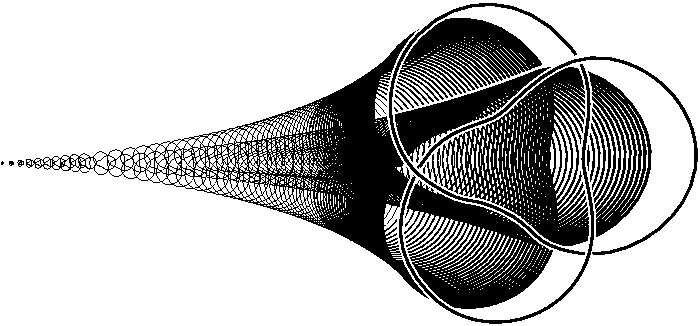}
\caption{The link of the singularity determines the topological type.}
\end{figure}

The idea of the proof is simple: consider the function $\R^N \buildrel{d} \over {\to} \R$ given by $d(x_1,...,x_N) = x_1^2 + \cdots x_N^2$, so that $d$ is
 the square of the function ``distance to $\0$''. The solutions of its gradient vector field $\nabla d$ are the straight rays that emanate from the origin. 
Let us adapt this vector field to $V$. 
For this, take the restriction $d_V$ of $d$ to $V$.
At each point $x \in V^*$ the gradient vector $\nabla d_V(x)$ is obtained by projecting $\nabla d(x)$
to $T_xV^*$, the tangent space of $V^*$ at $x$, so 
 $\nabla d_V(x)$ vanishes if and only if  $T_xV^* \subset T\s_x$. This means  that  a point $x \in V^*$ is a critical point of $d_V$
iff $V^*$ is tangent at $x$ to the sphere passing through $x$ and centered at $\0$.
Just as in~\cite[Corollary 2.8]{Mi1}, one has that $d_V$ has at most a finite number of critical values
corresponding to points in $V^*$, since it is the restriction of an analytic function on $\B_r(\0) $.
 Hence  $V^*$ meets transversally 
all sufficiently small spheres around the origin in $\R^N$. The gradient vector field of $d_V$ is
now everywhere transversal to the spheres around $\0$, and it can be assumed to be integrable. Hence it defines 
a 1-parameter family of local diffeomorphisms of $V^*$ taking each link into  ``smaller'' links, proving Theorem \ref{thm. cone-isolated case}.

\vv

Theorem \ref{thm. cone-isolated case} was extended in \cite{Bu-Ve}  to varieties with arbitrary singular locus  using Whitney stratifications (we refer to \cite{Gor-Mc, Dimca} for background material on stratifications). A more refined argument due to A. Durfee~\cite {Du3} (see also \cite {Le-Te2}) and based on the ``Curve Selection Lemma'' of 
\cite {Mi1}, shows that in fact the diffeomorphism type of the manifold $V \cap \s_\e$ is also independent of the 
choice of the embedding of $V$ in $\R^N$.  
One has:

\begin{theorem} \label{Conical-structure}
Let $V$ be a real or complex analytic set in $\R^m$ and $P$ a singular point in $V$. Then there exists a Whitney stratification of $\R^m$ for which $V$ is a union of strata, $P$ is a point stratum, and one has:
\begin{enumerate}
\item There exists $\e>0$ sufficiently small, so that every sphere $\S_r$ in $\R^m$ of radius $r \le \e$ and center at $P$ meets transversally every stratum in $V$.
\item One has a homeomorphism of pairs: $(\B_\e, \B_\e \cap V) \cong  {\rm Cone} \, (\S_\e, \S_\e \cap V)$.
\item The homeomorphism type of $L_V := \S_\e \cap V$ is independent of the choice of the defining equations for $V$. 
\end{enumerate}
\end{theorem}

%%\midinsert
%%\vspace{6truecm}
%%\botcaption{Figure 3}
%%{\it The link of the singularity determines the topological type }
%%\endcaption
%%\endinsert

\nn
\begin{definition}\label{Definition 2.1}
The manifold $L_V := V \cap \s_\e$ is called {\it the link}
\index{link, of the singularity} of $V$ at $0$, and a sphere $\S_\e$ as in Theorem \ref{Conical-structure} is called {\it a Milnor sphere} for $V$. We also denote $L_V$ by $L_f$ when we want to emphasize the function rather than the space $V$.
\end{definition}

%\newpage

%%%%%%%%%%%%%%
%%%%%%%%%%%%%%%
\section{\bf The classical fibration theorems for complex singularities}\label{sec. classical fibration theorems}

The  first version of Milnor's fibration theorem says:

\begin{theorem}[Fibration Theorem, 1st version]\label{Fibration Thm., version1}
Let  $U$ be an open neighborhood of the origin $\0 \in \C^{n+1}$ and
$f:(U, \0) \to (\C,0)$ a complex analytic map. Set $V:= f^{-1}(0)$ and $L_V:= V \cap \s_\e$ where $\s_\e$ is a sufficiently  small sphere in $U$ centered at $\0$. Then,
$$ \varphi:= \frac{f}{|f|} : \s_\e \setminus L_V \longrightarrow \s^1 \;,$$
is a $C^\infty$ fiber bundle.
\end{theorem}

The proof by Milnor in \cite{Mi1} uses the Curve Selection Lemma, first to show that the map $ \varphi$ has no critical points  and then 
to construct an appropriate vector field on $\s_\e \setminus L_V$ that shows the local triviality, {\it i.e.}, that each fiber of $\varphi$ has a neighborhood which is a product.
Nowadays, the most common proof of Theorem \ref{Fibration Thm., version1} follows the original approach sketched by Milnor himself in a previous unpublished article \cite{Milnor:ISH}. The starting point is the following:

\begin{theorem}[Fibration Theorem, 2nd version]\label{Fibration Thm., version2}
With the above hypothesis and notation,   let $\delta >0$ be sufficiently small with respect to $\e$, so that for every $t \in \C$ with $|t| \le \delta$  the fiber $f^{-1}(t)$ meets the sphere $\s_\e$ transversally. Let $\D_\delta$ be the disc in $\C$ of radius $\delta$ and center at $0$; let $\partial \D_\delta \cong \s^1$ be its boundary  and set  $N(e,\delta) := f^{-1}(\partial \D_\delta) \cap \B_\e$, where $\B_\e$ is the open ball in $\C^{n+1}$ bounded by $\s_\e$. Then,
$$ f_{|_{N(e,\delta)}} : N(e,\delta) \longrightarrow \partial \D_\delta \;,$$
is a $C^\infty$ fiber bundle, (essentially) isomorphic to that in Theorem \ref {Fibration Thm., version1}.
\end{theorem}

The word essentially in the last statement  is because the fibers in \ref{Fibration Thm., version2} are compact, while those in 
Theorem \ref {Fibration Thm., version1} are open manifolds. To  have an actual isomorphism of the two fibrations one must  restrict the fibration in \ref{Fibration Thm., version2} to the open ball.

 Milnor  proved this theorem in \cite{Milnor:ISH} when the map-germ $f$ has an isolated critical point. In the general case, Milnor proved in \cite {Mi1} that the fibers in Theorem \ref{Fibration Thm., version2} (restricted to the open ball) are diffeomorphic to those in Theorem \ref{Fibration Thm., version1}. In order to prove that one actually has a fiber bundle  in Theorem \ref{Fibration Thm., version2} one must grant that given $\e>0$ as above, there exists a $\delta$ as stated, such that  all fiber $f^{-1}(t)$ with $|t| \le \delta$ meet the sphere $\s_\e$ transversally.
 This was not known till 1977 when Hironaka proved in {\cite{Hironaka:SF} that all complex valued holomorphic maps have a Thom Stratification.

 To prove Theorem \ref{Fibration Thm., version2}  we restrict $f$ to a sufficiently small open ball $\B_r$  around $\0$ so that  $0 \in \C$ is its  only critical value. 
We  equip  $\B_r$ with a Thom stratification \cite{Hironaka:SF}, so that $V$ is a union of strata and we  assume that $\0$ itself is a stratum. 
Now let $\s_\e \subset \B_r$ be a
 Milnor sphere for $f$, so that every sphere of radius $\le r$ meets transversally each stratum in $V$; this is possible by Theorem \ref{Conical-structure}. By compactness, this implies that  there exists $\delta >0$ such that for each $t \in \C$ with $0 < |t| \le \delta $, the fiber $f^{-1}(t)$ meets $\S_\e$ transversally. 
Hence all  fibers in \ref{Fibration Thm., version2} are compact smooth manifolds with boundary. The proof of the local triviality is  like in the usual proof of Ehresmann's fibration lemma, by lifting via the Jacobian of $f$ appropriate vector fields in $\C$ to vector fields in $\B_\e$ which are normal to the fiber. The only additional thing is that we must choose the liftings so that the vector fields are  also tangent to $\S_\e$, which is possible because the fiber is transversal to the sphere (see \cite{Seade-librosing, Wolf}).

The next step to prove  Theorem \ref{Fibration Thm., version1} is implicit in \cite{Mi1}: we inflate the Milnor tube, carrying the fibration in the tube into the fibration in the sphere as stated.
 This relies on the Curve Selection Lemma. The key for this is constructing a  vector field as stated in the following lemma.

\begin{lemma}\label{equivalence-of-fibrations} 
There exists an integrable vector field 
 $\xi$ on $\B_\e \setminus V$ such that:

\begin{figure}
\centering
\includegraphics[height=7cm ]{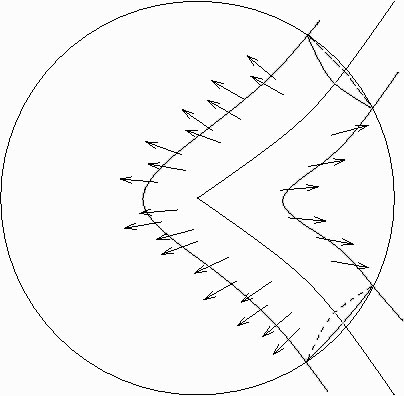}
\caption{The vector field that carries a Milnor tube into the sphere.}
\end{figure}

\begin{enumerate}
\item  Its integral lines are transversal to all Milnor tubes $f^{-1}(\S^1_r)$;
\item Its integral lines are transversal to all spheres centered at $\0$;
\item Its integral lines travel along points where  $f$ has constant argument. That is, if $z, w$ are points in $ \B_\e \setminus V$ which are in the same integral line of $\xi$, then $f(z)/|f(z)| = f(w)/|f(w)|.$
\end{enumerate}
\end{lemma}

We remark that constructing a vector field that satisfies the first two conditions is easy and can be done also in the real analytic category. This allows one to inflate the tube to the sphere so that we get a homeomorphism 
$$h:  N(\e,\eta) \longrightarrow  \overline{\S_\e \setminus N(\e,\eta)} \,,$$
in the obvious way: for each $z \in N(\e,\eta)$ we consider the unique integral line of $\xi$ passing by $z$; we then travel along this integral line till it hits the sphere $\S_\e$. We thus get a fiber bundle $\varphi: \S_\e \setminus L_V   \to \S^1$ with projection map $\varphi := f \circ h^{-1}$. 
The hard part is having one such vector field that further satisfies the third condition. This grants that 
the projection map in  (\ref {Fibration Thm., version1}) can be taken as  $\varphi = {f}/{|f|}$ and the two fibrations are equivalent. \qed

\vv

Having these two equivalent fiber bundles associated to a map-germ brings great richness. The first fibration is  interesting for topology and differential geometry. This has important relations with knot theory, open-book decompositions, contact and symplectic geometry. 
The second fibration lends itself more naturally to generalizations, and this has strong relations with algebraic geometry, as it exhibits the special fiber $V$ as the limit of a flat family of complex manifolds that degenerate to $V$. 

%%%%%OJO: Milnor number

%%%%%%%%%%%%%%
%%%%%%%%%%%%%%%
%%%%%%%%%%%%%%
%%%%%%%%%%%%%%%

%%%%%%%%%%%%%%%%
%%%%%%%%%%%%%%%%

\section{\bf Extensions and refinements of Milnor's fibration theorem}\label{sec. pencil}

%%%%%%%%%%%%%%%%%%%%
A natural extension of Milnor's theorem is due to  Hamm in \cite{Hamm} (see also \cite{Lo2}) for isolated complete intersection singularities, ICIS for short. This means a local complete intersection germ 
$$ f = (f_1,...,f_k) : \C^{n+k} \to \C^k \,,$$
such that $V:= f^{-1}(0)$ has an isolated singularity at the origin. One has a fibration over the regular values of $f$ sufficiently near the origin, where the fiber $F_t$ is the intersection of the complex manifold $f^{-1} (t)$ with a small ball around the origin in $\C^{n+k}$. 
Notice that in this case the set of critical values has in general complex codimension $1$ in $\C^k$ (see for instance \cite[Theorem 1.1]{Ples-Wall}).
%\footnote{referencia}

It is noticed in \cite{Lo2} that given an ICIS as above, one can always find good representatives of these singularities, which means that the first $k-1$ equations define an ICIS $\mathcal W$ of one dimension more, and  the last equation defines an isolated singularity  hypersurface germ in $\mathcal W$. Therefore we may see this as a special case of  the following theorem from \cite{Le1} (see  \cite{CSS1} for the equivalence of the two fibrations):

\begin{theorem}\label{Fibration-general}{\rm [L\^e D\~ung Tr\'ang]}
Let $X$ be an analytic subset of an open neighborhood  $U$ of the
origin $\0$ in $\C^n$.
Given $f\colon(X,\0) \to (\C,0)$ holomorphic with a critical point
at $\0 \in X$ {\rm (in the stratified sense \cite{Gor-Mc})}, let $V:= f^{-1} (0)$,  
$\B_\e$ a
closed ball of sufficiently small radius $\e$ around $\0 \in \C^n$
and  $\s_\e$  its boundary. Then:

\begin{enumerate}
\item Let $L_X=X\cap\s_\e$ be
the link of $X$ and let $L_V= f^{-1}(0)\cap\s_\e$ be the link of
$V$ in $X$. One has a fiber bundle:
\begin{equationth}\label{Fibration1-general}
\varphi = \frac{f}{|f|}\colon L_X \setminus L_V \longrightarrow \s^1.
\end{equationth}

\item Now choose $\e >> \delta
>0$ sufficiently small and consider the {\it Milnor tube}
$$ N(\e,\delta) = X \cap \B_\e \cap f^{-1}(\partial \D_\delta) \,,$$
where $\D_\delta \subset \C$ is the disc of radius $\delta$ around $0
\in \C$. Then
\begin{equationth}\label{Fibration2-general}
f \colon N(\e,\delta)  \longrightarrow \partial \D_\delta\,,
\end{equationth}
is a fiber bundle, $C^\infty$-isomorphic to the previous bundle.
\end{enumerate}
\end{theorem}

Notice that the fibers in \ref{Fibration1-general} are subsets of the link $L_X:= X \cap \s_\e$ while the fibers in  \ref{Fibration2-general} are contained in the interior of $X \cap \B_\e$,  in analogy with the classical Milnor fibrations \ref {Fibration Thm., version1} and \ref {Fibration Thm., version2}. These statements can be refined by giving a fibration on the whole ball $\B_\e$ minus the variety $V:= f^{-1}(0)$ which has the two fibrations in  Theorem \ref {Fibration-general} as subfibrations. For this we need:

\begin{theorem}[\bf The Canonical Pencil]\label{Thm:Can.Dec}
For each $\theta \in [0,\pi)$, let $\cL_\t$ be the line through $0$ in
$\R^2$ with an angle $\t$ with respect to the positive orientation of $x$-axis. Set $\,V
= f^{-1}(0)\,$ and $X_\theta =f^{-1}(\cL_\t)$. Then:

\vskip.2cm
\noindent
{\bf i)}
The $X_\theta$ are all homeomorphic real analytic hypersurfaces of $X$
with singular set $\sing(V) \cup (X_\theta \cap \sing(X))$. Their
union is the whole space $X$ and they all meet at $V$, which
splits each $X_\t$ in two homeomorphic halves.\label{it:decomposition}

\vskip.2cm
\noindent
{\bf ii)}
 If $\{S_\a\}$ is a Whitney stratification of $X$ adapted to
$V$, then the intersection of the strata with each $X_\t$
determines a Whitney stratification of $X_\t$, and   for each
stratum $S_\a$ and each $X_\t$, the intersection $S_\a \cap X_\t$
meets transversally every sphere in $\B_\e$ centered at $\0$.\label{it:stratification}

\vskip.2cm
\noindent
{\bf iii)}
There is a uniform conical structure for all $X_\t$, {\it i.e.}, there is a  homeomorphism
\begin{equation*}
 h\colon (X \cap \B_\e, V \cap \B_\e) \to \bigl(\hbox{Cone}(L_X), \hbox{Cone}(L_f)\bigr),
\end{equation*}
which restricted to each $X_\t$ defines a homeomorphism
\begin{equation*}
 (X_\theta \cap \B_\e) \cong \hbox{Cone}(X_\theta \cap \s_\e).
\end{equation*}\label{it:cones}

\end{theorem}

The next theorem implies that the fibrations over the circle in Milnor's theorem actually are liftings of fibrations over $\R \mathbb P^1$:

\vskip.2cm

\begin{theorem}[\bf Fibration Theorem]\label{Thm:fib.thm}
One has a commutative diagram of  fiber bundles
\begin{equation*}
\xymatrix{
(X \cap \B_\e) \setminus V\ar[r]^-{\varphi}\ar[rd]_{\Psi}  &\s^1\ar[d]^{\pi}\\
& \RP^1  }  \end{equation*}
where $\Psi(x) = (\hbox{Re}(f(x)):\hbox{Im}(f(x)))$ with fiber
$(X_\t\cap \B_\e)\setminus V$, $\pi$ is the natural two-fold covering and $\varphi(x)=\frac{f(x)}{{|f(x)|}}$. The restriction of
$\varphi$ to the link $L_X \setminus L_f $ is the usual Milnor
fibration (\ref{Fibration1-general}), while the restriction to the
Milnor tube $f^{-1}(\partial \D_\eta) \cap \B_\e$ is the
fibration (\ref{Fibration1-general})   (up to multiplication by a
constant), and these two fibrations are equivalent.
\end{theorem}

\vskip.2cm

The proof of Theorem \ref{Thm:fib.thm} follows the same line as in the case when $X$ is non-singular.
The key point is constructing an appropriate integrable vector field in the vein of Lemma \ref{equivalence-of-fibrations} 
above. When the ambient space $X$ is singular we must consider stratified vector fields and use either 
Mather's controlled vector fields \cite{Math}, or  Verdier's  rugose vector fields \cite{Ver}, which are all continuous and integrable.
The proof in \cite{CSS1} of \ref{Thm:fib.thm} also shows:

\vskip .2cm

\begin{corollary}\label{Cor:Mil.fib}
Let $f\colon (X,\0) \to (\C,0)$ be as above, a holomorphic map
with a critical point at $\0 \in X$, and consider its Milnor
fibration
$$ \varphi = \frac{f}{|f|}\colon L_X \setminus L_f \longrightarrow
\s^1 \,.$$ If the germ $(X,\0)$ is irreducible, then 
every pair of fibers of $\varphi$ over antipodal points of $\s^1$ are
glued together along the link $L_f$ producing the link of a real
analytic hypersurface $X_\t$, which is homeomorphic to the link of
$\{Re \,f = 0\,\}$. Moreover, if both $X$ and $f$ have an isolated
singularity at $\0$, then   this homeomorphism is in fact a
diffeomorphism and the link of each $X_\t$ is diffeomorphic to the
double of the Milnor fiber of $f$ regarded as a smooth manifold
with boundary $L_f$.
\end{corollary}

\section{\bf  The topology of the fibers: Milnor and L\^e numbers}\label{sec. topology}

 This section surveys known results  about the topology of the Milnor fibers. We refer to  Massey's expository articles \cite{Massey2, Massey3},  and to \cite {Massey}, for more complete accounts of the subject.

 We start with an example. 
 
 \begin{example}\label{ex. Le numbers}
 Consider the homogeneous map $f: \C^2 \to \C$ defined by $(x,y)  \mapsto xy$; this has a unique critical point at $x=0=y$. Its zero locus $V(f)$ consists of the two axis $\{x= 0\} \cup \{y=0\}$ with the origin as an isolated singularity. So its link $L_f :=  V(f) \cap \s^3$  is the Hopf link.  
 By  \cite[Lemma 9.4]{Mi1}, the Milnor fiber  $F_f$ is diffeomorphic to the whole fiber $f^{-1}(1)$, which consists of the points where $x \ne 0$ and $y = 1/x$. Hence $F_f$  is diffeomorphic to a copy of $\C^*$; in fact it is an open cylinder $\s^1 \times \R$, and it can be regarded as being the tangent bundle of the circle. In particular $F_f$ has the homotopy-type of $\s^1$.
 
 We extend these considerations to higher dimensions in two different ways. First notice that we can make the change of coordinates  $(x,y) \mapsto (z_1,z_2)$ with $z_1= (x+iy)$ and $z_2= (x-iy)$. In these new coordinates the above map becomes $z_1^2  + z_2^2$ and we may consider, more generally, the homogeneous polynomial:
$$f(z_0,\cdots,z_n) = z_0^2 + \cdots + z_n^2 \,.$$
The link $L_f$ consists of the points where one has:
$${\rm Re}(z_0^2 + \cdots + z_n^2 ) = 0  \quad , \quad {\rm Im}(z_0^2 + \cdots + z_n^2 ) = 0  \quad \hbox{and} \quad |z_0|^2 + \cdots + |z_n|^2  = 1 \;.$$
Hence $L_f$ is diffeomorphic to the unit sphere bundle of the $n$-sphere $\s^n$.
By Milnor's Lemma 9.4 in \cite{Mi1}, the Milnor fiber  $F_f$ is diffeomorphic to the set of points where $z_0^2 + \cdots + z_n^2 =1$, {\it i.e.}, 
$${\rm Re}(z_0^2 + \cdots + z_n^2 ) = 1  \quad \hbox{and}  \quad {\rm Im}(z_0^2 + \cdots + z_n^2 ) = 0 \;.$$
This describes the tangent bundle of the $n$-sphere $\s^n$, and  $F_f$ actually is the corresponding open unit disc bundle. In particular $F_f$ has the sphere $\s^n$ as a deformation retract and therefore $F_f$  has non-trivial homology only in dimensions $0$ and  $n$;  in these dimensions its integral homology is isomorphic to $\Z$.

Starting again with the initial example,  consider now the map $f: \C^3 \to \C$ defined by $(x,y,z)  \mapsto xyz$.  Its zero set $V(f)$ consists of the coordinate planes $\{x= 0\} \cup \{y=0\} \cup \{z=0\}$, with the three axes as singular set. The Milnor fiber $F_f$  is diffeomorphic to $\{xyz= 1\}$, {\it i.e.},  $x \ne 0$, $y \ne 0$  and $z = 1/xy$. Therefore $F_f$ is diffeomorphic to $\C^* \times \C^*$ and it has the torus $\s^1 \times \s^1$ as a deformation retract. So $F_f$ now has non-trivial homology in dimensions $0, 1, 2$.
\end{example}

We know from \cite{Mi1} that the Milnor fiber $F_f$ of an arbitrary holomorphic map-germ $(\C^{n+1},\0) \to (\C,0)$  has 
the homotopy-type of a finite CW-complex of middle-dimension $n$. This follows too from \cite{AF} since $F_f$ is a Stein manifold and, perhaps moving the origin $\0$ slightly if necessary, the square of the function distance to $\0$ is a strictly plurisubharmonic Morse function on $F_f$, so one has severe restrictions on the possible Morse indices.

If we further assume that  $f$ has an isolated critical point at $\0$, then Milnor  used Morse theory  to show \cite[Lemma 6.4]{Mi1} that $F_f$ is $(n-1)$-connected. 
Lefschetz duality, together with  the above observations about the homology of $F$,  implies that in this setting the fiber $F_f$ has the homotopy-type of a bouquet of $n$-spheres. 
 The number  $\mu= \mu(f)$  of  such spheres is  now called the Milnor number of $f$.  
  This statement can be made stronger in at least two ways.

 On the one hand, we know from \cite{Le-polyhedros} that there exists in $F_f$ a
 polyhedron $P$ of middle dimension, with the homotopy-type of a bouquet of spheres of dimension $n$, which is a deformation retract of $F_f$, and there is a continuous map $F_f \to V$ that carries $P$ into $\0$ and is a homeomorphism in the complement of $P$. So $F_f$ can be thought of as being also a topological resolution of the singularity. This is implicit in \cite{Pham}  for the Pham-Brieskorn polynomials; the polyhedron $P$ is the so-called ``join of Pham''.  The theorem in general is in \cite{Le-polyhedros} with a sketch of its proof. A  complete proof is given  in \cite{Le-Me}. 
 
 On the other hand, as proved in \cite{Mi1,Le-Pe}, the fiber $F_f$ is in fact diffeomorphic to a $2n$-ball with $\mu$ $n$-handles attached.
 To explain this we
first recall S. Smale's  classical process of  ``attaching handles''.
To attach a $p$-handle to an $m$-manifold $M$  we assume one has a smooth embedding $\iota$ of 
$S^{p-1}\times D^{m-p}$ into the boundary $\partial M$. Set  $H^{p}=D^{p}\times D^{m-p}$ and define a manifold 
$  M\cup _{f}H^{P}$ by taking the disjoint union of $M$ and $H^{p}$ and identify $S^{p-1}\times D^{m-p}$ with its image by  $\iota$. We think of $  M\cup _{f}H^{P}$ as being obtained from $M$ by attaching a $p$-handle; the integer $p$ is the index of the handle.

 Milnor noticed in \cite {Mi1} that in high dimensions,  Smale's h-cobordism theorem and the fact that the fiber $F_f$  has the homotopy-type of a bouquet of $n$-spheres, actually imply that   $F_t$  is  diffeomorphic to a $2n$-ball with $\mu$ $n$-handles attached. This claim also holds for $n=1$.    The only case left open was $n=2$; this was done 
by L\^e and Perron in \cite{Le-Pe} by a different method and  their proof actually works in all dimensions. This introduces an important technique  which in fact is a first step towards the celebrated {\it L\^e's carroussel}. The idea is to consider an auxiliary function $\ell : \C^{n+1} \to \C$, which is linear and ``sufficiently general'' with respect to $f$. The 2 maps together determine a map-germ
$$ \varphi= (f, \ell) : (\C^{n+1},\0) \to (\C^2,0) \,,$$
and the Milnor fiber of $f$ corresponds to the inverse image of an appropriate line $\ell =c$. This allows us to reconstruct $F_t$ by looking at the slices determined by the level hyperplanes of $\ell$ (see Theorem \ref{Le attaching} below).  This brings us  to the remarkable theory of ``polar varieties'' developed by Bernard Teissier and L\^e D\~ung Tr\'ang in the 1970s.
We  define first the relative polar curve of $f$ with respect to a linear form, $\Gamma^1_{f,\ell}$  (see for instance 
\cite{Le-01, Le-02, Le:VCCAS, Te1, Te4}).

Given  $f$ and $\ell$ as above, as a set the curve $\Gamma^1_{f,\ell}$ is the union of those components in the critical set of $(f,\ell)$  which are not  in $\Sigma f$, the critical points of $f$. In other words, assume we have coordinates $(z_0, \cdots, z_n)$ so that the linear function is $\ell = z_0$ is  ``sufficiently general''.  Then the critical locus of $(f, \ell)$ is $V({\partial f}/{\partial z_1}, \cdots,{\partial f}/{\partial z_n})$, the set of points where ${\partial f}/{\partial z_i} = 0$ for all $i= 1, \cdots,n$. Now write the cycle represented by $V({\partial f}/{\partial z_1}, \cdots,{\partial f}/{\partial z_n})$ as a formal sum over the irreducible components:
$$\Big [ V \big (\frac{\partial f}{\partial z_1}, \cdots,\frac{\partial f}{\partial z_n}\big )\Big ]  \, =\, \sum n_i [V_i]
$$
Then $\Gamma^1_{f,\ell}$, as a cycle, is defined by:
$$\Gamma^1_{f,\ell} \, = \, \sum_{V_i \nsubseteq \Sigma f} n_i [V_i]
$$

Now, given $f$, we need to choose the linear form $\ell$ to be  ``general enough''. For this, let us equip $f$ with a {\it good stratification} at the singular point $\0$ in the sense of \cite{Ha-Le}. This means an analytic  stratification of a neighbourhood $U$ of $\0$, such that $V(f)$ is a union of strata, the regular set of $V(f)$ is a stratum and one has Thom's $a_f$-condition with respect to $U \setminus V(f)$.  We further need $\ell$ to define a {\it prepolar slice} for $f$ at $\0$ with respect to the good stratification $\{S_\alpha\}$. This means that the hyperplane $H = \ell^{-1}(0)$ meets transversely  all the strata in $U$, except perhaps the stratum $\0$ itself.

If $H$ is a prepolar slice for $f$ at $\0$ defined by a linear form $\ell$, one has that the intersection number  $\Big ( \Gamma^1_{f,\ell} \, \cdot \, V(f) \Big )$ is finite and well-defined.
We may now state L\^e's attaching theorem: 

\begin{theorem}\label{Le attaching}
Given  a map-germ $f$ as above, with a possibly non-isolated critical point at $\0$,  consider its
local Milnor fiber $F_f$ at $\0$. If $H$ is a prepolar slice for $f$ at $\0$, then $F_f$ is obtained from the Milnor fiber of the slice $V(f) \cap H$ by attaching a certain number of $n$-handles. The number of such handles is the intersection number  $\Big ( \Gamma^1_{f,\ell} \, \cdot \, V(f) \Big )$.
\end{theorem}

 This leads naturally to the definition of the L\^e numbers, introduced by  Massey in \cite {Massey0, Massey}. The idea  is  similar to the above discussion. A linear functional $\C^{n+1} \to \C^r$ gives rise to a polar variety relative to $f$, determined by the points of non-transversality of the fibers of $\ell$ and $f$. Massey showed that this gives rise to certain local analytic cycles, that he called the L\^e cycles, that depend on the choice of the linear functional $\ell$, but they are all equivalent when the form is ``general enough''. These cycles encode deep topological properties of the Milnor fibration. We denote these by $\Lambda^k_{f,\ell}(\0)$. Up to equivalence, these are independent of the choice of the prepolar slice.

 Each of these analytic cycles has a certain local multiplicity: these are the (generic)  L\^e numbers 
 $\lambda^k_{f,\ell}(\0)$. If the singularity is isolated, then there is only  one generic  L\^e number and it coincides with the Milnor number. Massey's  theorem  (see \cite[Theorem 3.1]{Massey3}; also \cite {Massey0,  Massey1, Massey})  tells us how to build up the Milnor fiber by successively attaching handles of various dimensions. This is the second statement in  Theorem \ref {topology-fiber} below, that essentially summarizes the above discussion:

\begin{theorem}\label{topology-fiber} Let  $f:(\C^{n+1}, \0) \to (\C,0)$ be a holomorphic  map-germ and let $F_f$ be its Milnor fiber.
Then  $F_f$  is a parallelizable complex Stein manifold that can be regarded as the interior of the compact manifold obtained by attaching to it  its boundary; $F_f$ has the homotopy-type of a CW-complex of middle dimension and:
\begin{enumerate}
\item If the defining function $f$ has an isolated critical point, then:
\begin{itemize}
\item The boundary of $F_f$ is isotopic to the link $L_f$. 
\item $F_f$ has the homotopy type of a bouquet $\bigvee \S^n$ of spheres of middle dimension, and it actually is diffeomorphic to a closed ball $B^{2n}$ to which we attach $\mu$ $n$-handles. The number $\mu$ is called the Milnor number of $f$.
\item The number $\mu$ can be computed as the intersection number:
$$\mu \, = \, {\Large \rm dim}_\C  \, \frac{{\mathcal O}_{n+1,\0}}{\rm{Jac} \,f} \,,$$
where ${\mathcal O}_{n+1,\0}$ is the local ring of germs of holomorphic functions on $\C^{n+1}$ at $\0$;  $\rm{Jac} \,f$  is the Jacobian ideal, generated by the derivatives $\big ({\partial f}/{\partial z_0},..., {\partial f}/{\partial z_n}\big )$.

\item There is in $F_f$  a polyhedron $P$ of middle dimension, which is a deformation retract of $F_f$, and there is a continuous map $F_f \to V$ that carries $P$ into $\0$ and is a homeomorphism in the complement of $P$.
\end{itemize}

\item If the defining function $f$ has a non-isolated critical point:
\begin{itemize}
 \item If  the complex dimension $s$ of its critical set is $s \le  n -2$, 
then $F_f$ is obtained up to diffeomorphism,  from a $2n$-ball by successively attaching $\lambda^{n-k}_{f,\ell}(\0)$ $k$-handles, where  $n-s \le k \le n$ and $\lambda^{n-k}_{f,\ell}(\0)$  is the $(n-k)^{th}$ L\^e number.
\item  If  the complex dimension of its critical set is $s = n -1$, then $F_f$ is obtained up to diffeomorphism,  from a real $2n$-manifold with the homotopy type of a bouquet of 
$\lambda^{n-1}_{f,\ell}(\0)$ circles, by successively attaching $\lambda^{n-k}_{f,\ell}(\0)$ $k$-handles, where  $2 \le k \le n$.
\end{itemize}
\end{enumerate}
\end{theorem}

The literature about this topic is vast and includes  important bouquet theorems by D. Siersma  \cite{Siersma}  and M. Tib\u ar  \cite{Tib1} for functions defined on singular spaces,  both of these theorems being reminiscent of L\^e's attaching theorem \ref{Le attaching}.

\begin{remark}[Vanishing cycles, monodromy and the Milnor Lattice]
 The fibration theorem \ref{Fibration Thm., version2} 
 tells us that 
 the Milnor fibers can be regarded as a 1-parameter flat family $\{F_t\}$ of complex manifolds that degenerate to the special fiber $F_0 = V := f^{-1}(0)$. Since $V$ is a topological cone, this means that all the homology groups of $F$ vanish in the limit. In particular, if the critical point of $f$ at $\0$ is isolated, then 
 Theorem \ref{topology-fiber} says that  the only interesting homology group of the
Milnor fiber $F$ is in dimension $n$ and it is generated by $\mu(f)$  cycles of dimension $n$, which are called {\it the vanishing cycles}. This group,
$$\mathcal L(f) := H_n(F, \Z) \cong \Z^{\mu(f)} \,,$$
is naturally equipped with a $(-1)^n$-symmetric bilinear form $\langle\,,\, \rangle$ coming from the
intersection of cycles. This group $\mathcal L(f)$, together with this additional structure is called the Milnor lattice of the singularity. The literature about it is vast; we refer to Dimca's book \cite {Dimca} for an account of this subject.

The $E_8$ lattice is an example of a unimodular lattice which is rather  famous in singularity theory as well as in low dimensional topology. This is the Milnor lattice of the singularity $z^2 + z_2^3 + z^5_3 = 0$ in $\C^3$, whose link is the celebrated Poincar\'e homology 3-sphere. 

A cornerstone in the study of Milnor lattices is the monodromy of the fibration. For this, it is useful to  consider   distinguished bases of the vanishing cycles. We refer for this to  \cite{Dimca, Gab, Ebel,Ebel2, Saito, Sebas-Thom}.
Another  viewpoint for studying the monodromy is via its zeta-function; see for instance \cite{AC1, AC2, Var0, Oka-zeta, GLM1} for some classical material, or \cite{ACLMe, Cau-Veys, Bor-Veys} and the bibliography in these for more recent work on the subject. 
\end{remark}

\section{\bf On  the Milnor number}\label{Milnor number}

Given a  map-germ $(\C^{n+1},\0) \buildrel {f} \over {\to} (\C,0)$  with an isolated critical point at $\0$, we know from the previous section that $F$ has the homotopy type of a bouquet $\bigvee \s^n$ of spheres of middle dimension. Therefore all its homology groups vanish except in dimension $0$ where it is $\Z$,  and in dimension $n$ where it is free abelian of rank $b_n$.  The Milnor number $\mu := \mu(f)$  is defined as the $n^{th}$- Betti number $b_n(F)$ .

Theorem 8.2 in Milnor's book says that $\mu(f)$  equals the multiplicity of the map-germ $f$ at $\0$, which equals the local Poincar\'e-Hopf index of its gradient vector field $\nabla f(z) := \big({\partial f}/{\partial z_0}(z), \cdots, {\partial f}/{\partial z_n}(z)\big)$. 
Since the vector field $\nabla f$ is holomorphic, standard arguments in algebraic geometry say this index is the intersection number in Theorem \ref{topology-fiber}:
\begin{equationth} \label{formula Milnor-number}
\mu \, = \, {\Large \rm dim}_\C  \, \frac{{\mathcal O}_{n+1,\0}}{\rm{Jac} \,f} \,. \end{equationth}
It follows that the Euler characteristic of the fiber $F$ is
$\,\chi(F) \,=\, 1 + (-1)^n \mu \;.$
If we think of $F$ as being a compact manifold with boundary the link $L_f \cong \partial F$, then
the theorem of Poincar\'e-Hopf for manifolds with boundary says that  $\chi(F)$ is the total Poincar\'e-Hopf index of a vector field in $F$ that  points outward at each point of the boundary.

We know from \cite {Hamm} that one has a Milnor fibration for isolated complete intersection germs (ICIS):
$$ f:= (f_1,\cdots,f_k) : (\C^{n+k},\0) \to (\C,0)\;.$$
Hamm  proved  that in this setting the Milnor fiber also has the homotopy type of a bouquet of spheres of middle dimension. So ICIS germs have also a  well-defined Milnor number, 
defined as the rank of the middle-homology of the Milnor fiber. 
It is thus natural to search for an algebraic expression for the Milnor number 
of ICIS germs in the vein of \ref {formula Milnor-number}. This is known as the  L\^e-Greuel formula for the Milnor number:

\begin{theorem}\label{Le-Greuel}
 If
$f_1,\cdots,f_k$ and $g$ are holomorphic map germs $(\C^{n+k+1},\0) \to (\C,0)$ such that $f=(f_1,\cdots,f_k)$ and $(f,g)$ define isolated complete intersection germs, then their Milnor numbers are related by:
\begin{equation}\label{Le-Greuel}
\mu(f) + \mu(f,g) = {\rm dim}_\C \, \frac{{\mathcal O}_{n+k,\0}}{(f, \mbox{Jac}_{k+1}({f,g}))}\,,
\end{equation}
where $\mbox{Jac}_{k+1}({f,g})$ denotes the ideal generated by the determinants of all  $(k+1)$ minors of the corresponding Jacobian matrix.
\end{theorem}

This formula was  proved independently  by L\^e   \cite{Le2} and G.-M. Greuel  \cite{Greuel}. 
At about the same time  Teissier  proved \cite[Proposition II.1.2] {Te1},  a {\it ``formule de restriction''}, which is the same theorem in the case where one of the two functions is linear;  this is known as ``Teissier's Lemma''. 
We also have the celebrated {Teissier
sequence} of numbers of an ICIS, which are the Milnor numbers of the
corresponding complete intersection germs one gets by taking linear
slices of various dimensions.
This is briefly discussed in the following section.
We refer to  \cite{CMSS, CGS, Du-Gru} for other recent viewpoints on the L\^e-Greuel formula for the Milnor number. 

\vv
There are two questions that arise  naturally:

\begin{question}\label{question mu-for-non-iso}
 {What is}  {(or what ought to be)}   the Milnor number of a non-isolated hypersurface or complete intersection singularity?
\end{question}

\begin{question}\label{question mu in gral}
 {What is}  {(or what ought to be)}  { the Milnor number of an isolated complex analytic singularity $(V,P)$ which may not be an ICIS?}
\end{question}

There is a vast literature about both of these questions, with different viewpoints. Concerning the first of these,  there are two particularly interesting viewpoints. One, due to D. Massey, is via the L\^e numbers  of $V$,  explained in Section \ref {sec. topology} above.  These describe the topology of the local Milnor fiber (cf. Theorem \ref {topology-fiber}) and when the singularity is isolated, there is only one L\^e number and it coincides with the Milnor number.

Another viewpoint is that of Parusi\'nski in \cite{Pa1}. He considers a compact complex manifold $M$ and a codimension 1 subvariety $V$ defined by a holomorphic section of some line bundle $L$ over $M$. There is a Milnor number $\mu(S;V)$ associated to each connected component of its singular set. This has important properties; some of these are:

i)  it coincides with the usual Milnor number when the singularity is isolated; 

ii)   its total index in $V$ equals the 0-degree Fulton-Johnson class (cf. Section \ref{sec. Milnor classes}); 

iii) if  we can approximate $V$ by a family of manifolds $\{V_t\}$ defined by sections of $L$ which are transversal to the zero section, then $\mu(S;V)$ measures the change in the Euler characteristic as the $\{V_t\}$ degenerate into $V$; 

iv) for complex polynomials $f:\C^m \to \C$, it measures the change  in the Euler characteristic of the  fibers at infinity   \cite{ALM}.

Concerning Question \ref{question mu in gral}, for 
 a reduced complex curve  $(X,P)$, its 
 Milnor number is defined in \cite{Buch-Gr}  by:
  $$\mu_{\rm BG} : = {\rm dim}\,  \frac{\omega_X}{d \mathcal O_X}\,,$$ 
 where $\omega_X$ is Grothendieck's dualizing module (cf.  \cite {Mond-Straten}). 
 When the germ $(X,P)$ is smoothable, it equals the first Betti number of the smoothing, so for ICIS it coincides with the usual Milnor number. We recall what smoothable means:
 
 \begin{definition}\label{def. smoothable}
 A normal isolated singularity germ $(X,P)$ of pure dimension $n$ 
 is {\it smoothable} if there exists 
 an $n+1$-dimensional complex analytic space $\mathcal W$ with an isolated normal singularity that we also denote $P$, and a 
flat morphism:
$$\mathcal G: (\mathcal W,P) \longrightarrow (\C,0) \,,$$
such that ${\mathcal G}^{-1}(0)$ is $X$ and  ${\mathcal G}^{-1}(t)$ is non-singular for all $t \ne 0$ with $|t|$ sufficiently small. 
 \end{definition}
 
 The fibers of a smoothing  give a similar picture to that of the Milnor fibration of a hypersurface singularity, the difference being that  $\mathcal W$ may now be singular.
 
 There is in \cite {EGS} an alternative definition of the Milnor number for curves  using indices of  1-forms, denoted $\nu$, that works in all dimensions and in the case of curves this gives:
   $$\nu : = {\rm dim}\,  \frac{\Omega_X}{d \mathcal O_X} \,,$$
where $\Omega_X$ is the canonical module of holomorphic 1-forms. This  invariant $\nu$ is by definition the difference between the GSV and the radial indices of a holomorphic 1-form on  $V$  with an isolated singularity at $\0$, and this is well-defined in all dimensions (see Subsection \ref {sec. indices} below).

In dimension 2,  Greuel and Steenbrink  proved in \cite {Gr-Sten} that if an isolated surface singularity $(V,\0)$  is normal and Gorenstein, then its 1st Betti number $b_1(V,\0)$ vanishes and $b_2(V,\0)$ is independent of the choice of the smoothing. Hence every such germ has a well-defined Milnor number $\mu_{\rm GS}(V,\0) := b_2(V,\0)$. One has the Laufer-Steenbrink formula for the Milnor number, proved in \cite{Laufer} for hypersurfaces (the same proof works for ICIS), and by Steenbrink  \cite{Sten} in general:

\begin{theorem} \label{Laufer formula} Assume that $(V,0)$ is a smoothable Gorenstein 
       normal surface singularity.
Then:
$$\mu_{\rm GS}  +1 \,=\, \chi (\widetilde V) + K^2 + 12 \rho_g(V) $$
where   $\chi (\widetilde V)$ is the usual Euler characteristic of  a good resolution, 
 $K^2 $ is the self-intersection number of the canonical class of $\widetilde V$ and 
$\rho_g : = {\rm dim} \, H^1(\widetilde V,\mathcal O) \,$ is the geometric genus.
\end{theorem}

For surfaces, the Milnor number $\nu$ defined in \cite {EGS} seems to coincide with $\mu_{\rm GS}$.

\section{\bf Indices of vector fields and the Milnor number}\label{sec. indices}
The Poincar\'e-Hopf local index of a vector field on a smooth manifold is a fundamental invariant that has given rise to a vast literature.
 In the case of  vector fields on singular varieties,  there are  several extensions of this concept, each having its own properties and characteristics (see \cite{BSS3} for a thorough account of the subject).

In the case of manifolds, up to sign issues,  it is essentially the same to consider indices of 1-forms or vector fields. In the case of singular varieties, this is no longer the case, though the two theories are parallel and have many similarities. The case of 1-forms has been developed mostly by W. Ebeling and S. Gusein-Zade (see for instance  \cite{EG1,EG4, EG5, EGS} or \cite[Chapter 9]{BSS3}). Here we focus on vector fields.

 Let $(V,\0)$ be an analytic variety of dimension $n$ (say irreducible and reduced), for simplicity  with a unique  singular point in an open set $U$ in some $\C^m$, and let $v$ be a vector field tangent to $V$, singular only at $\0$, restriction of a continuous vector field in $U$. Then one has several possible notions of its local index. Some of these are:
 \begin{enumerate}
 \item the radial index;
 \item the homological index;
  \item the local Euler obstruction.
 \end{enumerate}
 If we further assume that $V$ is an ICIS, then one has also the GSV index and the virtual index.

\vv
Let us say a few words about these. The  Euler obstruction is discussed on its own in the next subsection.

\vv

 \noindent {\bf i) The radial index}.
  This  is a reminiscent of the Schwartz index in \cite {Sch1}, defined only for radial vector fields. The radial index was introduced  in \cite{KT} and then, independently, in \cite{ASV, EG1}.  This   is defined for continuous stratified vector fields on arbitrary compact (real or complex) varieties equipped with a Whitney stratification (and even more generally, see \cite {KT}), and it measures the lack of radiality of the vector field:
$${\rm Ind}_{\rm rad}(v; \0, V) \, = \, 1 + d(v,v_{rad}) \,,$$
where $d(v,v_{rad})$ is the total Poincar\'e-Hopf index of a vector field in a cylinder in $V$ bounded by two copies of the link, and which is radial in the ``smaller'' link and coincides with $v$ in the ``larger'' link. 

The extension of this index for vector fields on varieties with non-isolated singularities is straight-forward (see \cite{ASV, BSS3}).

\begin{figure}
\centering
\includegraphics[height=6cm ]{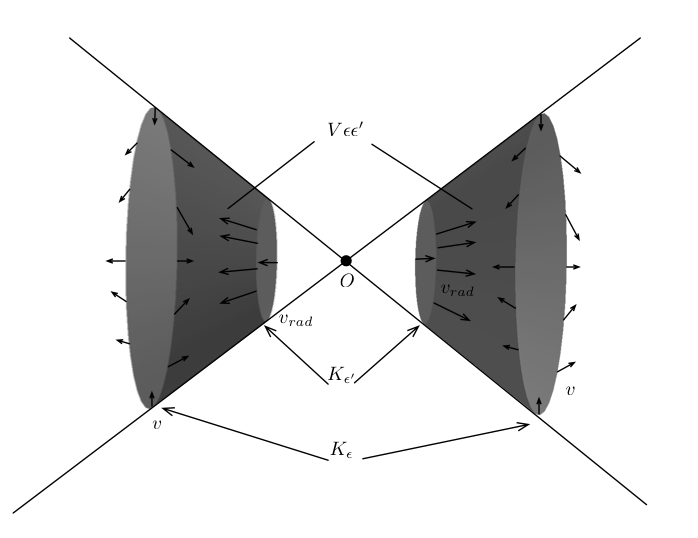}     \includegraphics[height=6.3cm ]{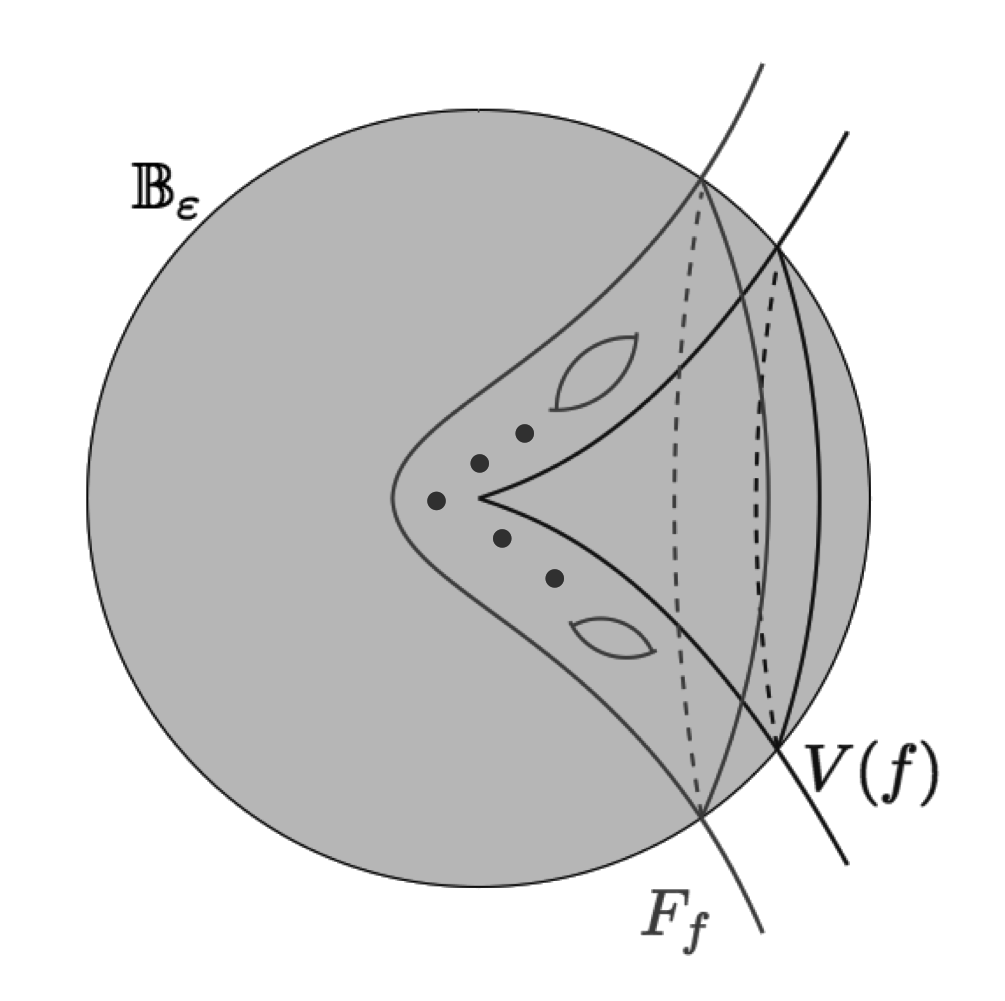}
\caption{The radial and GSV indices.}
\end{figure}

\vv
 \noindent {\bf ii) The GSV-index}. 
  The GSV-index was the first ``index of vector fields'' in the literature defined in general, since the Schwartz index and the local Euler obstruction were  defined only for radial vector fields. 
 To define the GSV-index we need $V$ to be a complete intersection. For simplicity we assume    $(V,\0)$ is a hypersurface defined by some map $f$ in an open neighborhood $U$ of $\0$ in $\C^{n+1}$. Given a  continuous vector field $v$ tangent to $V$ and singular only at $\0$, we notice that 
 $({\overline \nabla} f, v)$ determine a continuous map from the link $L_f$ to the Stiefel manifold of complex orthonormal 2-frames in $\C^{n+1}$. Such maps have a degree and this is  the GSV-index of $v$, ${\rm Ind}_{\rm GSV}(v; \0, V)$, see \cite{GSV, BSS3} for details. 
 
 The GSV-index can be interpreted as follows: we can always assume that $v$ is the restriction of  a continuous vector field $\tilde v$ on $U$, which is tangent
 to  the Milnor fibers $F_t$  for all $t$ sufficiently near to $0 \in \C$, and the restriction $v_t$ of $\tilde v$ to each $F_t$ has finitely many singularities. As $t$ tends to $0$, the fibers $F_t$ degenerate to the special fiber $V$ and the vector fields $v_t$ degenerate to $v$. If $V$ is irreducible, then ${\rm Ind}_{\rm GSV}(v; \0, V)$ equals the total Poincar\'e-Hopf index of each $v_t$.
 
\vv 
   One has (see \cite{BSS3}):
 \begin{proposition}\label{mu and indices}
 Given an isolated complete interscetion germ $(V,\0)$ and a continuous vector field $v$ on $V$ singular only at $\0$, the difference of the radial and GSV-indices is independent of the vector field and equals the Milnor number up to sign:
 $$\mu(f) \, = \, (-1)^{n} \big({\rm Ind}_{\rm GSV}(v; \0, V) - {\rm Ind}_{\rm rad}(v; \0, V) \big) \,.$$

 \end{proposition}

\vv
 \noindent {\bf  iii) The  virtual index.}
The virtual index was defined in \cite{LSS} for $v$ holomorphic and extended to continuous vector fields in \cite{SS3}. The idea is very simple. Recall first that the classical Chern-Weil theory tells us how to construct the Chern classes of complex manifolds out from a connection. Recall too 
that the Gauss-Bonnet formula allows us to identify the Euler characteristic $\chi(M)$ of a compact complex $m$-manifold $M$ with $c_m(M)[M]$, the top Chern class of $M$ evaluated on the orientation cycle. When we have a vector field with singularities on $M$, we can follow the classical Baum-Bott theory and construct a special connection, which around each connected component of the singular set $S$  of $v$, is $v$-trivial. This yields a curvature form on $M$ which determines a representative of $c_m(M)$ that vanishes
away from a regular neighborhood of $S$. We thus get an expression for $c_m(M)[M]$ which is localized at $S$. When $S$ consists of isolated points, the contribution of $v$ at each singularity is the  Poincar\'e-Hopf local index of $v$. The point is that exactly the same idea goes through for vector fields on a complete intersection $V$ in a compact complex manifold $M$, defined by a regular section of a rank $k$ vector bundle $E$ over $M$. Now we have the virtual tangent bundle of $V$, which by definition is $\tau V:= TM|_V - E|_V$. This virtual bundle has well defined Chern classes, and just as above, a continuous vector field $v$ on $V$ allows us to localize the top Chern class of this virtual bundle at the singular set $S$ of $v$.
What we get is the {\it virtual index} of $v$ associated to each connected component  of $S$. When the singularities of both $V$ and $v$ are isolated, the virtual index coincides with the local GSV-index (see \cite{LSS, BSS3}). One has \cite[Theorem 6.2]{BSS3}:

\begin{theorem} Let $V$ be a global complete intersection in a complex manifold $M$ and let $S$ be a connected component of the singular set of $V$. Let $v$ be a continuous vector field defined in a neighborhood of $S$ in $V$, with no singularities away from $S$. Then
the difference between the radial and the GSV-indices  is independent of $v$. Furthermore, if $V$ has codimension 1 in $M$, then this difference is Parusi\'nski's Milnor number up to sign: $\mu (S;V) = (-1)^n \Big({\rm Ind}_{\rm Vir} (v; S, V) - {\rm Ind}_{\rm rad}(v; S, V) \Big) \;.
$
\end{theorem}

\vv

 \noindent {\bf iv) The  homological index.}
Using the fact
 that  when the ambient space is smooth, the
 Poincar\'e-Hopf local index    can be interpreted as the
Euler-characteristic of a certain Koszul complex,
G\'omez-Mont  introduced in \cite {Gom}
 the notion of a {\it homological index} of holomorphic vector fields.
  Let us explain this
invariant.
Let $(V,\0)\subset ({\C}^m, \0)$ be a germ of a complex analytic
variety of pure dimension $n$, which is regular on $V \setminus
\{\0\}$.  A vector field $v$
  on $(V,\0)$ can always be defined   as the restriction to $V$ of a vector field
  $\widehat v$ in the
ambient space which is tangent to $V \setminus \{\0\}$;
   $v$ is holomorphic if $\widehat v$ can be chosen to be holomorphic. So we may write
  $v$ as $v = (a_1,\cdots, a_m)$ where the $a_i$ are restriction to $V$ of holomorphic functions on a
neighborhood of $\0$ in $\C^m$.

 A (germ of) holomorphic $j$-form on $V$ at $\0$
 means the restriction to $V$ of a holomorphic $j$-form on a neighborhood of
  $\0$ in $\C^m$.
  We denote by $\Omega^j_{V,\0} $ the space of all such forms (germs); these are the K\''ahler
differential
 forms on $V$ at $\0$.
  So $\Omega^0_{V,\0} $ is the local structure ring $\O_{(V,\0)}$
  of holomorphic functions on $V$ at $\0$ and each
  $\Omega^j_{V,\0} $ is an $\Omega^0_{V,\0} $-module. 
Now, given a holomorphic vector field $\widehat v$ at $\0 \in \C^m$
with an isolated singularity at the origin, and a differential
 form $\omega \in \Omega^j_{\C^m,\0}$, we can always contract $\omega$ by $v$ in the usual way,
 thus getting a differential  form $i_v(\omega) \in \Omega^{j-1}_{\C^m,0} $. If $v = \widehat v\vert_V$
is tangent to $V$, then contraction is well defined at the level
of differential
 forms on $V$ at $\0$ and one gets a
complex $(\Omega^\bullet_{V,\0}, v)$:
\begin{equationth}\label{1.1}
0 \longrightarrow \Omega^n_{V,\0}\longrightarrow \Omega^{n-1}_{V,\0}\longrightarrow \cdots\longrightarrow
\O_{V,\0}\longrightarrow 0\, ,
\end{equationth}\hskip -1pt
where the arrows are contraction by $v$ and $n$ is the dimension
of $V$.
 We consider the homology groups of this complex:
$$
H_j(\Omega^\bullet_{V,\0}, v) = {\rm Ker}\,(\Omega^{j}_{V,\0} \to
\Omega^{j-1}_{V,\0})/{\rm Im}\,(\Omega^{j+1}_{V,\0} \to
\Omega^j_{V,\0}) \;.
$$
Since the
contraction maps are $\O_{V,\0}$-module maps, this implies that if
$V$ has an isolated singularity at the origin, then the homology
groups of this complex are concentrated at $\0$, and they are
finite dimensional because the sheaves of K\"ahler differentials on $V$
are coherent. Hence it makes sense to define:   

\begin{definition}\index{Homological index !}
  The {\it homological index}
$\,{\rm Ind}_{\rm hom}(v,\0;V)$\index{$\ind_{\rm hom}$   homological index} of a holomorphic vector field $v$
on $(V, \0)$ with an isolated singularity at $\0$ is  the Euler characteristic of the above complex:
\[
{\rm Ind}_{\rm hom}(v,\0;V) = \sum_{i=0}^n (-1)^{i}
h_i(\Omega^\bullet_{V,\0},v)\,,
\]
where $h_i(\Omega^\bullet_{V,\0},v)$ is the dimension of the
corresponding homology group as a vector space over $\C$.
\end{definition}

The homological index coincides with the  GSV-index when the germ of $V$ is a complete intersection,  by \cite{BEG}. Hence in this case its difference with the radial index is the Milnor number. Yet, the homological and the radial index are defined for holomorphic vector fields on arbitrary normal isolated complex singularity germs. It is an exercise to show that their difference does not depend on the choice of vector field, so it is an invariant of the germ $(V,\0)$ and one has:

\begin{question}
If $(V,\0)$ is a normal isolated complex singularity germ which is not a complete intersection, and $v$ is a holomorphic vector field tangent to $V$ and non-singular away from $\0$, what is the difference between its homological and radial indices?
 \end{question}

\subsection{$\; \,$The  local Euler obstruction}\label{Euler obs}

MacPherson defined in \cite{MP} the notion of the  local Euler obstruction of a complex  analytic space at each of its points.  This was a key ingredient for 
constructing Chern classes for singular varieties.

The idea is that whenever we have a singular variety $V$ in some complex manifold $M$, we can consider its  Nash transform $\widetilde V$. This is the analytic space obtained
by removing from $V$ its singular set, and replacing this by all limits of tangent spaces over the regular part  $V_{\rm reg}$. So we may think of it as being a blow-up.
The space 
$\widetilde V$ itself is again singular in general, but it has the nice property that it has a natural projection $\widetilde V \buildrel {\nu} \over \to V$ which is a biholomorphism over $V_{\rm reg}$, and $\widetilde V$ is naturally equipped with a vector bundle $\widetilde T$ called the Nash bundle, which over 
$\nu^{-1}(V_{\rm reg})$ is isomorphic to the tangent bundle.

 Recall that if we equip $V$ with  a Whitney stratification, a stratified vector field means the restriction to $V$ of a continuous vector field $v$ in a neighborhood of $V$ in $M$, such that for each $
x \in V$, the vector $v(x)$ is contained in the space tangent at $x$ to the corresponding stratum. We know from \cite{BS} that the Whitney $a$-condition implies that
every stratified vector field on $V$ has a natural lifting to a section of $\widetilde T$.

We may now consider a local index ${\rm Eu}_V(v,x) $ of vector fields as follows: given $V$ as above (and equipped with a Whitney stratification), a point $x \in V$ and a stratified vector field $v$ in a neighborhood $U$ of $x$,  which is non-singular away from $x$, we lift it to a nowhere-zero section $\tilde v$ of the Nash bundle $\widetilde T$ over $\nu^{-1}(U \setminus \{x\}) $. Then ${\rm Eu}_V(v,x) $ is an integer which is the obstruction to extending $\tilde v$ as a nowhere-zero section of $\widetilde T$ over $\nu^{-1}(U) $.

A stratified  vector field $v$ in a
neighborhood of a point $x_o$ in $V$ is radial at $x_o$ if  there is small ball $\mathbb
B_\varepsilon$ in $M$, such that  $v(x)$ is pointing outward the ball at each point $x \in \mathbb S_\varepsilon \cap V$, where 
 $\mathbb S_\varepsilon:=\partial
\mathbb B_\varepsilon$ is the boundary sphere.

\begin{definition}\label{def. Euler obstruction}
 The {\it local Euler obstruction}  ${\rm Eu}_V(x) \in \Z$ of $V$ at a point $x$ is the local index ${\rm Eu}_V(v_{\rm rad},x) $ of a stratified vector field which is radial at $x$.
 \end{definition}
 
 Besides its importance  for MacPherson's proof of the Deligne-Grothendieck conjecture that we discuss in Section \ref{sec. Milnor classes},  this index has important relations with the Milnor number when $V$ is an ICIS (see \cite {BLS, Dubson}) and with various important invariants of singular varieties (cf. \cite{BMPS, EG4, EG5, Gon, Gru, Le-Te1, STV1}).  In particular one has Dubson's theorem \cite{Dubson}:
 \begin {theorem} Let $f: (\C^{n+1},\0) \to (\C,0)$ be holomorphic vith an isolated critical point at $\0$, and let $H$ be a generic hyperplane in $\C^{n+1}$ passing through $\0$. Set $V = f^{-1}(0)$. Then ${\rm Eu}_V(\0) $ is determined by the Milnor number of the hyperplane section $V \cap H$. More precisely:
 $$ {\rm Eu}_V(\0)  = 1 + (-1)^n \mu(V \cap H) \;.$$
 \end{theorem}
 This generalizes to varieties with arbitrary singular locus \cite{BLS}. If we cut down $V$ by a linear form which is not general enough, or else by a function with an isolated critical point on $V$ at $\0$, the above formula gets a correction term: this is the {\it Euler defect} in \cite{BMPS}.

%%%%%%%%%%%%%
%%%%%%%%%%%%%
%%%%%%%%%%%%%

\section{\bf Milnor classes}\label{sec. Milnor classes}  
Milnor classes  are another generalization of the classical Milnor number to 
varieties with arbitrary singularities. These measure the difference 
between two natural extensions for singular varieties of the classical Chern classes, namely 
the total Schwartz-MacPherson class
$c_*^{SM}(X)$ (see \cite{MP,Sch1}) and the total Fulton-Johnson class $c_*^{FJ}(X)$ (see \cite{Fu}).  By definition, the total Milnor class is:
\begin{equationth}\label{def Milnor class}
{\mathcal M}(X):=(-1)^{{\rm dim} X}\left(c^{FJ}(X)-c^{SM}(X) \right).
\end{equationth}

Given the category of compact algebraic varieties, one has the functor ${\mathcal F}(X)$ that assigns to each such variety $X$ the abelian group of $\mathbb Z$-valued constructible functions on $X$. One also has the funtor $H_*(X)$ that assigns to $X$ its singular homology with integer coefficients. 
 In 1974
MacPherson \cite{MP} settled a conjecture of Deligne and Grothendieck  stating that there is a unique natural transformation from the  functor ${\mathcal F}$ 
 to homology,
associating to the constant function $\mathbbm 1_X$ on a nonsingular variety $X$ the Poincar\'e
dual of the total Chern class of the tangent bundle $T X$ of $X$ (cf. \cite[p. 168]{Sul0} and Grothendieck's own comments in 
\cite[Note $87_1$, p. 376]{Groth}). Therefore, for a possibly singular variety $X$, the  homology class corresponding
to the function $\mathbbm 1_X$  is 
a natural candidate for a notion of Chern class. 
 It was proved in \cite {BS} that these classes agree, up to Alexander
duality, with certain classes defined earlier by M.-H. Schwartz in analogy with the classical definition of Chern classes via obstruction theory. Hence these  are known as Schwartz-MacPherson
classes. In particular, the 0-degree class is an integer and it equals the Euler characteristic, essentially by definition. A key ingredient in MacPherson's proof of the Deligne-Grothendieck conjecture is the local Euler obstruction introduced in Section \ref{Milnor number}.

On the other hand, whenever we have a  subscheme $X$ of a nonsingular variety $M$, its Fulton class   \cite{Fu} is the cap product
of  the total Chern class of $M$ with  the Segre class
of $X$:
$$c_F(X) := c(TM|_X) \cap s(X,M) \,.$$
When $X$ is an $(n-k)$ local complete intersection in $M$, defined by a regular section of some rank $k$ holomorphic bundle $E$ over $M$, the Fulton class is the total Chern class of the virtual bundle $\tau X:= TM|_X - E|_X$ and it coincides with the Fulton-Johnson class previously defined in \cite{FJ}. Fulton proved in \cite {Fu} that  $c_F(X)$  is independent of $M$.

In the complex analytic context, Milnor  classes
 are elements in the
homology group $H_{2*}(X, \Z)$ and in the algebraic context these can be lifted to 
elements in the Chow group $A_*(X)$.  

Just as Chern classes are related to the local index of Poincar\'e-Hopf, so too, the Schwartz-MacPherson and the Fulton Johnson classes are related to the radial and virtual indices defined in Section \ref{sec. indices} (see \cite {BSS3}).

We know from \cite {BLSS2, Su1} that for  local complete intersections 
the total Milnor class actually has
support in the singular set ${\rm Sing}(X)$ of $X$, and there is a
Milnor class in each dimension, from 0 to that of ${\rm Sing}(Y)$. Presumably this happens too in general for varieties with arbitrary singular locus, but to my knowledge, this has not yet been proved.

It follows that if  $X$  is a 
 complete intersection and its singularities  are all isolated, then there is only a $0$-degree Milnor
class, which is an integer, and we know from  \cite{SS3} 
 that this integer is the sum of the local Milnor
numbers. This is a consequence of Proposition \ref{mu and indices}, the fact that the 0-degree Schwartz-MacPherson
class can be regarded as being the total radial index of a continuous vector field  on $X$, and that for complete intersections, the 0-degree Fulton-Johnson class equals the total GSV-index of  a vector field on $X$. This justifies the name {\it Milnor classes}.

Milnor  classes spring from
\cite {Al1} and they  are an active field of current research
with significant applications to other related areas.  There is a large  literature on Milnor classes,   for instance \cite {Aluffi2, Alu-Mar, Beh, BSS3, BSY, BLSS2,  CMS1, CMS-3, CMaSS, Max, MSS, O-Y, Pa-Pr, Schu2}. 
Milnor classes   encode much
information about the varieties in question, and this is being
studied by various authors from several points of view. 

Most of the literature on Milnor classes is for hypersurfaces, though there are some recent works that throw light on the subject in the case of complete intersections. In \cite{CMS-2} the authors study the total
 Milnor class of  complete intersections $Z(s)$ defined by a regular
section $s$ of a rank $r$ holomorphic bundle $E$ over a compact
manifold $M$. It is  noticed that  $s$ determines a
hypersurface $Z(\tilde s)$ in the total space of the
projectivization $\mathbb{P}(E^{\vee})$ of the dual bundle
$E^{\vee}$, and  one has  a formula  expressing
 the total Milnor class of  $Z(s)$ in terms of the
 Milnor classes of the  hypersurface $Z(\tilde s)$. 
This means that morally, everything known for  the Milnor classes of  hypersurfaces is also  known for complete intersections.

In \cite {CMS-3} there is a surprisingly simple formula for  the total Milnor class when $X$ is defined by a finite number of hypersurfaces $X_1,\cdots,X_r$ in a complex manifold $M$, satisfying  certain transversality conditions:
$${\mathcal M}(X)=(-1)^{{\rm dim} X} \, c\left( \left(
TM|_{X}\right)^{\oplus r-1} \right)^{-1} \, \cap \,\;\Big(
c^{FJ}(X_{1}) \cdots c^{FJ}(X_{r}) - c^{SM}(X_{1}) \cdots c^{SM}(X_{r})\Big).$$

  Notice that for $r=1$ this is just the definition of the class ${\mathcal M}(X)$. There is also \cite{Sch5} where the author gives a general transversality formula that throws  light on the theory of Chern classes for complete intersections.
  
For varieties which are not complete intersections, essentially nothing is known regarding their Milnor classes. In fact, there is even an ambiguity in the definition, since it is not clear if one should  consider  the Fulton or the Fulton-Johnson class \cite{FJ}, which coincide for complete intersections; perhaps both and get two different interesting concepts.  Here is an easy to state question (cf. Question \ref{question mu in gral}):
 
 \begin{question}
 {\rm Consider an algebraic normal isolated singularity germ $(V,P)$ of dimension $n \ge 1$ which is not a complete intersection. Consider a projective compactification of it and resolve its singularities at infinity. We get a compact variety $\overline V$ with an isolated singularity. What is the difference between its Fulton and Schwartz-MacPherson classes?
 }
 \end{question}        

\section {\bf Equisingularity}\label{sec. equisingular}    

The vast field of equisingularity theory emerged after the seminal work of O. Zariski  \cite{Zar1,Zar2, Zar3, Zar4}, B. Teissier (see  for instance \cite{Te1, Te2, Te3}),   L\^e D\~ung Tr\'ang  (see  for instance \cite{Le-equi, Le-Ram, Le-Ma1}) and many others,  for instance  \cite{Bobadilla1, Bobadilla2, Ga1, Ga2, Ga3, Massey1, NOT1}. 

  The term ``equisingular''  refers to a relation of equivalence which   formalizes the intuitive idea of singularities of ``the same type''  in some sense. In general one would like to find  conditions, and perhaps numerical invariants, that grant some type of equisingularity. 
There are several important notions of equisingularity.  The most basic one is:

\begin{definition}
Two germs of reduced complex analytic hypersurfaces $(V_1,q_1)$ and $(V_2,q_2)$ in $\C^{n+1}$ are {\it topologically equisingular} if they have the same embedded topological type, {\it i.e.},  if there exist representatives of these germs $(V_i,q_i) \subset (U_i,q_i)$, $i =1,2$, with $U_i$ an open set in $\C^{n+1}$, and a homeomorphism of pairs $(U_1,q_1) \cong (U_2,q_2)$ taking $V_1$ into $V_2$.
\end{definition}

One has the celebrated Question A by Zariski in \cite{Zar5}, known as Zariski's multiplicity conjecture:

\begin{question}\label{Zariski conjecture}
Does topological equisingularity of two hypersurface germs $V_1$ and $V_2$ at $q_1$ and $q_1$ imply that 
these have the same multiplicity,  $\nu(V_1, q_1) = \nu(V_2 , q_2)$?
\end{question}

Recall that if $f: (\C^{n+1},\0) \to (\C,0)$ is a holomorphic map-germ, we can write it in terms of its Taylor series expansion:
$$ \qquad \qquad f(z) \,= \, \sum_{j=1}^{\infty} \, f^j(z) \quad, \quad \hbox{with} \quad f^j(z)  \,=  \; \sum_{\alpha_0+\cdots+\alpha_n = j} a_\alpha \, z^\alpha  \;, $$
where each $f^j$ is 
homogeneous of degree $j$, 
$\alpha = (\alpha_0,\cdots,\alpha_n) \in \mathbb N^{n+1}$, $z^\alpha = z_0^{\alpha_1} \cdots z_n^{\alpha_n}$ and $a_\alpha \in \C$.
{\it The order} of $f$ at $\0$ is the smallest degree $j$  such that $f^j $ is not identically $0$. If the germ is reduced, this coincides with its {\it multiplicity} $\nu(f,\0)$, which by definition is  the number of points of intersection, near $\0$, of $V= f^{-1}(0)$ with a generic complex line in $\C^{n+1}$ passing arbitrarily close to $\0$ but not through $\0$.

There is a vast  literature about this question, known as Zariski's multiplicity conjecture, with many partial answers. We refer to \cite{Ey} for a survey on the topic, and to \cite{Fer-Sam} for the answer to a metric version.

%%%%%%
Another important notion of equisingularity  is   Whitney regularity, which has already appeared several times in this work. 
The existence of Whitney stratifications for every analytic space
$X$ was proved by Whitney in \cite[Theorem~19.2]{Whitney:TAV} for
complex varieties, and  by Hironaka \cite{Hironaka:SubAnal} in a more 
general setting.
Thom \cite {Th1} and Mather \cite {Math}  proved that Whitney equisingularity implies local topological triviality; this is essentially a consequence of the  first Thom Isotopy Lemma
(cf. \cite{AC3, Ver}).

\begin{definition} Let $V$ be a reduced complex analytic hypersurface  in $\C^{n+1}$ and $V_\alpha \subset V$ a stratum of some Whitney stratification. We say that $V$ is {\it topologically trivial} along $V_\alpha$ at a point 
 $x \in V_\a$ if there is a neighborhood $W$ of $x$ in ${\C}^{n+1}$, homeomorphic to 
$\Delta \times U_\a$, where $U_\a$ is a  neighborhood of $x$ in $V_\a$ and 
$\Delta$ is a small closed disk through $x$ of complex dimension 
$(n+1)-\dim_{\C}V_\a$, transverse to all the strata of $V$, and such  that $W \cap V_\beta= (\Delta\cap V_\beta )  \times U_\a$ for
each stratum $V_\beta$ with $x\in \overline{V_\beta}$.
\end{definition}

In his remarkable ``Carg\`ese'' article \cite{Te1}, B.  Teissier introduced a decreasing sequence of numbers $\mu^*$, today known as the Teissier sequence, 
and  proved  that if a family has constant $\mu^*$-
sequence, then it is Whitney equisingular. The converse, that  Whitney equisingularity implies constant $\mu^*$-
sequence was proved later by Brian\c con-Speder in \cite{Bri-Sp2} (cf. \cite{Bri-Sp}).
The $\mu^*$-sequence of an isolated hypersurface singularity $(X, \0)$ in $\C^{n+1}$ is:
$$ \mu^* \, := \{ \mu^{n+1}(X), \cdots, \mu^{i}(X), \cdots, \mu^{0}(X)  \}\,,$$
where $\mu^{n+1}(X) = \mu(X)$ is the usual Milnor number and $ \mu^{i}(X)$ is the Milnor number of the intersection of $X$ with a general plane in $\C^{n+1}$ of dimension $i$, passing through $\0$.

There are several other concepts of equisingularity, for instance the Milnor equisingularity studied in \cite{Le-Ma1}, $\mu$-constant families as in L\^e-Ramanujam's theorem \cite{Le-Ram}, bi-Lipschitz equisingularity (see \cite {Mostow1, Pa2}), etc. The concept of equisingularity is also closely related to that of ``simultaneous resolution'' and work by Hironaka, Lipman and others (see \cite[p. 595]{Te2}).

 %%%%%%%%%%%%%%%%%%%

Equisingularity has proved to be a very subtle subject, with a myriad of different aspects and  open questions. For more on this subject, we refer to the literature, which is vast;  see for instance  the  survey articles \cite{Greuel-Shu, Lip} or  \cite {Bobadilla1, Bobadilla2, Bobadilla3, Massey0, Massey1, Massey, Ga5}.

  %%%%%%%%%%%%

  %%%%%%%%%

 \subsection{$\, \; \;$Lipschitz geometry of singularities}    
Recall that a  continuous map $f : Y \to Z$  between metric spaces  $(Y, d_Y)$ and $(Z, d_Z)$ is Lipschitz
if there is a constant $L \ge 1$ such that
$$L^{-1} d_Y (a, b) \le  d_Z(f(a), f(b)) \le  L \, d_Y (a, b)$$
for each pair of points $a, b \in X$. 
A metric space $Y$ admits local bi-Lipschitz
parameterizations by $\R^n$
if every point in $Y$ has a neighborhood that is bi-Lipschitz
homeomorphic to an open subset of $\R^n$. Two metric spaces $Y,Z$ are bi-Lipschitz equivalent if there a bi-Lipschitz homeomorphism $Y \to Z$.
  bi-Lipschitz classification is stronger than topological and weaker than $C^\infty$  classifications.

A topological $n$-manifold, $n \ge 2$, is a Lipschitz
manifold if it can be 
equipped with an atlas where all transition functions  are bi-Lipschitz. Such a manifold is smoothable  
if there is  further a subsystem of charts where the transition functions
are diffeomorphisms. A deep theorem of D. Sullivan in \cite{Sul2}  guarantees that in  dimensions $\ne 4$, there is  a unique Lipschitz structure on every
topological manifold [57], and this can be used to study the smoothability of  topological manifolds, which is a deep theory explored in \cite{Sieb-Sull}.

bi-Lipschitz homeomorphisms have good properties that can be of interest in various fields of research, in particular for studying  the geometry and topology of analytic spaces.

The study of bi-Lipschitz geometry of complex spaces started with Pham and Teissier  in \cite{Pham-Teiss}. Later, Mostowski in \cite {Mostow, Mostow1} studied Lipschitz equisingularity and Lipschitz stratifications in analytic sets, a notion that  grants the constancy of the Lipschitz type of the stratified set along each stratum. The 
 existence of Lipschitz stratifications for complex analytic sets was established in \cite {Mostow1}, and in the real analytic case, including semi-analytc sets, this was done  by Parusinski in \cite{Pa2, Pa3}  (see also \cite{N-Valette}).

In \cite{Bir-Mostow} the authors look at Lipschitz properties of semialgebraic sets with singularities and study the concept of normal embeddings, which has opened an important line of research. For this, notice that given an analytic subset $X$ of $\R^n$, we
have two natural metrics on $X$: one is the metric induced from the ambient space; this is called the {\it outer metric}. The other is the
{\it inner}, or  length, {\it metric} defined in the usual way in differential geometry, as the infimum of lengths of piecewise smooth curves connecting two given points. The Lipschitz equivalence in terms of the outer metric is more rigid and it implies the equivalence in inner metric, but not inversely. The set $X$ is  {\it normally embedded} if these two metrics are equivalent. 

The main result in \cite{Bir-Mostow}  states that every compact semialgebraic set is bi-Lipschitz equivalent to some normally embedded semialgebraic set. The article \cite{Bir-Alex} started the study of the bi-Lipschitz geometry of 
      complex surface singularities and in recent years 
 there has been remarkable progress in this and other related topics thanks to the work of L. Birbrair, A. Fernandes, A. Pichon,  W. Neumann, G. Valette, D. Kerner,  T. Gaffney, J. E. Sampaio and many others  (see for instance \cite{Bir1, Bir-Alex,  BFLS, Bir-Alex-Neu,  Bir-Neu-Pi, Fer-Sam, Ga4, Ke-Pe-Ru, Ne-Pi1, Ne-Pi2, Ne-Pi3, Ne-Pe-Pi, Valette} and the references therein).

Some of the important recent results in this area are the complete classification of the inner metrics of surfaces in \cite{Bir-Neu-Pi}, the proof in \cite {Ne-Pi3} that Zariski equisingularity is equivalent to bi-Lipschitz triviality in the case of surfaces, the proof  in \cite{BFLS}  that outer Lipschitz regularity implies smoothness, and the  important  partial answer in \cite{Fer-Sam} of Zariski's multiplicity conjecture (\ref{Zariski conjecture}), proving that if 
$ f, g : \C^n \to \C$ are irreducible homogeneous polynomials such that there is a bi-Lipschitz homeomorphism $h: (\C^n,V(f),0) \to (\C^n,V(g),0)$, then $f$ and $g$ have the same multiplicity at $0$.

\section{\bf Relations with other branches}\label{sec. relations}

\subsection { Fibered knots and  open-books}
The concept of {\it open-books}  was introduced by 
E. Winkelnkemper and
 we refer to his  appendix in \cite {Rn} for a clear account of the subject. An open-book decomposition of a smooth $n$-manifold $M$ 
  consists of a codimension 2 
submanifold $N$, called {\it the binding}\index{binding}, embedded in $M$ with trivial normal bundle, 
together with a fiber bundle decomposition of its complement:
$$\theta: M  - N \to \s^1\,,$$ 
satisfying  that on a tubular neighborhood of $N$, diffeomorphic to  $N \times \D^2$,  the restriction of $\theta$
to $N \times (\D^2  - \{0\})$ is the map $(x,y) \mapsto y/\|y\|$.  The fibers of $\theta$ are called 
{\it the pages}\index{pages, of an open-book} of the open-book. These 
 are all diffeomorphic and each page $F$ can be compactified 
by attaching the binding $N$ as its boundary, thus getting a compact manifold with  boundary.  

Milnor's fibration theorem grants that in the isolated singularity case, we get open-book decompositions. In fact, given 
a complex analytic space $X$ which is non-singular away from a point, say $\0$, and a holomorphic map-germ 
$$f: (X,\0) \longrightarrow (\C,0)$$ which is regular away from $\0$, by the  fibration theorem 
 one has an open-book decomposition on  the link $L_X$ of $X$:
$$ \varphi:= \frac{f}{|f|} \colon \, L_X \setminus L_f  \longrightarrow \S^1 \,, $$
where  $L_f$ is the link of $f$ in $X$ and the pages of the open-book are the Milnor fibers.

Recall that if $M$ is a smooth, closed, connected manifold, a {\it knot} in $M$ means a smooth
 codimension 2 closed, connected submanifold $N$ of $M$. If $N$ has several connected components then it 
is called {\it a link} in $M$; so the binding of an open-book is a knot (or link). The name ``algebraic knot''  was coined by L\^e D\~ung Tr\'ang in \cite{Le1} to characterize the  knots defined by an algebraic (or  analytic) equation with complex values. 

In the literature there are also real algebraic knots. For instance Perron  proved in \cite{Perron}  that the figure eight knot is real algebraic where, in analogy with the complex case, this means that it is defined by the link of a polynomial map  $\R^4 \to \R^2$ with an isolated critical point. In \cite{AK} the authors prove that ``all knots are algebraic''; yet, one  must be careful that the statement here is slightly different: although the codimension 2 real algebraic varieties that define the knots have an isolated singularity, the functions that define them may have non-isolated critical points. When this happens,  the corresponding knot is not fibered: a knot (or link) $N \subset M$ is {\it fibered} if it is the binding of an open-book
decomposition of $M$. The
concept  of fibered knots  was introduced by A. Durfee and B. Lawson  in \cite {{Du-La}}, where they use Milnor fibrations
to construct codimension 1 foliations on odd-dimensional spheres.  By Milnor's theorem, every (complex) algebraic knot is  fibered.

The classical  theory of 1-dimensional knots in the 3-sphere which are algebraic  goes back  to Brauner \cite{Brauner}. The literature on this topic is vast and we refer to
the excellent book \cite{Bri-Kn} for a clear account (see also  \cite {Wall-book} and \cite [Chapter 10]{Mi1}). Just a few words: every map germ
$$f: (\C^{2},\0) \longrightarrow (\C,0) \,,$$
has an essentially unique prime factorization $f = f_1^{a_1} \cdots f_r^{a_r}$. The zero locus $V(f)$ of $f$ is the union of the zero-loci $V(f_i)$ of the $f_i$. Each of these defines a branch of $f$,  an irreducible component of $V(f)$, and we would like to describe these branches.

%%%%%%%%%%%
As an example consider the complex polynomial $ f(z_1,z_2) \to z^p + z^q \,\,,$
for some $p,q >1$. Let $k$ be the greatest common divisor of $p,q$ and set $p' = p/k$, $q'= q/k$. Set $V= f^{-1}(0)$. Then $V$ has $k$ branches and the intersection of each branch with the unit sphere is a torus knot of type $(q',p')$, {\it i.e.}, it is wrapped in a torus  $S^1 \times S^1$ so that it goes around a parallel $q'$-times, and around a meridian $p'$-times.

It is known in general that each branch of an analytic plane curve admits a particular type of parameterization called a Newton-Puiseux parametrization. This describes the link as an iterated torus knot, and the corresponding Puiseux pairs tell us exactly how to construct it (see \cite{Bri-Kn} for details).

For $n=1$ the link of the singularity is either a circle, or one circle for each branch, and the interesting point is knowing how these knots are embedded in the sphere. For $n > 1$ the link itself has interesting topology. 
It is also interesting to study the links of holomorphic functions on complex surface singularities. These give knots and open books in the link of the surface singularity, which is a 3-manifold (see for instance \cite{Pi1, Pi2, PS2}).

\subsection { Open-books and contact structures}

Recall   that if $ M$ is an oriented $(2n - 1)$-dimensional manifold, a {\it contact structure} on $M$ is a hyperplane distribution $\zeta$ in its tangent bundle $TM$, which is locally given  by a 1-form $\alpha$ such that $\alpha \wedge  (d\alpha)^{(n-1)} \ne 0$. In this case we say that the pair $(M, \zeta)$ {\it  is a contact manifold and $\alpha$ is a contact form. }
The contact structure is called oriented if the vector bundle $\zeta$   is oriented.
If  $\alpha$ is a contact form, it is called positive if the volume form $\alpha \wedge (d \alpha)^{(n-1)}$ defines the orientation of $M$. 

There is a natural way in which contact manifolds arise in complex geometry: start with a complex manifold $X$ and a real hypersurface $M$ in it. At each point $z \in M$ we have the tangent spaces $T_z M \subset T_zX$. Multiplication 
of $T_z M $ by the complex number $i$ gives another real hyperplane $i(T_z M ) \subset T_zX$. The intersection $ \zeta_z :=T_z M  \cap 
i(T_z M ) $ is a real hyperplane in $T_z M $.

This hyperplane distribution in $M$  may or may not be a contact structure: if the real hypersurface $M$ in the complex manifold $X$  is strongly pseudoconvex, then the distribution $\zeta$ defined above  is a naturally oriented contact structure (see for instance \cite[Proposition 5.11]{PPP-1}).  Pseudoconvex means that   $M$ can be defined locally, in a neighborhood of each of its points, as a regular level of a  strictly plurisubharmonic function.

Consider now a complex analytic variety $X$ of pure dimension $n$ in some complex space $\C^m$, and let $p$ be an isolated singularity in $X$. 
 Consider the intersection of $X$ with the  spheres in $\C^m$ centered at $p$.  
 Varchenko in  \cite{Var}  noticed that the square of the distance function restricted to $X$ still is a strictly  plurisubharmonic function, so the link 
  $L_X$ is pseudoconvex and it can be equipped with a natural contact structure as above. Furthermore, this contact structure on the link is independent of the choices of the embedding and of the (sufficiently small) spheres, up to contactomorphisms (well-defined up to isotopy). This leads to  the following definition:
 
 \begin{definition}
The oriented contact manifold associated in this way to every isolated singularity germ $(X,p)$, 
 up to contactomorphisms isotopic to the identity, is called {\it the contact boundary} of $(X,p)$, and is denoted $( \partial(X,p), \xi(X,p))$. 
 \end{definition}

 An oriented contact manifold which is contactomorphic to the contact boundary of an isolated singularity is called {\it Milnor fillable}, a
 name  introduced  in \cite{CNP} in reference to  Milnor's work \cite{Mi1}.

 Every Milnor fillable contact manifold $(M,\xi)$ is holomorphically fillable, since every resolution of a singularity whose contact boundary $( \partial(X,p), \xi(X,p))$ is contactomorphic to $(M,\xi)$ gives a holomorphic filling of it. Moreover, if there is a singularity germ $(X,p)$ with contact boundary $( \partial(X,p), \xi(X,p))$ which  is smoothable (see Definition \ref{def. smoothable}), then it is easy to construct Stein representatives of its Milnor fibers and these  are Stein fillings of the contact boundary $( \partial(X,p), \xi(X,p))$
 (see \cite[Proposition 6.8]{PPP-1}). 
 
 As pointed out in \cite[p. 62]{PPP-1} there is a remarkable difference  between complex dimension 2 and higher dimensions which is  highlighted by the following two theorems. First recall the Pham-Brieskorn singularities envisaged in Section \ref{sec. initial}. One knows from \cite {Brieskorn} that in the special case
 
\begin{equationth}\label{Ustilovsky}
  z_0^{2} + z_1^{2} + \cdots +  z_{2m-1}^{2}  + z_{2m}^{p}  \, =\, 0 \; , \; m \ge 2 \;,
   \end{equationth}
   \hskip-4pt
 with $ p \equiv  \pm 1 (\hbox{mod} \; 8) \,, $ the link is always diffeomorphic to the standard $(4m-1)$-sphere. One has Ustilovsky's theorem in \cite{Ustilovsky}:

 \begin{theorem}
Varying $p \in \mathbb N$ in \ref{Ustilovsky} above, the contact boundaries of the corresponding isolated hypersurface singularities  are pairwise non-contactomorphic. One thus gets in this way infinitely many different contact structures on the $(4m-1)$-sphere which are contact boundaries of some isolated complex hypersurface singularity.
 \end{theorem}
 
 In complex dimension 2 the situation is radically different, as was proved  
 by Caubel, N\'emethi and Popescu-Pampu in \cite{CNP}:
 
  \begin{theorem}\label{theorem CNP}
Every Milnor fillable oriented 3-manifold admits a unique Milnor fillable contact structure up to contactomorphism.
 \end{theorem}

The proof of Theorem \ref{theorem CNP} uses  work by E. Giroux in \cite{Giroux}  that motivates the title for this subsection: the notion of contact structures carried by an open-book:
 
 \begin{definition}
 A positive contact structure   $\xi$ on a closed oriented manifold $M$ {\it is carried by an open book} $(N,\theta)$ if it admits a defining contact form  $\alpha$ 
 which verifies the following:
 
 $\bullet$  $\alpha$ induces a positive contact structure on $N$;
 
$\bullet$  $d \alpha$ induces a positive symplectic structure on each fiber of  $\theta$.
 
 \noindent
  If a contact form  $\alpha$ satisfies these conditions, then it is said {\it to be  adapted to} $(N,\theta)$.
 \end{definition}
 
 Giroux  proved in \cite{Giroux}  that on each 3-dimensional closed oriented manifold, every contact structure is carried by
some open book and two positive contact structures carried by the same open book are isotopic. Thus,  in order to describe a positive contact structure on a 3-dimensional closed and oriented manifold, it is enough to describe an open book which carries it. This is the strategy adopted in \cite{CNP} to prove Theorem \ref{theorem CNP}.

 There is another important line of research that springs from the following theorem of  Eliashberg  \cite{Eliashberg}:
 
 \begin{theorem}\label{thm. Eliash}
 Every Stein filling of the natural contact structure on $\mathbb S^3$ is diffeomorphic to the 4- dimensional compact ball.
 \end{theorem}

We may  naturally ask when and how this theorem extends to the contact boundaries of normal complex surface singularities. For this we may
want to characterize the Milnor fibers of a given isolated singularity (up to diffeomorphisms) amongst the fillings of the boundary of the singularity. Since the Milnor fibers of every smoothing can be choosen to be Stein fillings, it is natural to restrict  our attention to the Stein fillings of the contact boundary. In complex dimension 2 one has Theorem 6.3, so in this case one is led to asking the following questions (see \cite[Section 6.2]{PPP-1}):
\vv

\begin{itemize}

\item {\it Is it possible to characterize the Milnor fibers of the various isolated surface singularities with a given topological type among the Stein fillings of the associated Milnor fillable contact 3-manifold? 

\item Are there situations in which one gets all the Stein fillings up to diffeomorphisms as such Milnor fibers?   }
\end{itemize}

There are several important contributions to this line of research done by various authors and we refer to \cite[Section 6.2]{PPP-1} for an account of this. We finish this subsection with the following theorem from \cite{Nem-PPP} that provides a  generalization of Eliashberg's theorem \ref {thm. Eliash}:

\begin{theorem} The Milnor fibers of a cyclic quotient singularity exhaust the Stein fillings of its contact boundary up to diffeomorphism.
\end{theorem}

We refer to \cite[Section 6]{PPP-1} for a clear account of the subject discussed in this subsection, including a fairly complete bibliography and a list of interesting open questions.

\subsection { Low dimensional manifolds}\label{3-manifolds}
If $(V,\0)$ is a normal surface singularity in some affine space $\C^N$, then its link $L_V:= V \cap \s_\e$ is a closed oriented 3-manifold  with a rich geometry. 
The   interplay between  3-manifolds  theory and  complex surface singularities goes back to  F. Klein \cite {Klein} and many others, for instance \cite{Dol1, Dol2, Mum, Mi1, Mi2, Ne1, Ne2}. We refer to  \cite[Chapters 3, 4]{Seade-librosing} for a thorough discussion on that subject.

 For instance  consider the polynomial map $(\C^3,\0) \buildrel {f} \over {\to} (\C,0)$ given by:
$$ (z_1,z_2,z_3) \mapsto z^{p}_1 + z^q_2 + z_3^r \quad, \;  \; {\rm with} \; \; p,q,r \ge 2 \,.$$
The link is a  3-dimensional Brieskorn manifold.   Klein in \cite {Klein} showed that  when $1/p +1/q +1/r  > 1$, the link is diffeomorphic to a quotient of $\S^3$  divided by a discrete subgroup. For instance for the triple $(2,2,r)$ we get the quotient $\S^3/ \mathbb Z_n$ which is the lens space $L(n,1)$. For the triple $(2,3,5)$,  the group is the binary icosahedral group and  $L_f$ is Poincar\'e's homology 3-sphere.

Notice that if we order $p,q,r$ so that $p \le q \le r$ then the condition $1/p +1/q +1/r  > 1$ is satisfied only for the triples $(2,2,r)$, for every $r \ge 2$, $(2,3,3)$, $(2,3,4)$ and $(2,3,5)$. In all cases the singularity we obtain is a rational double point, also called a Klein or Du Val singularity. 

For $1/p +1/q +1/r = 1$ the only possible triples up to permutation are $(2,3,6)$, $(2,4,4)$, $(3,3,3)$;  Milnor proved in \cite{Mi2} that the links of all these singularities are quotients of the 3-dimensional Heisenberg group of all $3 \times 3$ matrices which are upper triangular, with 1s in the diagonal, divided by appropriate discrete subgroups. In all other cases we have  $1/p +1/q +1/r < 1$ and the corresponding links are quotients of the universal cover ${\widetilde {{\rm SL}}(2,\R)}$ of ${\rm PSL}(2,\R)$ divided by the commutator of the lifting of the triangle group $\langle p,q,r \rangle \subset {\rm PSL}(2,\R)$ to ${\widetilde {{\rm SL}}(2,\R)}$. We refer to \cite{Mi2} or  \cite[Chapter 3]{Seade-librosing} for details; see also \cite{Dol1, Dol2, Ne2}.

In general, it follows from work by Grauert and Mumford and also (independently) Du Val, that an oriented closed 3-manifold $M$ is orientation preserving diffeomorphic to the link of an isolated complex surface singularity  if and only if it is a Waldhausen manifold with negative definite intersection matrix (cf. \cite {Ne1}).  The following  is a classical open problem and to my knowledge there has not been any significant improvement in more than two decades (cf. \cite{Yau2}):
 
 \begin{problem}
 Characterize the Waldhausen manifolds (with negative definite intersection matrix) that appear as links of surface singularities in $\C^3$.
 \end{problem}

By  \cite{Ne1} the orientation preserving homeomorphism type of the link $L_V$ of a normal surface singularity $(V,\0)$ depends only 
on the analytic type of $(V,\0)$.  Hence every 3-manifolds invariant is an invariant of singularities and,  conversely, 
whatever invariant of 3-manifolds we want to understand, the links of surface singularities are a great source of examples.  

It is interesting to study, in particular, relations between topological  invariants of the singularity, which can be determined from the link, and  analytic invariants, for instance the Milnor number, the geometric genus and the signature of the Milnor fiber.

This begins with M. Artin  in \cite{Artin1, Artin2}, proving  that the rational singularities  are ``taut'', {\it i.e.}, they are characterized by the topology of the link. Laufer proved in \cite {Laufer} a formula for the Milnor number  of hypersurface singularities in $\C^3$ using the Hirzebruch-Riemann-Roch theorem:
$$ \mu +1 = \chi(\widetilde V) + K_{\widetilde V}^2 + 12 \rho(V)\;,$$
where  $\widetilde V$ is a good resolution, $\chi$ its topological Euler characteristic, $K_{\widetilde V}^2$ is the self-intersection number of the canonical class of the resolution, and  $\rho(V)$ is the geometric genus. This formula was extended in \cite{Sten} to all smoothable Gorenstein surface singularities.

Notice that the right hand side in Laufer's formula is well defined for non-smoothable singularities; this  is called the Laufer invariant in   \cite{Seade-Laufer}. In fact that invariant splits in two parts, $\chi(\widetilde V) + K_{\widetilde V}^2$ and $12 \rho(V,q)$. The first of these is topological, depending only on the topology of the link $L_V$. Notice that if $L_V$ is a rational homology sphere, then the first betti number $b_1(\widetilde V)$ of $\widetilde V$ vanishes and $\chi(\widetilde V) + K_{\widetilde V}^2$ essentially coincides with the invariant in a conjecture by N\'emethi-Nicolaescu  (see \cite[Remark 4.8]{Ne-Ni1}; also \cite[Subsection 2.4]{Ne-Ni2}) related to the Casson invariant conjecture that we state below.

Laufer's formula inspired Durfee's theorem in \cite {Du1}, which states that if a normal surface singularity $(V,p)$ is numerically Gorenstein, smoothable and the complex tangent bundle of the Milnor fiber $F$ of a smoothing is trivial, then the signature of $F$ can be expressed as:
$$\sigma(F) = - \frac{1}{3} \Big( 2 (\chi(V_t) - 1) + K_{\widetilde V}^2 + 2b_1(\widetilde V) + b_2(\widetilde V) \Big) \;.
$$
The proof is based on Hirzebruch's signature theorem for closed oriented 4-manifolds. This formula was completed in \cite{Se2} by proving   \cite[Conjecture 1.6]{Du1},  that the tangent bundle  of the Milnor fiber of every smoothing of a Gorenstein surface singularity is trivial. Hence  Durfee's formula applies to all Gorenstein smoothable singularities.

Just as  Laufer's formula, Durfee's signature formula can be regarded as  a consequence of the general  Atiyah-Singer index theorem. It is then natural to ask what information about invariants such as the Milnor number, the signature of the Milnor fiber and the geometric genus can be determined topologically. This brings us to \cite{ESV, Se2} and several important articles by N\'emethi {\it et al} for instance \cite{Ne-2, Ne-3, Ne-Ni1, Ne-Ni2, NeSigu}. 

Consider
 the minimal resolution $\pi: \widetilde V \to V$ of a Gorenstein surface singularity $(V,\0)$ and let $K:= K_{\widetilde V}$ be a divisor of the canonical bundle ${\mathcal K}:= {\mathcal K_{\widetilde V}}$. Assume further (with no loss of generality) that the divisor $K$ is vertical, {\it i.e.},  the support of the divisor is contained in the exceptional curve.
For each vertical divisor $D \ge 0$ set $W= 2D - K_{\widetilde V}$. Such a divisor $W$ was called characteristic in \cite{ESV},  in analogy with the classical theory of characteristic vectors and submanifolds, see for instance \cite{Fr-Kir, Ki, La-Se}, because  $W$ represents an integral homology class whose reduction modulo 2 essentially is the 2nd Stiefel-Whitney class of $K_{\widetilde V}$. Then we 
have the equalities:
$$\qquad \qquad{\rm dim} \, H^1(\widetilde V, {\mathcal O}_{\widetilde V}) \;  = \;  {\rm dim} \, H^0(-K, {\mathcal K}|_K)   \, = \,   {\rm dim} \, H^0(W, {\mathcal D}|_W) + \frac{1}{8} (W^2 - K_{\widetilde V}^2) \;,$$
where ${\mathcal D}$ is the line bundle of the divisor $D$.
Since the first term above is the geometric genus, and the last term in the right obviously is topological, it follows that $\rho_g $ is topological whenever we can find a characteristic divisor $W$ for which the integer $ {\rm dim} \, H^0(W, {\mathcal D}|_W)$ is topological.

The 3-manifolds which are links of complex surface singularities  carry canonical contact and ${\rm Spin^c}$ structures inherited from the holomorphic structure on $V$. If the singularity is Gorenstein, which includes all hypersurface and ICIS germs, then the link also has a  Spin structure, canonical up to homotopy (by \cite{Seade-Laufer}). Assume the link is  an integral homology sphere $\Sigma$, and  let $$R(\Sigma) = {\rm Hom}^* \big(\pi_1(\Sigma), {\rm SU}(2)\big) \big/ ({\rm ad \;SU}(2) \big)$$
 be the space of irreducible ${\rm SU}(2)$-representations of its fundamental group modulo conjugation. The space $R(\Sigma)$  is nondegenerate if it satisfies a certain condition for every $\alpha \in R(\Sigma )$ (see \cite{Taubes}).   In this case  $R(\Sigma)$ has finite cardinality 
  and its Casson invariant  is defined via a signed count of its points:
  $$ \qquad  \lambda(\Sigma) = \frac{1}{2} \sum_{\alpha \in R(\Sigma)} \e_{\alpha} \quad , \quad \hbox{with} \; \e_{\alpha} = \pm 1 \;,
  $$
The integers $\e_{\alpha} = \pm 1$  are determined from an intersection theory associated with a Heegaard splitting of  $\Sigma$.
If $R(\Sigma)$  is degenerate,  then it needs to be perturbed first to make it finite and then $\lambda(\Sigma)$ can be defined  similarly, see \cite{Taubes}.
 In ~\cite {FS} Fintushel and Stern  proved  that the Casson invariant of the Brieskorn homology spheres $\Sigma(p,q,r)$ is: 
$$ \lambda(\Sigma(p,q,r)) \, = \, \frac{1}{8} \; \sigma (F(p,q,r))\,,$$
 where $F(p,q,r)$ is the Milnor fiber  (see  \cite {Col-Save1} for a geometric proof). This led to the
{\it Casson invariant conjecture} stated 
in \cite {Ne-Wa}:

\begin{conjecture} If $\Sigma$ is an integral homology sphere which is the link of an isolated complete intersection surface singularity, then its Casson invariant equals $1/8$ the signature of the Milnor fiber:
$$ \lambda(\Sigma) \,=\, \frac{1}{8}\sigma (F) \,.$$
\end{conjecture}

\noindent
This conjecture was proved 
in \cite{Ne-Wa}  for all 
 weighted homogeneous surface singularities and for other families too, including the Brieskorn-Hamm complete intersections.
Yet, the conjecture remains open. 

This can also be 
 regarded from the viewpoint of the Seiberg-Witten monopole equations (see \cite[p. 282, 3.1]{Ne-Ni1}). 
The Seiberg-Witten invariant of a closed oriented 3-manifold $M$ is a function $SW$ from the set 
${\mathcal S}(M)$ of ${\rm Spin}^c$-structures on $M$ to the integers $\Z$. Roughly speaking, this invariant  counts the gauge equivalence classes of solutions to the Seiberg-Witten equations. 
In the case of homology spheres $\Sigma$ there is only one ${\rm Spin}^c$-structure and therefore 
a single Seiberg-Witten invariant, that we denote ${\rm SW}(\Sigma)$.
In \cite{MOY} is proved that for the Brieskorn homology spheres, this coincides with the Casson invariant:
$${\rm SW}(\Sigma(p,q,r)) = \lambda(\Sigma(p,q,r)) \,.$$
It was then conjectured by P. Kronheimer that these two invariants coincide for all 3-dimensional homology spheres. That conjecture is proved in \cite {Lim-1}.
Hence the Casson invariant conjecture can be studied in terms of Seiberg-Witten invariants. This line of research is being done by A. N\'emethi {\it et al}
 (cf. \cite{La-Ne, Ne-3, Ne-Ni1, Ne-Ni2, NeSigu}). They  also study conditions under which the geometric genus is a topological invariant.

 In the previous discussion the surface singularities in question are  isolated, so the link is a smooth 3-manifold and if the singularity is an ICIS, then the link  is isotopic to the boundary of the Milnor fiber. In 
   \cite{Ne-Sz} the authors consider    the boundaries of Milnor fibers of non-isolated surface singularities in $\mathbb{C}^3$. These turn out to be Waldhausen manifolds as well, by   \cite{MPi1, MPi2, Ne-Sz}. In \cite{Ne-Sz} the authors  give a 
way to determine the Waldhausen decomposition of these manifolds, and they study thoroughly their geometry and topology. The fact that the boundaries of the Milnor fibers of non-isolated complex hypersurfaces in $\C^3$ are Waldhausen manifolds is proved also in \cite {FM}, where the authors extend that theorem to singularities defined by functions of the type $f \bar g$ with $f,g$ being holomorphic (see Section \ref{Sec. iso. crit. value} below).

 %\newpage
 
\section*{}

%%%%%%%
%

%%%%%%%%%%%%%
%%%%%%%%%%%
%%%%%%%%%%%%

\vskip2cm

%\centerline{\LARGE \bf Part II: The real analytic case}
\specialsection*{\bf  PART II: BEYOND THE HOLOMORPHIC REALM}
We now look at Milnor fibrations  for real analytic singularities. This  emerged too from Milnor's seminal work in \cite{Milnor:ISH, Mi1}. We also look at meromorphic functions, and remark that the much of the discussion below goes through for semi-algebraic and subanalytic maps (see \cite{Rafaella-Men, Dut1, Dut-Ar-Ch-An}).

\section{\bf Foundations and first steps}\label{foundations} 
 The following was stated as Hypothesis 11.1 in \cite{Mi1}.

\begin{definition} \label{Milnor condition} A real analytic map-germ $f: (\R^n, \0) \to (\R^p,0)$, $n \ge p > 0$ 
 {\it satisfies the Milnor condition at } $\0$ if the derivative $Df(x)$ has rank $p$ at every
point $x \in U\setminus \0$, where $U$ is some open neighborhood of $\0 \in \R^n$.
\end{definition} 

\vv
The following   extends  the fibration
theorem to the real setting:

\begin{theorem}~\label{Milnor-thm-real}{\rm [Milnor]} Let $f $ satisfy the  condition {\rm (\ref{Milnor condition})} at $\0$.
For every $\e > 0$ sufficiently small, let $\delta> 0$ be sufficiently 
small with respect to $\e$ and consider the Milnor tube $N(\e,\delta) := f^{-1}(\partial \D_\delta) \cap \B_\e$, where $\D_\delta$ is the disc in $\R^p$ of radius $\delta$ and center at $0$,  $\partial \D_\delta$ is its boundary and $ \B_\e$ is the closed ball in $\R^n$
of radius $\e$ and center $\0$.
Then 
\begin{equationth}\label{fibration-iso-song-real}  f|_{N(\e,\delta)}: N(\e,\delta)   \longrightarrow\, \partial  \D_\delta \;,\end{equationth}  \hskip-4pt
is a fiber bundle.  Moreover, the tube $N(\e,\delta)$ is diffeomorphic to $\s_\e^{n-1} \setminus \big(f^{-1}({\buildrel {\circ} \over {\D}}_\delta) \cap \s_\e^{n-1}\big)$, where  ${\buildrel {\circ} \over {\D}}_\delta$ is the interior, and {\rm (\ref{fibration-iso-song-real})} determines an equivalent fiber
bundle:
\begin{equationth}\label{Milnor-real}
\varphi \colon \s_\e^{n-1} \setminus L_f \to \s^{p-1}\,,
\end{equationth} \hskip-4pt
where $L_f = f^{-1}(0) \cap \s_\e^{n-1}$ is the link. The projection  $\varphi$ is $f/\Vert f \Vert$ in a tubular neighborhood of  $L_f$. 
 \end{theorem}

The statement that $\varphi = f/\Vert f \Vert$ in a tubular neighborhood of the link $L_f$ is implicit in Milnor's book and it was made explicit in \cite{CSS2,PS2}.
The proof of (\ref{fibration-iso-song-real}) is  an easy extension of the proof of Ehresmann's Fibration Lemma. As in the complex case, 
 one  then constructs an
integrable vector field $v$ in the ball $\bar{\B}_\e$, which is
transverse to all spheres in this ball centered at $\0$, and
transverse to all Milnor tubes. The integral curves of $v$ allow
us to carry $N(\e,\delta) $  diffeomorphically into the complement
of $f^{-1}(\D_\delta)\cap\s_\e^{n-1}$ in the sphere $\s_\e^{n-1}$,
keeping its boundary fixed, and  one extends the induced fibration
to all of $\s_\e^{n-1} \setminus L_f$ using for instance that
the normal bundle of the link is trivial.

%\begin{figure}
%\centering
% \includegraphics[height=7cm ]{fibra-Milnor-Mi-2}
%\caption{Milnor's Figure 4 in \cite[p. 54]{Mi1}: pushing the  fiber from the interior of the ball toward the sphere.}
%\end{figure}

Yet, we cannot in general  inflate the tube in such a way that the projection $\varphi$ is $f/\Vert f \Vert$ everywhere. In fact this theorem has two weaknesses:

\vv

\smallskip

1) It is much too stringent: map-germs satisfying 
condition (\ref{Milnor condition}) are highly non-generic.

\smallskip
2) One has no control over the projection map $\varphi$ outside a neighborhood of the link.

\medskip

Of course that every complex valued holomorphic function with an isolated critical point  satisfies  (\ref{Milnor condition}), and so does if we compose
such a function with a real analytic local diffeomorphism of either the target or the source. The interesting point is  finding examples which are honestly real analytic.
Milnor exhibited the following examples in his book, suggested to him by 
N. Kuiper.
Let $A$ denote either the complex numbers, the quaternions or the Cayley numbers, and define 
\[h : A \times A \to A \times \R  \;,\]
by $h(x,y) \,=\,(2x \bar y, |y|^2 - |x|^2)$. Milnor first proves 
(\cite[ 11.6]{Mi1}) that this mapping carries the unit sphere of $ A \times A$ to the unit
sphere of $ A \times \R$ by  a {\it Hopf fibration}. Then he defines, more generally,
\[f : A^n \times A^n \to A \times \R \;,\] 
by 
\begin{equationth}\label{Milnor-mixed}
f(x,y) \,=\,(2 \langle x, y \rangle, \|y\|^2 - \|x\|^2)\,,
\end{equationth} \hskip-4pt
where $\langle \cdot, \cdot \rangle$ is the Hermitian inner product in $A$. This map is a local submersion on a punctured neighborhood 
of $(0,0) \in A^n \times A^n$. The link  of the corresponding singularity is the Stiefel manifold of 2-frames in $A^n$ and the Milnor fibre  
 is a disc bundle over the unit sphere of $ A^n$.

\v
For $p=1$, condition (\ref{Milnor condition}) is always satisfied (see for instance \cite{Ver}). 
For maps into $\R^2$,  generically the critical values are 
real curves converging to $(0,0)$, though there are several families of singularities satisfying (\ref{Milnor condition}), see for instance Section \ref{Section mixed}.  For $p > 2$, few examples are known of map-germs satisfying  (\ref{Milnor condition}) and having a ``non-trivial Milnor fibration", where non-trivial means that the fibers are not discs.

There are in fact pairs $(n,p)$ as above for
 which no such examples exist, as stated in Theorem \ref{thm. classification of pairs} below, proved 
   by Church and Lamotke in  \cite{CL}, completing previous work  by Looijenga in \cite{Lo1}:

\begin{theorem} \label{thm. classification of pairs}
Let $n, p$ be positive integers. 
\begin{enumerate}
\item  If
$ \;  2 \ge  n - p  \ge  0$, then such examples exist for the pairs $\{(2, 2), (4, 3), (4, 2)\}$.
\item  If $ n - p = 3$, non-trivial examples exist for $(5, 2)$ and $(8, 5)$, and perhaps for  $(6, 3)$.
\item If  $ n - p   \ge 4$, then such examples exist for all $(n, p)$.
\end{enumerate}
\end{theorem}

In particular, if $p=2$, such examples exist for all $n \ge 4$. 
The case $(6,3)$ was left open and it was recently settled affirmatively  in 
 \cite{RHSS}. 

The proof in  \cite{CL} follows the line  in \cite{Lo1} and consists of an inductive process to decide
for which pairs $(n,p)$ there exists a 
 codimension $p$ submanifold $K$ of the sphere $\s^{n-1}$ with a tubular neighborhood $N$ which is a product $N \cong K \times D^p$, such that the natural projection 
 $K \times (D^p \setminus \{0\}) \to \s^{p-1}$ given  by $(x,y) \mapsto y/\Vert y \Vert$ extends to a smooth fiber bundle projection $\s_{n-1} \setminus K \to \s^{p-1}$. No explicit singularities satisfying the Milnor condition (\ref {Milnor condition}) were given.

 The first  explicit non-trivial example of a real analytic singularity with target $\R^2$ 
satisfying the condition (\ref{Milnor condition}), other than those in \cite{Mi1}, was given by A'Campo
\cite {AC1}. This is
the map $ \C^{m+2} \to \C$ defined by
\[(u,v,z_1,...,z_m) \longmapsto uv(\bar u + \bar v) + z_1^2 +...+ z^2_m \,.\]
  
%%%%%%%% IJM

The following notion was introduced in \cite {RSV}:

\begin{definition}\label{def. strong-milnor}
Let 
$ f : (\R^n ,\0) \to (\R^p,0)\,, $ $n > p \ge 2$,  be analytic and  satisfy  the Milnor condition at $\0$, and  let $L_f$ be its link. We say that $f$ satisfies {\it the strong 
Milnor condition at}
 $\0$ if for every sufficiently small sphere $\s_\e$  around $\0$, 
 \[ \frac{f}{|f|}  : \s_{\e} - L_f\,\to
\,\s^{p-1}\, \] is  a fiber bundle. 

\end{definition}

Jacquemard in \cite{Jaq} studied sufficient conditions to insure that  maps into $\R^2$  satisfying the condition (\ref {Milnor condition}) actually satisfy the strong Milnor condition. He gave two conditions which are sufficient  but not necessary, and he constructs several examples of maps satisfying these conditions. The first of Jacquemard's conditions is geometric,   the second
condition is algebraic. These are:

\vv  {\bf Condition  A)}: 
There exists a neighborhood $U$ of the origin in
$\R^n$ and a real number  $0 < \rho < 1$ such that  for all $x \in U-{0}$
one has:
\[\frac{|\langle grad \,f_1(x), grad \, f_2(x) \rangle|} {||grad \,
f_1(x)||\,\cdot \, ||grad \, f_2(x)||} \, \leq \, 1-\rho \, ,\]
where $\langle \cdot,\cdot \rangle$ is the usual inner product in $\R^n$.

\vv  {\bf  Condition   B)} If  $\varepsilon_n$ denotes the local ring of analytic
map-germs at the origin in $\R^n$, then
the  integral closures in $\varepsilon_n$ of the
ideals  generated by the partial derivatives 
\begin{equationth}\label{2.3}
\Big(\frac{\partial
f_1}{\partial x_1}\,,\,\frac{\partial f_1}{\partial x_2}\,,
\cdots,\, \frac{\partial f_1}{\partial x_n}\Big)  \quad \hbox
{and} \quad \,
\Big( \frac{\partial f_2}{\partial x_1}\,,\,\frac{\partial
f_2}{\partial x_2}\,, \cdots,\, \frac{\partial f_2}{\partial
x_n}\Big)  
\end{equationth} \hskip-4pt
coincide, where $f_1, f_2$ are the components of $f$.

One has:

\begin{theorem}\label{Jacquemard} {\rm [Jacquemard]}
 Let $f:(\R^n,\0) \to (\R^2,0)$, $n > 2$, be an analytic map-germ.  If the component
functions $f_1$ and $f_2$ of $f$  satisfy the previous two conditions A and B, then   for  every sphere
$\s_\e^{n-1}$ of radius $\e>0$ sufficiently small and centered at $\0$, one has that the projection map in {\rm (\ref {Milnor-real})},
$$\varphi \colon \s_\e^{n-1} \setminus L_f \to \s^{1}\,,$$
 can be taken to be $f/\Vert f \Vert$ everywhere.
\end{theorem}

\begin{problem} {\rm
What is the equivalent of Theorem \ref{thm. classification of pairs} for the strong Milnor condition? That is, for which pairs $(n,p)$ there exists an analytic map-germ $ f : (\R^n ,\0) \to (\R^p,0)\,, $ $n > p \ge 2$,  satisfying the strong 
Milnor condition at
 $\0$?}
\end{problem}

When $p=2$, such maps exist for all $n \ge 4$. There are also several  examples with $p=3$ in \cite{Jaq}.

It was noted in~\cite {RSV} that the above condition (B) can be relaxed and still have sufficient conditions to guarantee the strong Milnor 
condition. For this we recall the notion of the real integral closure of an ideal as given in~\cite{Ga1}:

\begin{definition}  Let $I$ be an ideal in the ring ${\varepsilon_m}$. {\it The real integral closure}\index{integral closure} of $I$, denoted by
$\overline I_{\R},$ is the set of $h \in \varepsilon_m$ such that for all analytic $\varphi: (\R,0)\rightarrow (\R^m,0)$, we
have $h\circ \varphi\in (\varphi^{\ast} ({ I})){\varepsilon}_{1}$.
\end{definition}

\vv
Given  $f:(\R^n,0) \to (\R^2,0)$ as above, let us set: 

\vv \n {\bf Condition  ${\rm B}_\R$:} The real integral closures of the Jacobian ideals in (\ref{2.3}) coincide.

\vv
 For complex analytic germs both conditions (B) and (${\rm B}_\R$) are equivalent (see~\cite{Te1},~\cite{Ga1}).
As pointed out in~\cite {RSV}, essentially 
the same proof of Jacquemard in~\cite {Jaq} gives:

\vv
\begin{theorem}\label{Jaq-RSV}
{  Let $f:(\R^n,0) \to (\R^2,0)$ be an analytic map-germ
that satisfies the Milnor condition.
 If its components $f_1$, $f_2$ satisfy the condition (A) above and condition (${\rm B}_\R$), then $f$ satisfies the strong Milnor condition.\index{Milnor condition, strong}
}
\end{theorem}

This improvement of \ref{Jacquemard} was used in~\cite {RSV} to prove a stability theorem for 
real singularities with the strong Milnor
 condition. This was also used in~\cite {AR} to find 
a theorem in the vein of (\ref{Jaq-RSV}) but using  ``regularity conditions'' instead of the
Jacquemard's conditions. 
This  inspired \cite {Se5, RSV} and the use of the canonical pencil described in Section \ref {d-regularity}. We refer for instance to \cite{ADD, RHSS, Dut-Ar-Ch-An, DRC, Oka1, Oka2, Oka3, Oka5, PS1} for recent work on the topology of the Milnor fibers.

\vv

D. Massey in \cite{Mas} improved (\ref {Jaq-RSV}) using a different viewpoint,  via a generalized \L ojasiewicz inequality (see the next section). We remark that Massey's viewpoint applies in a larger setting, not requiring $f$ to have an isolated critical point, and it relaxes significantly condition (${\rm B}_\R$).

%%%%%%%%%%%
%%%%%%%%%%%%
%%%%%%%%%%%%% 

\section{\bf On functions with a non-isolated critical point.}\label{Sec. iso. crit. value} 

As noted before, considering map-germs $(\R^n,\0) \to (\R^p,0)$ with an isolated critical point is very stringent. We now discuss the general case of arbitrary critical locus, starting with the slightly more general case of functions with an isolated critical value.

\subsection{Functions with an isolated critical value}
$\qquad$ Every holomorphic   map-germ   \\ $(\C^m,\0) \to (\C,0)$ with a critical point at $\0$ have an isolated critical value, and the fibration theorems \ref {Fibration Thm., version1} and \ref {Fibration Thm., version2} hold in this setting. It is thus natural to look for extensions  of Milnor's theorem \ref {Milnor-thm-real}
for analytic map-germs $f:=(f_1,\cdots,f_p): (\R^n,\0) \to (\R^p,0)\,$ with a possibly non-isolated critical point at $\0$, such that $0 \in \R^p$ is the only critical value, {\it i.e.}, all critical points of $f$ are in the special fiber $V:=f^{-1}(0)$.
This was first done  in \cite{PS2}.

\begin{definition}
 If $f$  admits a fibration in tubes of the type (\ref{fibration-iso-song-real}), then we call this the (local)  {\it Milnor-L\^e fibration} of $f$ (or the Milnor fibration in tubes). If it admits a fibration on the spheres of the type (\ref{def. strong-milnor}), then we call this the  (local) {\it  Milnor fibration} of $f$ (or the Milnor fibration on spheres).
\end{definition}

Given a real analytic map-germ $(\R^n,\0) \buildrel{f} \over {\to} (\R^p,0)$, $n > p \ge 1$, with an isolated critical value at $0 \in \R^p$, we want to know when is there a local Milnor-L\^e fibration. That is, we want conditions to insure that given a   ball $\B_\e$ bounded by a Milnor sphere $\s_\e$  for $f$
         (see Definition \ref{Definition 2.1}), 
there exists  
 a ball $\D_\delta$  of some radius $\delta$ in $\R^p$, centered at $0$, such that if we 
set $N_f(\e,\delta) := \ f^{-1}(\D_\delta \setminus \{0\})  \cap \mathbb B_\e $, then
$$ f|_{N_f(\e,\delta)}: N_f(\e,\delta)  \longrightarrow \Delta_\delta \setminus \{0\}$$
is  a $C^\infty$ fiber bundle.

We know from \cite{Mi1} that this   is always satisfied when $f$ has an isolated critical point  (that  is  immediate from the implicit function theorem). Yet, when the critical  point is not isolated the situation is more delicate. 
In \cite{PS2} ii was noticed that if the map-germ $f$ is such that $V(f) = f^{-1}(0)$ has dimension $>0$ and $f$ has the Thom $a_f$-property,  then $f$  has a local Milnor-L\^e fibration.

The study of Milnor fibrations for real analytic map-germs was also addressed by D. Massey in \cite{Mas}. Recall that 
in \cite{Ha-Le} Hamm and  L\^e
used the complex analytic \L ojasiewicz inequality   to show that Thom stratifications exist. Massey gives the appropriate generalization for the real analytic setting:

\begin{definition}\label{Massey-def}
An analytic germ  $(\R^n,\0) \buildrel {f} \over {\to} (\R^2,0)$ satisfies the strong \L ojasiewicz inequality at $\0$ 
if there exists a neighborhood $\mathcal W$ of $\0$ and constants $c, \theta \in \R$ with $c >0$, $0 < \theta < 1$, such that for all $x \in \mathcal W$ one has:
$$ \Vert f(x) \Vert^{\theta} \le \, c \, \min_{|(a,b)|=1}  \vert a \nabla g(x) + b \nabla h(x)    \vert \;.
$$
In this case the germ  $f$ is said to be {\it \L -analytic}. 
\end{definition}

The main theorem in \cite{Mas} says:

\begin{theorem}[Massey] If $f$ is \L -analytic, then for every Milnor sphere $\s_\e$ there is a Milnor tube $N_f(\e,\delta)$ where $f$ is a proper stratified submersion and the projection of a $C^\infty$ fiber bundle. That is, \L -analytic maps have  Milnor-L\^e fibrations.
\end{theorem}
  
 Now we need the following definition from \cite{CisMSS}:

\begin{definition}\label{transversality property}
Let $f: (\R^n,\0) \buildrel{f} \over {\to} (\R^p,0)$, $n > p \ge 1$, have an isolated critical value at $0 \in \R^p$. We say that 
$f$ has {\it the transversality property} if   for every sufficiently small sphere  $\s_{\e}$ around the origin in $\R^n$, there exists $\delta > 0$ such that all the fibers 
$f^{-1}(t) $ with $|| t || \le \delta$ 
meet transversally the sphere  $\s_{\e}$. 
\end{definition}

The transversality property appears already in \cite{Ha-Le}, and  in \cite{Oka5} this is called the Hamm-L\^e condition. Maps with the Thom $a_f$-property and non-empty link have the transversality condition, but not conversely:
there are examples by M. Oka in \cite{Oka4}  of maps with the transversality property which do not have the Thom $a_f$-property (see also \cite{MS2}).   

The theorem below  is Theorem 3.4 in \cite{CSS3}. This improves  \cite[Theorem 1.3]{PS2}. 

\begin{theorem}\label{prop.fibr.tube}
Let $f: (\R^n,\0)   {\to} (\R^p,0)$, $n > p \ge 1$, have an isolated critical value at $0 \in \R^p$. Assume further that $f$ has the transversality property and $V(f) := f^{-1} (0)$ has dimension greater than 0.  Then $f$ has a local Milnor-L\^e fibration:
$$   f|_{N_f(\e,\delta)}: N_f(\e,\delta)  \longrightarrow \Delta_\delta \setminus \{0\} \,,$$  
with $N_f(\e,\delta) := f^{-1}(\D_\delta \setminus \{0\}) \cap \B_\e$ for some ball $\D_\delta \subset \R^p$, $0 < \delta  \ll \e$.
This determines an equivalent fiber bundle:
 $$\varphi: (\s_\e \setminus L_f) \longrightarrow \s^{p-1} \,,$$
 where the projection map $\varphi$ is $f/\Vert f \Vert$ restricted to $[\s_\e \cap N(\eta,\delta)]$.
 \end{theorem}

\medskip  

The way to pass from the fibrations in tubes to that on the sphere is as before: one constructs a smooth vector field $\zeta$ in the ball $\B_\e$ minus $V$, satisfying:
$$
\bullet  \quad \hbox{Each integral line is transversal to all spheres in } \B_\e \, \hbox{centered at} \, \0 \;.\qquad \qquad 
\qquad \qquad \qquad \qquad 
$$
$$
\bullet  \quad \hbox{Each integral line is transversal to all tubes}  \, f^{-1}(\partial  \D_\delta) \; \hbox{contained in}  \, \B_\e \,,\qquad \qquad \qquad \; \; \,
\qquad \quad  
$$

The difference with the holomorphic setting is that we cannot guarantee now a third condition: that the vectors $f(z)$ are collinear for all points in each integral line (cf. Theorem \ref{Th d-regular} below). We discuss this in the following section.

For maps of the type $f \bar g$ that we envisage below, 
 there is simple criterium in \cite{FM} and \cite[Proposition 3.5]{MS2} to decide whether or not the map has the transversality property. This is called  CT-regularity in \cite{MS2}. The advantage  of this criterium is that it is easy to use in practice.

\subsection{Functions with arbitrary discriminant}

We now consider the general setting and study Milnor fibrations for map-germs $(\R^n,\0)   {\to} (\R^p,0)$ with arbitrary critical points, following \cite{CGS, CisMSS}. We start with an example studied in \cite{LopezdeMedrano:SHQM} by L\'opez de Medrano. 

\begin{example}\label{ex. Santiago}
Consider   
maps $(f,g): \R^n \to \R^2$  the form:
\[
 (f,g) = \left( \sum_{i=1}^n a_i x_i^2 \, , \, \sum_{i=1}^n b_i x_i^2 \right) \, ,
\]
where the $a_i$ and $b_i$ are real constants in generic position in the Poincar\'e domain. This means that the origin is in the convex hull of the points $\lambda_i := (a_i,b_i)$ and 
no two  $\lambda_i$ are linearly dependent.

A simple calculation shows that $(f,g)$  is a complete intersection and the corresponding link is a smooth non-empty manifold of real codimension 2 in the sphere. 
The  critical points $\Sigma$  of $(f,g)$ are the coordinate axis of $\R^n$ and the set  $\Delta(f,g)$ of critical values is the union of the $n$ line-segments in 
$\R^2$ joining the origin to the points $\lambda_i$. Hence $\R^p$ splits into various connected components, and it is proved in \cite {LopezdeMedrano:SHQM} that the topology of the fibers over points in different components changes. Yet, we know from \cite{CisMSS} that these map-germs have the transversality property (\ref{transversality property}), and away from the critical set  they have a Milnor-L\^e fibration. 
In fact these maps are $d$-regular too, a concept that we discuss in the next section and implies that  away from the discriminant, they have also a Milnor fibration on small spheres  with projection map $(f,g)/\Vert (f,g) \Vert$. 
\end{example}

More generally, 
consider   now an open neighbourhood $U$ of $\0\in\R^n$ and 
 a $C^\ell$ map  $f\colon (U,\0)\to (\R^p,0)$, $n > p \ge 2$, $\ell \ge 1$,  with a critical point at $\0$.
 Denote by $\mathcal{C}_{f}$ the set of critical points of $f$ in $\B_\e$ and let $\Delta_f$ be the image
$f(\Delta_f)$. These are the critical values of $f$; we call $\Delta_f$  the discriminant of $f$.

\begin{definition} \label{defi_tp}
We say that the map-germ $f$ has the \textit{transversality property} at $0$ if there exists a real number $\e_0>0$ such that, for every $\e$ with $0<\e \leq \e_0$, there exists a real number $\delta$, with $0<\delta \ll \e$, such that for every $t \in \B_\delta^k \setminus \Delta_f$ one has that either $f^{-1}(t)$ does not intersect the sphere $\S_\e^{n-1}$ or $f^{-1}(t)$ intersects 
$\S_\e^{n-1}$ transversally in $\R^n$.
\end{definition}

The transversality condition of the fibers with small spheres ensures having a Milnor-L\^e fibration, even for $C^\ell$ maps with non isolated 
critical values. Of course that as in  Example \ref{ex. Santiago}, if the base of the fibration has several connected components (sectors), then the topology of the fibers can change from one sector to another. We have the following result from \cite{CisMSS}.

\begin{proposition}\label{prop:MLF}
 Let $f\colon (\R^n,\0) \to (\R^p,0)$, $n\geq p\geq 2$, be a map-germ of class $C^\ell$ with $\ell\geq 1$. If $f$ has the transversality property, then the restrictions:
 %\begin{align}
 
 $\bullet$ $f|\colon \B_\e^n \cap f^{-1}(\mathring{\D}_\delta^p \setminus \Delta_f)  \longrightarrow  (\mathring{\D}_\delta^p \setminus \Delta_f) \cap \hbox{Im}(f).$
%\intertext{and}

$\bullet$  $f|\colon \B_\e^n \cap f^{-1}(\S_\delta^{p-1} \setminus \Delta_f) \longrightarrow (\S_\delta^{p-1} \setminus \Delta_f) \cap \hbox{Im}(f)$

\noindent
are  $C^\ell$ fiber bundles, where $\e$ and $\delta$ are small enough positive real numbers,  $\B_\e^n \subset \R^n$ and  $\D_\delta^p \subset \R^p$ are the closed balls of radius $\e$ and  $\delta$ centered at $\0$ and $0$, respectively, $\mathring{\D}_\delta^p$ is the interior of the ball $\D_\delta^p$ and $\S_\delta^{p-1} $ is its boundary. If $f$ is analytic, then the fibrations above are $C^\infty$.
\end{proposition}

%%%%%%%%%%%%

In \cite{CGS} the authors study
 the topology of the fibres of  real analytic maps $\R^n \to \R^{p+k}$, $n > p+k $,
 inspired by the classical L\^e-Greuel formula for the
Milnor number of isolated, complex, complete intersection germs. The idea is that if the map germ
is defined by analytic functions $(f_1,...,f_p, g)$, then we can study the topology of its fibres  by
comparing it with the topology of the germ we get by dropping down $g$.
We  require for this that the map $f:= (f_1,...,f_p)$ actually satisfies the Thom $a_f$-property with
respect to some Whitney stratification $\{S_\alpha\}$, and  that its zero-set $V(f)$ has dimension
$\ge 2$ and it is union of strata. The map-germ $(\R^n,\0)  \buildrel{g}\over{\to} (\R^k,0)$ is assumed to have an
isolated critical point in $\R^n$ with respect to the stratification $\{S_\alpha\}$. By
Proposition~\ref{prop:MLF} the map-germs $f$ and $(f,g)$ have associated local  Milnor-L\^e  fibrations.
Then one has  the corresponding L\^e-Greuel formula \cite[Theorem 1]{CGS}:

\begin{theorem}\label{Thm:LG}
Let $F_f$ and $F_{f,g}$ be  Milnor fibres of $f$ and $(f,g) $ (any Milnor fibres, regardless of the  fact that the topology of the fibers may depend on the connected component of the base once we remove the discriminant). Then
one has:
\begin{equation*}
 \chi(F_f) = \chi(F_{f,g}) + {\rm Ind}_{\mathrm{PH}} {\nabla} \tilde{g}|_{F_f\cap\B_{\e'}} \,,
\end{equation*}
where $\tilde{g}\colon\R^n\to\R$ is given by $\tilde{g}(x)=\Vert{g(x)-t_0}\Vert^2$ with $t_0\in\R^k$ such
that $F_{f,g}=g|_{F_f}^{-1}(t_0)$ and $\B_{\e'}$ is a small ball in $\R^n$ centered at the origin.

The term ${\rm Ind}_{\mathrm{PH}}  {\nabla} \tilde{g}|_{F_f\cap\B_{\e'}}$ on the right, which by
definition is the
total Poincar\'e-Hopf index in $F_f$ of the
vector field ${\nabla} \tilde{g}|_{F_f}$, can be expressed also in the following equivalent ways:
\begin{enumerate}
\item As the Euler class of the tangent bundle of $F_f$ relative to the vector field  ${\nabla}
\tilde{g}|_{F_f\cap\B_{\e'}}$ on its boundary;
\item As a sum of polar multiplicities relative to $\tilde{g}$ on $F_f\cap\B_{\e'}$.
\item As the index of the gradient vector field of a map $\tilde{g}$ on $F_f$
associated to $g$; 
\item As the number of critical points of $\tilde{g}$ on $F_f$;  
\item When $p=1 = k$, this invariant can also be expressed algebraically, as the
signature of a certain bilinear form that originates from \cite{Arnold, EL, EG1, GM1, GM2}.
\end{enumerate}
\end{theorem}

When $n = 2m, p=2q$, $k=2$ and  $\C^m \buildrel {(f,g)} \over{\longrightarrow} \C^{q+1}$  is holomorphic, this is a reformulation of the classical L\^e-Greuel formula (\ref {Le-Greuel}).

We remark that when $k=1$ and the germs of $f$ and $(f,g)$ are both isolated complete intersection germs, there is a  L\^e-Greuel type formula  in \cite{Du-Gru} expressed in terms of   normal data of $f$ with respect to an appropriate Whitney stratification. See also \cite{MS1} for refinements of the above discussion.

%%%%%%%%%%%
%%%%%%%%%%
%%%%%%
%%%%%%%%%%%
%%%%%%%%%%

\section{\bf $d$-regularity and Milnor fibrations} \label{d-regularity}   

The concept of d-regularity introduced in \cite{CSS2} is inspired by  \cite{AR,Bekka,RSV, Se5} and it   is a key for understanding  the difference between  real and  complex singularities concerning Milnor fibrations.

\subsection{The case of an isolated critical value}
Let $U$ be an open neighborhood of the origin $\0 \in \R^n$, and
consider a real analytic germ  $f\colon (U,\0) \to (\R^p,0)$
which  is a
submersion for each $x \notin  V:= f^{-1}(0)$  and has a  critical point at $\0$. 

\begin{definition} 
The   {\it canonical pencil} of $f$ is a family $\{X_\ell\}$ of real analytic spaces  
parameterised by $\R \mathbb P^{p-1}$, defined as follows:  for each $\ell \in
\R \mathbb P^{p-1}$, consider the line $\mathcal L_{\ell} \subset \R^p$
that determines $\ell$, and set
$$X_{\ell} = \{x \in U \, \vert \, f(x) \in \cL_\ell\}\,.$$
\end{definition}
Note that  every two distinct elements of the pencil $X_{\mathcal L}$ and $X_{\mathcal L'}$ satisfy
\begin{equation*}
 X_{\mathcal L} \cap X_{{\mathcal L'}} = V \,.
\end{equation*}
Each $X_{\mathcal L}$ has dimension
$n-p+1$, is non-singular outside
$V$ and their union covers all of $U$.

Each line ${\mathcal L}$ intersects the sphere $\s^{p-1}$ in two
antipodal points $\theta ^-$ and $\theta ^+$. We 
decompose the line ${\mathcal L}$  into the corresponding  half lines
accordingly:
\begin{equation*}
 {\mathcal L} = {\mathcal L}^- \cup \{ 0 \} \cup {\mathcal L}^+.
\end{equation*}

If we define $E_{{\theta}^\mp}$ to
be the inverse image $f^{-1}({\mathcal L}^\mp)$, respectively, then we can express each element of the
canonical pencil as the following union:
\begin{equationth}\label{XL=union-Eteta-link}
X_{\mathcal L} = E_{{\theta}^-} \cup V \cup E_{{\theta}^+}.
\end{equationth} \hskip-4pt
If $L_V$ is the link of $f$, we can describe the fibers of the map $\varphi = {f}/{\Vert{f} \Vert}:
\s_\e^{n-1} \setminus  L_V \rightarrow \s^{p-1}$ as :
\begin{equation*}
 \varphi^{-1}(\theta ^-) = E_{{\theta}^-} \cap \s^{n-1}_\e \qquad , \qquad 
 \varphi^{-1}(\theta ^+) = E_{{\theta}^+} \cap \s^{n-1}_\e \;.
\end{equation*}
Then
we can write:
\begin{equation*}
 X_{\mathcal L} \cap \s^{n-1}_\e = (E_{\theta^-}\cap \s_\e^{n-1}) \cup L_V \cup (E_{\theta
^+}\cap \s_\e^{n-1}) = \varphi^{-1}(\theta^-) \cup L_V \cup
\varphi^{-1}(\theta^+).
\end{equation*}

We now assume that $f: U \rightarrow \R^p$ is real analytic, with an isolated critical value at $\0$ and it is {\it locally surjective}, {\it i.e.} the restriction of $f$ to every
neighborhood of   $\0 \in U$ covers a neighborhood of $0 \in \R^p$.

\begin{definition}
We say that $f$ {\it  is
$d$-regular}  if there exists $\e_0
> 0$ such that for every $\e \leq \e_0$ and for every line ${\mathcal L}$ through the origin in
$\R^p$, the sphere $\s^{n-1}_\e$ and the
manifold  $X_{\mathcal L}\setminus V$ are transverse.
\end{definition}

As examples of d-regular maps one has:
\begin{itemize}

\item All  holomorphic maps $\C^n \to \C$, all polar weighted homogeneous polynomials and real weighted homogeneous maps with an isolated critical value, are d-regular maps. 

\item If $f$ and $g$ are holomorphic maps $\C^2 \to \C$ such that the product $f \bar g$ has an isolated critical value at the origin, then the map $f \bar g$ is d-regular, by \cite{PS2}.

\item The strongly non-degenerate mixed functions in \cite{Oka2} are all d-regular, by  \cite{Oka2, CSS2}. 

\item Direct sums of $d$-regular maps. That is, if $f$ is $d$-regular in the variables $(x_1,\cdots,x_n)$ and $g$ is $d$-regular in the variables $(y_1,\cdots,y_m)$, then $f + g$ is $d$-regular in the variables $(x_1,\cdots,x_n, y_1,\cdots,y_m)$, by  \cite{CSS2}.
\end{itemize}

%%%%%%%%%%%%

The following is a fundamental property of d-regularity. We refer to \cite{CSS2} for its proof.

\begin{theorem}\label{Th d-regular}
The real analytic map  $f$ is $d$-regular if and only if  there exists a smooth vector
field $\zeta$ such that its integral
lines are transverse to all spheres around
$\0$, transverse to all Milnor  tubes $f^{-1}(\partial {\D}_\eta)
\cap \B_\e \,,$ and tangent to each element $X_{\mathcal L}$ of
the canonical pencil.
\end{theorem}

Such a vector field allows us to inflate the tube and get a fibration on the sphere minus the link, granting that the projection map is $f/\Vert f \Vert$. Hence we get a slight refinement of \cite[Theorem 1]{CSS2}:

\begin{theorem}\label{equivalence of fibrations}
Let $f := (f_1,...,f_p) \colon(\R^n,\0) \to (\R^{p},0)$ be a locally surjective real analytic map with an
isolated critical value at $0 \in \R^p$ and assume $V=f^{-1}(0)$ has dimension $>0$. Then $f$ admits a Milnor-L\^e fibration if and only if it has the transversality property. If this is so, then $f$ is $d$-regular at $0$ if and only if 
one has a commutative diagram
of smooth fibre bundles:
\begin{equation*}
\xymatrix{
\s_\e^{n-1} \setminus L_V \ar[r]^-{\phi}\ar[rd]_{\psi}  &\s^{p-1}\ar[d]^{\pi}\\
& \R \mathbb P^{p-1}
}
\end{equation*}
where $L_V$ is the link, $\psi= (f_1(x): \cdots :f_p(x))$ and
$\phi=\frac{f}{\Vert{f}\Vert}\colon\s_\e^{n-1} \setminus K_\eta\to
\s^{p-1}$ is the    Milnor fibration. Furthermore, 
if  the two fibrations exist (one
on a Milnor tube, another on the sphere minus the link), then
these fibrations are smoothly equivalent. That is, there exists a
diffeomorphism between their corresponding total spaces, carrying
fibers into fibers.
\end{theorem}

This answers affirmatively a question  raised by R. Araujo dos Santos in \cite{Ra2}, where the author proved it for $p=2$ and $f$ weighted homogeneous. The proof in \cite{CSS2} of the equivalence of the two fibrations has a small gap that has been filled in \cite{CisMSS} where the theorem is proved in the more general setting of real analytic maps with arbitrary linear discriminant.
In \cite{CSS2, CSS3} there are other criteria to determine $d$-regularity which can be useful in practice.

The following corollary is an immediate consequence of the previous theorem:

\begin{corollary}\label{Cor:Mil.fib}
Given $f\colon (\R^n,\0) \to (\R^p,0)$ as in the  theorem above,
consider its   Milnor fibration
\begin{equation*}
\phi = \frac{f}{\Vert{f}\Vert}\colon \s_\eta^{n-1} \setminus K_\eta \to \s^{p-1}.
\end{equation*}
Then the union of the link $K_\eta$ and each pair of fibres of
$\phi$ over antipodal points of $\s^{p-1}$
corresponding to the line $\cL_\ell$, %is naturally glued together along
is the link of the real analytic variety $X_\t$.
\end{corollary}

For instance,  if $f:(\C^n,0) \to (\C,0)$ is
holomorphic and it has an isolated critical point at $0$, then
$\{Re \, f = 0\}$ is a real hypersurface and its link is the
double of the Milnor fiber of $f$ with the link $L_f$ as an equator. If $n=2$, then
the link of $Re\,(f) = 0$ is a compact Riemann surface of genus
$2g_f + r -1$ where $g_f$ is the genus of the Milnor fibre of $f$
and $r$ the number of connected components of the link of $f$.
Thus for instance, we know from \cite{CSS1} that for the map $(z_1,z_2)
\buildrel{f} \over {\mapsto} z_1^2 + z_2^q$ one gets that the link
of $Re \,f$ is a closed oriented surface  in the 3-sphere, union
of the Milnor fibres over the points $\pm i$; an easy computation
shows that it has genus $q-1$. This provides an explicit way to determine 
closed surfaces of all genera $\ge 1$ in the 3-sphere by a single analytic equation. 

It would be interesting to study the
geometry and topology of the 4-manifolds one gets
in this way, by considering the link of the  hypersurface in
$\C^3$ defined by the real part of a holomorphic function with an
isolated critical point. For example, for the map $(z_1,z_2,z_3)
\buildrel{f} \over {\mapsto} z_1^2 + z_2^3 + z_3^5$, the
corresponding 4-manifold is the double of the  $E_8$
manifold, whose  boundary is Poincar\'e's homology 3-sphere.

\subsection{The general case}

In \cite{CisMSS} the authors continue the work begun in \cite{CGS} and extend the above discussion on $d$-regularity  to differentiable functions
 $ (\R^n,0) \buildrel {f} \over {\to} (\R^p,0)$  of class $\C^{\ell}$, $\ell \ge 1$, $n\geq p\geq2$, with a critical point at  $\0 \in \R^n$, arbitrary critical values  $\Delta_f$ and non-empty link $L_V$.  This is immediate when the discriminant $\Delta_f$ 
 is linear:

\begin{definition}
The map-germ $f\colon (\R^n,0) \to (\R^p,0)$ has \textit{linear discriminant} 
%if the germ $(\Delta_f,0)$ of the discriminant of the germ $f$ at $0 \in \R^p$ is linear. That is, 
if for some representative $f$ there exists $\eta = \eta(f) >0$ such that the intersection of $\Delta_f$ with the closed ball $\D^p_\eta$ is a union of line-segments, {\it i.e.}:
\[
 \Delta_f \cap \D_\eta^p = {\rm Cone} \big( \Delta_f \cap \S_\eta^{p-1} \big) \, .
\]
We call $\eta$ a \textit{linearity radius} for $\Delta_f$. (The case when $f$ has $0\in\R^p$ as isolated critical value is considered to have linear discriminant with arbitrary linearity radius.)
 \end{definition}

Let $f\colon (\R^n,0) \to (\R^p,0)$ be as above, with linear discriminant, and consider a representative $f$ with linearity radius $\eta>0$. Set 
$\partial \Delta_\eta := \Delta_f \cap \S_\eta^{p-1} \, .
$
For each point $\theta \in \S_\eta^{p-1}$, let ${\mathcal L}_\theta \subset \R^p$ be the open segment of line that starts at the origin and ends at the point $\theta$ (but not containing these two points). 
Set $
 E_{\theta} := f^{-1}({\mathcal L}_\theta)  \, ,$ so
 each $E_{\theta}$ is a manifold of class $C^\ell$ for every $\theta$ in $\S_\eta^{p-1} \setminus \partial \Delta_\eta $.

\begin{definition}
Let $f\colon (\R^n,0) \to (\R^p,0)$ be a map-germ of class $C^\ell$ with $\ell\geq 1$ and linear discriminant. We say that $f$ is \textit{$d$-regular} at $0$ if for some representative $f$ there exists 
$\e_0>0$ small enough such that $f(\B_{\e_0}^n)\subset \mathring{\D}_{\eta}^p$, where $\eta$ is a linearity radius for $\Delta_f$, and such that every $E_{\theta}$ intersects the sphere $\S_\e^{n-1}$ 
transversally in $\R^n$, for every $\e$ with $0<\e \leq \e_0$ and for all $\theta \in \S_\eta^{p-1} \setminus \partial \Delta_\eta $ such that the intersection is not empty. 
\end{definition}

\begin{example} \label{ex. Santiago-gral}
Let $\K$ be either $\R$ or $\C$. 
Let $(f,g): \K^n \to \K^2$ be a $\K$-analytic map of the form:
\[
 (f,g) = \left( \sum_{i=1}^n a_i x_i^q \, , \, \sum_{i=1}^n b_i x_i^q \right) \, ,
\]
where $(a_i, b_i) \in \K$ are constants in generic position as in Example \ref{ex. Santiago}, and $q \geq 2$ is an integer.
By \cite{CisMSS} the  discriminant $\Delta$ is linear and $(f,g)$ is $d$-regular.
 \end{example}
 
It is proved in  \cite{CisMSS} that the fibration theorems  \ref{Th d-regular}   and \ref{equivalence of fibrations} extend to this general setting of $C^\ell$ maps with linear discriminant which have the transversality property and are $d$-regular. Also, there are  in \cite{CisMSS}  examples of non-analytic maps for which these fibration theorems apply.

Consider now the following  example that generalizes \ref{ex. Santiago-gral}.
Let $(f,g)\colon \R^n \to \R^2$ be:
\[
 (f,g) = \left( \sum_{i=1}^n a_i x_i^p \, , \, \sum_{i=1}^n b_i x_i^q \right) \, ,
\]
with $p,q \geq 2$ integers and the $(a_i,b_i)$ as above. If $p \ne q$, the discriminant $\Delta_{(f,g)}$ is not linear. Yet we can always linearize it with a homeomorphism $h$ in $\R^2$. Moreover, these maps have the transversality property and they are $d_h$-regular in an appropriate sense that depends on the homeomorphism $h$. The fibration theorems in \cite{CisMSS} extend to this setting, and in fact to all 
$C^\ell$-maps  that admit an appropriate ``conic modification", a condition that seems to be always satisfied.

 %%%%%%%%%%%
 %%%%%%%%%%%
 %%%%%%%%%%%%
 
 \section{\bf Mixed Singularities} \label{Section mixed}
 
 A mixed function is a real analytic function $\C^n \to \C$ in the complex variables $z=(z_1,É,z_n)$ and their conjugates $\bar z=(\bar z_1,É,\bar z_n)$. This type of functions  appeared in singularity theory already in \cite[Chapter 11]{Mi1} as well as in the work of N. A'Campo  \cite{AC1} and  Lee Rudolph \cite{Ru}. The modern study of mixed functions in singularity theory springs from   \cite{RSV, Seade-librosing}.
 The term ``mixed singularity'' was coined by M. Oka in \cite {Oka2}. 
 
 %%%%%%%%%%
 %%%%%%%%

 \subsection{Twisted Pham-Brieskorn singularities}\label{Twisted Ph-Br}

The paradigm of real analytic functions with target $\R^2$ and satisfying the strong Milnor condition (\ref{def. strong-milnor})
are the {\it Twisted Pham-Brieskorn singularities} \cite [Chapter VII]{Seade-librosing}:
\begin{equation*}%\label{d:twisted Pham-Brieskorn}
(z_1,...,z_n) \buildrel {f} \over {\mapsto} z_1^{a_1} \bar z_{\sigma(1)} + ... +
z_n^{a_n} \bar z_{\sigma(n)} \; , \; a_i \ge 2\,,
\end{equation*}
where $\sigma$ is a permutation
of the  set $\{1,...,n\}$. 
It was noticed in \cite {RSV, Se5} that there exists a smooth action of  $\s^1 \times \R^+$  on $\C^n$ of the form
$$(\lambda,r)\cdot (z_1,\dots,z_n)=(\lambda^{d_1} r^{p_1}
z_1,\dots,\lambda^{d_n} r^{p_n} z_n) \;, \quad \lambda \in
\s^1\,,\, r \in \R^+ \;,
$$
where the $d_j,p_j$ are positive integers such that $
\gcd(d_1,\dots,d_n) = 1 = \gcd(p_1,\dots,p_n)$, and one has:
\begin{equation*}
 f ( (\lambda,r)\cdot (z_1,\dots,z_n)) \,=\, \lambda^d \, r^p \, f(z_1,\dots,z_n)\;,
\end{equation*}
for some positive integers $d, p$. So these are reminiscent of weighted homogeneous complex polynomials. 

The  basic properties are :

\begin{lemma}\label{circle-action}
The above  $\s^1$-action  on $\C^n$ has $V:= f^{-1}(0)$ as an invariant set  and it  permutes the elements 
$X_{\theta}$ of the canonical pencil defined in Section \ref{d-regularity} above. \end{lemma}

\begin{lemma}\label{R-action}
The  orbits of the $\R^+$-action    are transversal to every sphere around $0 \in \C^n$, except the orbit of
$\0 \in \C^n$ which is a fixed point, and all the orbits converge to $\0$ 
when the
time tends to $-\infty$. This action leaves invariant $V$ and each $X_{\theta}$.
\end{lemma}

 The explicit weights for the above mentioned action of $\s^1 \times \R^+$ are easily
computable from the exponents $a_i$ and the permutation $\sigma$ (see lemmas  4.3 and 4.4 in Chapter VII of \cite{Seade-librosing}).

%%%%%

The name {\it twisted Pham-Brieskorn singularities}
comes from the similarity that these  have with the classical Pham-Brieskorn
polynomials (cf. Section \ref{sec. Brieskorn}) and the fact, proved in \cite{Oka3, RSV}, that if the twisting 
$\sigma$ is the identity, then the corresponding open-books are  equivalent
to those of the usual Pham-Brieskorn singularities.

It was proved in \cite {Se5, RSV} that these functions have an isolated critical point at $\0$, so they have a Milnor-L\^e fibration in a  tube, and in fact the local triviality of this bundle follows easily from the  $\s^1$-action in Lemma \ref{circle-action}. The first return map of this action is the monodromy map of the fibration. Then the  flow given by  Lemma \ref{R-action} carries this fibration into one in the sphere with projection map $f/|f|$.

\medskip
%%%%%%%
 \subsection{Polar weighted  and radial weighted singularities}\label{Sec-Polar}

Consider, more generally,  non-constant polynomial maps $f: \R^{2n} \to \R^2$, $n > 1
$, that  carry the origin $\0 \in \R^{2n}$ into the origin $0 \in
\R^2$.  We identify the plane $\R^2$ with the complex line $\C$
and equip $\C^*$ with polar coordinates $\{r \, e^{i \t} \, | \, r
> 0 \,; \, \t \in [0, 2 \pi)\}$.  We also identify $\R^{2n}$
with $\C^n$ in the usual way.
The following concept is introduced in \cite{Cis}.

\begin{definition}\label{d:polar weighted}
The map $f$ is a
 \emph{polar weighted homogeneous polynomial}
if  there exists an action of $\s^1 \times \R^+$  on $\C^n$ of the form
\begin{equation*}\label{eq:polar.action}
(\lambda,r)\cdot (z_1,\dots,z_n)=(\lambda^{d_1} r^{p_1}
z_1,\dots,\lambda^{d_n} r^{p_n} z_n) \;, \quad \lambda \in
\s^1\,,\, r \in \R^+ \;,
\end{equation*}
where the $d_j,p_j$ are positive integers such that $
\gcd(d_1,\dots,d_n) = 1 = \gcd(p_1,\dots,p_n)$, and one has:
\begin{equation*}
 f ( (\lambda,r)\cdot (z_1,\dots,z_n)) \,=\, \lambda^d \, r^p \, f(z_1,\dots,z_n)\;,
\end{equation*}
for some positive integers $d, p$. The $\s^1$-action is called {\it a polar action} while that of $ \R^+$ is a {\it radial action}.
\end{definition}

This includes the twisted Pham-Brieskorn singularities as well as the weighted
homogeneous complex polynomials. There are many other families  (see \cite{Oka1, Cis-Rom}). We remark that the polar action implies that the critical value at $0 \in \C$, unless $f$ is constant.

Assuming that the dimension of $V = f^{-1}(0)$ is more than
$0$, the above Properties 1-4 in Section \ref{sec. Brieskorn} continue to hold in this setting.
One  has the canonical pencil as in  Section \ref{d-regularity}. Its elements $X_\t$ are real algebraic hypersurfaces,
 smooth away from $V:= f^{-1}(0)$, they fill out
the whole ambient space and meet exactly at $V$. The $\R^+$ action
leaves invariant every element of the pencil, which therefore is
transverse to all spheres around $\0$. And the orbits of the
$\s^1$-action are tangent to all the spheres around $\0$ and
permute the hypersurfaces $X_\t$.  Therefore one has a global
fibration as in equation (\ref {global fibration}), which
restricts to the Milnor fibration  (\ref{Fibration Thm., version1}) on
each sphere, and by Property 2 this is equivalent to the fibration
on a tube as in  (\ref{Fibration Thm., version2}). The monodromy  is the first return map of the $\s^1$-action.

In \cite {Cis-Rom} the authors give a classification  of the mixed homogeneous polynomials in three variables which are  polar weighted   with an isolated critical point. This  generalizes classical work of Orlik and Wagreich on complex weighted homogeneous polynomials. In \cite{Inba1, Inba2, Inba3} there are interesting relations with contact structures and the enhaced Milnor number. And in \cite{Bode}, a special class of mixed singularities is used to show that  certain braids can be compactified to become fibered real algebraic knots.

%%%%%%%%%%%%%%%%

%%%%%%%%%%%%%%%%

\subsection{Meromorphic germs and the case $f \bar g$.} 

Given two holomorphic map-germs $\C^n \buildrel {f,g} \over \longrightarrow \C$, we can associate to them: 

$\bullet$ the meromorphic germ $f/g$; 

$\bullet$  the real analytic map $f \bar g$ where $\bar g$ denotes complex conjugation. 

\noindent
Notice that the zero locus $V(f \bar g)$ of $f \bar g$ equals, as a set, the zero locus of $fg$ and consists of $V(f) \cup V(g)$, the union of the zero loci of $f$ and $g$. Away from $V(f \bar g)$  one has:
\begin{equationth}\label{mero-vs-bar}
\frac{f/g}{|f/g|} \, =\, \frac{f \bar g}{|f \bar g|}\;. \
\end{equationth} \hskip-4pt
Therefore the two cases are equivalent when we look at Milnor fibrations with this map as projection;  but they are very different when we look at the Milnor tubes. We discuss  first  the  case $f \bar g$.

\begin{figure}\
\centering
\includegraphics[height=8cm ]{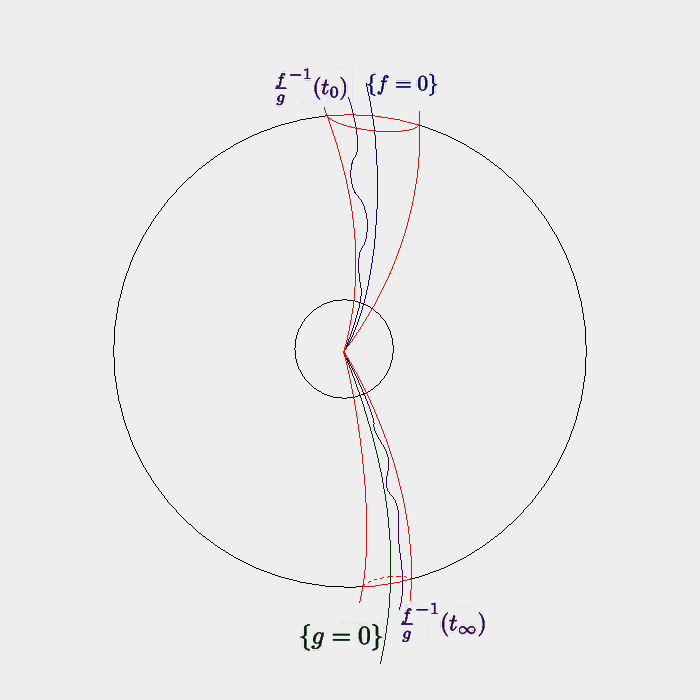}
\caption{A Milnor tube $f/g = |t_0|$ for a meromorphic map.}
\end{figure}

%%%%%%%%

The singular set of $V(f \bar g)$ contains all points in $V(f) \cap V(g)$, which necessarily  are critical points of $f \bar g$. Hence if $n >  2$,  this type of functions must have non-isolated  critical points. Yet, these may have an isolated critical value as in the following examples from \cite{PS2}. 
Define $h: \R^{2n} \to  \R^2$, $n > 1$, by $h = f \bar g$ where 
$\,g (z_1,\cdots,z_n) = z_1 \cdots  z_n\,$ and 
$f$ is the Pham-Brieskorn polynomial 
$$f(z_1,\cdots,z_n) = (z_1^{a_1} + \cdots z_n^{a_n})\,\,\,,\, a_i \geq 2. $$
 A straightforward computation shows that
  $h$ has an isolated critical value at $0 \in  \R^2$ if and only if the 
$a_i$ satisfy 
$ \, \sum_{i=1}^n {1}/{a_i } \, \ne \, 1\,.$
 For $n =2$ this means that at least one $a_i$ is more than 2. For n = 3 the condition is that  
the unordered triple $(p,q,r)$ is not $(2,3,6)$, $(2,4,4)$ or $(3,3,3)$.

In \cite {Pi2, PS2} the authors study  functions  $f \bar g$ restricted to a complex analytic surface $X$ with (at most) an isolated singularity; of course $X$ can be $\C^2$.  It is proved that these functions have the so-called Thom property, and this implies that one has a fibration of the Milnor-L\^e type  (\ref{Fibration Thm., version2}). Moreover, they prove that in this case one has a Milnor fibration:
\begin{equationth}\label{PS2}
f\bar g /|f \bar g| :  \mathcal{L}_X  \setminus L_{f \bar g}
\to \s^1 \,. \end{equationth}
 In \cite [Theorem 5.2] {PS2} they further look at the geometry of the fibration 
near the multilink $L_{f \bar g}$ in $\mathcal{L}_X $ and 
they prove the equivalence of the following three statements (we refer to \cite{Ei-Ne} for background on fibered multilinks):

 %\begin{description}

\nn {\bf (i)}
 The real analytic germ $f \bar  g : (X,p) \to (\mathbb R^2,0)$ has $0$ as an isolated critical value;

\nn {\bf (ii)}
 the multilink $L_f -L_g$ is fibered; and 

\nn {\bf (iii)}
 if $\pi : \tilde{X} \to X$ is a resolution of the holomorphic germ $fg: (X,p) \to ({\mathbb C},0)$, then 
for each rupture vertex $(j)$ of the decorated dual graph of $\pi$ one has 
$m_j^f \neq  m_j^g$.

\noindent
Moreover,  it is proved that if these conditions hold, then the Milnor fibration in (\ref{PS2}) actually is a fibration of the multilink $L_f -L_g$ (cf. \cite{Ei-Ne}). This extends previous results in \cite{Pi2} for maps $\C^2 \buildrel{ f \bar g} \over \to \C$ to map-germs defined on surface singularities.

\vskip.2cm

The Thom property is used in this setting to grant that the germ $f \bar g$ has the Transversality Property \ref{transversality property} and therefore it has a Milnor-L\^e fibration. By \cite{PS2}  this holds for all maps $X \buildrel{ f \bar g} \over \to \C$ where $X$ is a complex surface which is either regular or has a normal singularity. The same proof was generalized in \cite{PS3} to higher dimensions, but this is mistaken and the erratum to that paper gives counterexamples. Yet, all the known counterexamples still satisfy the Transversality Property (see \cite{Oka4,  MS2}). This suggests:

\begin{question}
 Do there exist germs $\C^n \buildrel{ f \bar g} \over \to \C$ which do not have the Transversality Property \ref{transversality property}? Can we classify them?
\end{question}

%%%%% meromorphic case
We now look at meromorphic germs. 
Let  $U$ be an open neighborhood of $\0$ in $\C^n$ and let $f,g :  U \longrightarrow \C$
be two holomorphic functions without common factors such that $f(\0) = g(\0)=0$.

Let us consider the meromorphic function $F = f/g : U \rightarrow
{\mathbb C P}^1 $ defined by $(f/g) (x) =  [f(x)/ g(x)]$. As in
\cite{GLM1}, two such germs   at $\0$, $F= f/g$ and $F'  = f'/g'$,
are considered as equal (or equivalent) if and only if $f = h f'$
and $g = h  g'$ for some holomorphic germ $h:\C^n \to \C$ such
that $h(\0) \neq 0$. Notice that $f/g$ is not defined on the whole
$U$; its indetermination locus is
$$I= \big\{ z \in U \mid f(x)=0 \text{ and } g(x)=0 \big\} \,.$$
In particular,  the fibers of $F=f/g$ do not contain any point of $I$: for each $c \in \C$,
the fiber $F^{-1}(c)$  is:
$$F^{-1}(c)= \big\{x \in U  \mid f(x)-cg(x) = 0\big\} \setminus I \,.$$

In a series of articles, Gusein-Zade, Luengo and 
Melle-Hern\'andez studied
local Milnor-L\^e type fibrations in Milnor tubes  $f/g = |c|$
associated to every critical value of the meromorphic map $F =
f/g$. See for instance \cite{GLM1, GLM2}.  Of course these Milnor tubes 
 are in fact ``pinched
tubes'' that have $\0$ in their closure.

It is thus natural to ask whether one has for meromorphic
map-germs  Milnor fibrations on spheres, and if
so, how these are related to those of the Milnor-L\^e type, in Milnor tubes. The first of these questions was addressed in
\cite{BP, Pi2,  PS2} from two different viewpoints, while the answer to
the second question is the bulk of \cite{BPS} where 
 the authors   compare  the local
 fibrations in Milnor tubes of a meromorphic germ $f/g$, with the Milnor fibration
in the sphere.  They prove that if the germ $f/g$ satisfies two technical conditions (it is semitame
and (IND)-tame), then the  Milnor fibration
 for $f/g$ in the sphere is obtained from the Milnor-L\^e
fibrations of $f$ and $g$ in local  tubes at $\0$ and $\infty$, by a gluing process that is a 
``fiberwise" reminiscent of the classical connected sum of manifolds.

 \begin{remark}
 We know that if $(V,p)$ is a normal isolated complex surface singularity, then its link $L_V$ is a Waldhausen manifold, and there is a rich interaction between 3-manifolds theory and complex singularities  (see Section \ref {sec. relations}). When the germ $(V,p)$ is defined by a single equation $f: \C^3 \to \C$, then $L_V$ can be regarded as the boundary of the Milnor fiber $F_t = f^{-1}(t) \cap \buildrel{\circ}\over {\B}_\e$. If we now consider a holomorphic map-germ $\C^3 \to \C$ with a non-isolated critical point, then the corresponding link is no longer smooth, but we still have the Milnor fiber $F_t = f^{-1}(t) \cap \buildrel{\circ}\over {\B}_\e$, that can be regarded as a compact 4-manifold by attaching to it its boundary $\partial F_t = f^{-1}(t) \cap {\S}_\e$. It is proved in \cite {MPi1, MPi2, Ne-Sz} that in this setting, the boundary  of the Milnor also is a Waldhausen 3-manifold.  By  \cite{FM}, this statement actually extends to the boundary of the Milnor fibers of all map-germs in $\C^3$ of the form $f \bar g$ that have the transversality condition, which includes the holomorphic germs. 
 
 \end{remark}

 %%%%%%%%%%%%
 %%%%%%%%%%%

 \subsection {A glance at Oka's work on mixed functions}
  
 Inspired by the theory of complex singularities, 
 Oka in \cite {Oka2} introduces for mixed functions the useful notion of non-degeneracy with respect to a naturally defined Newton boundary.  He uses this to prove a   fibration theorem for strongly non-degenerate convenient mixed functions, and to study their topology. We say a few words about this.
 
  Let $f(z)$ be a mixed analytic function of the form
  $$f(z) = \sum_{\nu,\mu} c_{{\nu,\mu}} z^\nu \bar z^\mu \,, $$
  where  $\nu  + \mu$  is the sum of the multi-indices of $z^\nu \bar z^\mu$, {\it i.e.}, 
  $$\nu  + \mu = (\nu_1  + \mu_1,...,\nu_n  + \mu_n).$$
  Assume for simplicity  $c_{0,0} = 0$, so that $\0 \in V(f):= f^{-1}(0)$. Following Oka, we call $V(f)$ a mixed hypersurface, though in general it has real codimension 2. 
  The {\it (radial) Newton polygon (at the origin)} $\Gamma_+(f)$ is defined in the usual way: it is the convex hull of: 
  $$\bigcup_{c_{\nu,\mu} \ne 0} (\nu + \mu) + \R^{+n} \;.
  $$
In analogy with complex polynomials, define 
the {\it Newton boundary} $\Gamma(f)$ as the union of the compact faces of $\Gamma_+(f)$. 
To every given positive integer (weight) vector $P = (p_1,...,p_n)$ we associate a linear function $\ell_P$ on the Newton boundary $\Gamma(f)$ defined by:
$$\ell_P(\nu) = \sum_{j=1}^n p_j \nu_j$$
for $\nu \in \Gamma(f)$. Let $\Delta(P,f) = \Delta(P)$ be the face where $\ell_P$ attains its minimal value. Then, 
for a positive weight $P$, define the  face function $f_P(z)$ by:
$$f_P(z) = \, \sum_{\nu + \mu \in \Delta (P)} c_{\nu,\mu} z^\nu \bar z^\mu \;.
$$

\begin{definition} Let $P$ be a strictly positive weight vector. We say that $f(z)$ is {\it non-degenerate} for $P$ if the fiber $f^{-1}(0) \cap \C^{*n}$ contains no critical point
of the map $\C^{*n} \buildrel {f_P } \over {\to} \C$. The map $f$ is 
{\it strongly non-degenerate} for $P$ if the mapping  $\C^{*n} \buildrel {f_P } \over {\to} \C$ has no critical points at all, ${\rm dim} \,\Delta(P ) \ge 1$  and $f_P : \C^{*n}  \to \C$ is surjective. The 
function $f(z)$ is called non-degenerate (respectively strongly non-degenerate) if it is non-degenerate (respectively strongly non-degenerate) for every strictly positive weight vector $P$.
\end{definition}

%%%%%%% Convenient
For a subset $J \subset  \{1, 2, . . . , n\}$, we consider the subspace $\C^J$ and the restriction $f^J := f|_{\C^J}$ . Consider the set
$${\mathcal N \mathcal V}(f) = \{I \subset \{1,...,n\} \,| \,  f^I \ne 0  \}.$$
We call ${\mathcal N \mathcal V}(f)$ the set of non-vanishing coordinate subspaces for $f$.

\begin{definition}
We say that $f$ is k-convenient if $J \in  {\mathcal N \mathcal V}(f)$ for every  $J \subset  \{1,...,n\}$ with $|J| = n -k$. We say that $f$ is convenient if $f$ is (n-1)-convenient.
\end{definition}

We may now state the main result in \cite{Oka2} concerning Milnor fibrations. In \cite{Oka2} this is stated as theorems 29, 33 and 36. Here we combine them into a single statement.

\begin{theorem}
Assume the mixed polynomial 
 $f(z)$ is  convenient and strongly non-degenerate. Then one has a fibration of the Milnor-L\^e type in a Milnor tube as in (\ref{Fibration Thm., version2}), as well as a Milnor fibration on every sufficiently small sphere with projection map $f/\Vert f \Vert$, as in (\ref{Fibration Thm., version1}), and the two fibrations are smoothly equivalent.
\end{theorem}

%%%%%%
%%%%%%
%%%%%%
%%%%%%

In \cite {Oka2} Oka also uses toric geometry to get a resolution of the corresponding singularity, in analogy with the complex case (see for instance the  book \cite{Oka0}). He then uses this to study the topology of the links, as well as the topology of the Milnor fibers. 
Oka has subsequently studied and published various articles on the subject. We list some of them (\cite {Oka1, Oka2, Oka3, Oka4}), but there are several more,  with important results that cover a wide spectrum of topics, from intersection theory to contact structures.

\section{\bf Linear actions, intersections of quadrics and  LVM manifolds}

The motivation in \cite {Se5} for  studying the twisted Pham-Brieskorn singularities 
comes from the pioneering work \cite {LM1} by S. L\'opez de Medrano on intersections of quadrics and the space of Siegel leaves of holomorphic linear flows. In fact, notice that if 
$\xi$ is the holomorphic vector field in $\C^n$ defined by:
$$\xi(z) = (\lambda_1 z_{\sigma(1)}^{a_1} , \cdots, \lambda_n z_{\sigma(n)}^{a_n}) \quad; \; \lambda_i \in \C^*\,,
$$
where $\sigma$ is a permutation of the set $\{1,\cdots,n\}$, then  the set of points where the Hermitian product $\langle \xi(z), z \rangle$ vanishes is the mixed variety in $\C^n$ defined by:
\begin{equationth}\label{PH-variety}
V_\xi \; := \; \{ \lambda_1 z_{\sigma(1)}^{a_1} \bar z_{1} + ... + \lambda_n z_{\sigma(n)}^{a_n} \bar z_{n} 
 \; = \; 0 \} \,.
 \end{equationth}
\hskip-4pt Re-labelling the variables $(z_1,\cdots,z_n)$, and if we assume all  $a_i> 1$,  we arrive to the twisted Pham-Brieskorn polynomials (up to the $\lambda_i$s). Notice that the variety defined in (\ref{PH-variety}) is the set of points where the leaves (or solutions) of the vector field $\xi$ are tangent to the foliation in $\C^n$ by all the spheres centered at $\0$, union $\0$ itself.

If the exponents $a_i$ are all 1, the vector field is linear. And if the twisting $\sigma$ is the identity, then $\xi$ is a linear diagonal vector field with eigenvalues $\lambda_1,\cdots,\lambda_n$. In this case one has:
$$\langle \xi(z), z \rangle \,= \, \sum_{i=1}^n \lambda_i \vert z_i \vert^2 \,,$$ 
and
the variety  $V_\xi$ in (\ref{PH-variety}) is the intersection of the two real quadrics below, corresponding to the real and the imaginary parts:
$$ {\rm Re} (\sum_{i=1}^n \lambda_i \vert z_i \vert^2) = 0 \qquad \hbox{and} \qquad  {\rm Im}(\sum_{i=1}^n \lambda_i \vert z_i \vert^2) = 0 \,.
$$
We say that $\xi$ is in the Siegel domain if the convex hull of the $\lambda_i$ in $\C$ contains the origin $\0$. Otherwise $\xi$ is in the Poincar\'e domain.

\begin{figure}
\centering
\includegraphics[height=7cm ]{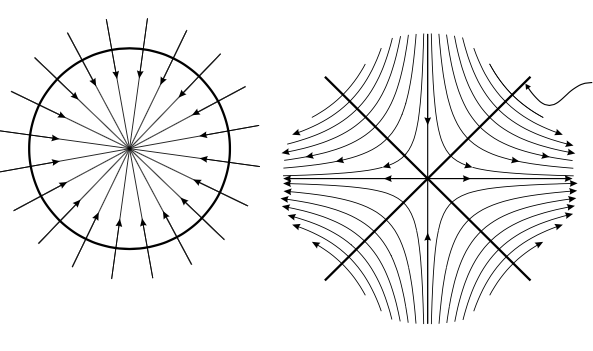}
\caption{Linear vector fields in the Poincar\'e and Siegel domains, respectively. The curled arrow points into the space of Siegel leaves.}
\end{figure}

When $\xi$ is in the Poincar\'e domain, it is easy to show that the variety  $V_\xi$ actually consists only of the point $\0$. 
So we will
restrict from now on to linear vector fields in the Siegel domain.

It is  known that the equality, and even real dependence, of two eigenvalues of $F$ complicates the topology of $\F$ and $V_F$ very much. 
Therefore one usually assumes the following generic {\it hyperbolicity hypothesis}: any two eigenvalues are independent over $\R$:
\begin{equationth}\label{hyperbolicity}
i \, \ne \, j \, \Rightarrow \lambda_i \notin \R \lambda_j\,, \quad \forall \; \; i,j = 1,\cdots,n \,,  
\end{equationth} \hskip-4pt
So we now let $F$ be a linear vector field in the Siegel domain satisfying (\ref{hyperbolicity}). It is clear that a point $z \in \C^n-0$ is 
in $V_F$ iff the restriction of the real function 
$d(z) = \|z\|^2 = \sum_{i=1}^n |z_i|^2$ to the leaf ${\mathcal L}_z$ through $z$ has a critical point at $z$. Furthermore,
as noted in~\cite[ Section 3]{CKP}, the fact that the solutions of $F$ are parametrized by exponential maps implies that the 
{\it leaves of $F$ are concave}. Thus, if a leaf ${\mathcal L}$ meets $V_F^*$, then it has a unique point in $V_F$ and it is the point in
${\mathcal L}$ of minimal distance to $0 \in \C^n$. Such a leaf is called {\it a Siegel leaf}. It is a  copy of $\C$ embedded in $\C^n$ and can be 
characterized by its unique point in $V_F$. Furthermore, the fact that the 
intersection
${\mathcal L} \cap V_F$ of each leaf that meets $V_F$ is at a local minimal point in 
${\mathcal L}$, implies that 
${\mathcal L} \cap V_F$ is a transverse intersection.
 By the flow-box theorem for complex
differential equations, this means that we have  at each $z \in V_F^*$ a neighborhood of the form $U_z \times \D^2$, where $U_z$ is a disc
of real dimension $2n-2$ and the second factor denotes small discs in the leaves. It follows that $V_F^*$ is a smooth real
submanifold of $\C^n$ of codimension 2. In fact $V_F$ is a real analytic complete intersection in $\C^n$ and 
 the union $W = V_F^* \times \C$  of all the Siegel leaves of $F$ is an open subset of $\C^n$ that can be identified with the total space of
 the normal bundle of $V_F^*$.

It is shown in~\cite {CKP} that $W$ is actually dense in  $\C^n$. It is an exercise to see that $V_F$ is globally embedded as a cone with vertex at $0$ and base the 
intersection $M = V_F \cap \s^{2n-1}$ with the unit sphere, which is the link of the corresponding singularity.

We summarize part of the above discussion in the following theorem, which  re-phrases results in ~\cite {CKP}:

\begin{theorem}\label{thm. LM}  
Let $F(z) = (\lambda_1 z_1, \cdots, \lambda_n z_n)$ be a linear vector field in the Siegel domain 
satisfying the hyperbolicity condition (\ref{hyperbolicity}). Then the real analytic variety:
\[V_F \,=\, \{\,z \in \C^n\,|\, \hbox{\rm Re}\,(\sum_{i=1}^n \lambda_i |z_i|^2 )\,= \,0\,\} \,\cap\,
\{\,z \in \C^n\,|\, \hbox{\rm Im}\,(\sum_{i=1}^n \lambda_i |z_i|^2 )\,= \,0\,\}\;,\]
is a complete intersection of real codimension 2 with an isolated singularity at $0$, and the regular points of $V_F$ parametrize the Siegel leaves of $F$. 
\end{theorem}

Notice that one has the mixed function   defined by $\big(\hbox{Re}\,(\sum_{i=1}^n \lambda_i |z_i|^2), 
\hbox{Im}\,(\sum_{i=1}^n \lambda_i |z_i|^2 )\,\big)$,  that determines the complete intersection germ in Theorem \ref{thm. LM}. It is easy to see that its discriminant $\Delta_F$ consists of the $n$ half-lines in $\C$ determined by the eigenvalues $\lambda_i$. These split $\C$ into $n$-sectors and just as in Example \ref {ex. Santiago},  one has a Milnor-L\^e fibration over each sector. The topology of the  Milnor fibers over the various sectors is determined in \cite{LopezdeMedrano:SHQM}.   In \cite {LM:singular} it is shown that the fibers over  points in the discriminant necessarily are singular. This is related  with 
 a problem studied by C. T. C. Wall
 ~\cite {Wall} in 
relation with the topological stability of smooth mappings.

The topology of the link of these singularities is studied in \cite{CKP} when $n=3$ and  in \cite{LM1, LM2} for $n>3$. For  $n=3$ the link is always a 3-torus $\S^1 \times \S^1 \times \S^1$. For $n>3$ the link is always a connected sum of products of spheres determined by the geometry of the polytope spanned by the eigenvalues of $F$. 

These constructions extend to sets of $m$  commuting vector fields on $\C^n$, $m \ll n$. This is studied in 
\cite {Me, MV1,  LV}.  A remarkable point is that in all these cases, the manifold $V_F^* := V_F \setminus \{0\}$ has a canonical complex structure. In fact, under the appropriate hyperbolicity condition,  $V_F^*$  is a smooth submanifold of $\C^n$ of real codimension $2m$; the leaves of $\mathcal F$ are transversal to $V_F^*$ everywhere and this equips $V_F^*$ with a canonical complex structure, by a theorem of Haefliger in \cite{Haef}. We remark that $V_F^*$ is not embedded in $\C^n$ as a complex submanifold.

This all is very surprising: for instance, consider the  case $n=3$ and $m=1$. One has that $V_F$ is a real analytic singularity with a canonical complex structure away from $0$ and the link is  $\S^1 \times \S^1 \times \S^1$. Then $V_F$ is a real analytic complete intersection with a canonical complex structure on its regular part  $V_F^*$, but   $V_F$  cannot be complex at $0$ because the 3-torus cannot be the link of a complex singularity, by \cite{Sul1}. 

The complex manifolds $V_F^*$ one gets in this way are equipped with a canonical $\C^*$-action and the quotient $V_F^*/\C^*$ is a compact complex manifold. These are known as LVM-manifolds (see \cite {MV2}) and they are a special class of the so-called moment-angle manifolds, with remarkable geometric and topological properties (cf. \cite{Gi-LM, BaLoVe}).

It is also interesting to
 determine how the topology of the singularities defined as above varies as we pass
from one component in the Siegel domain to another one, {\it i.e.}, as we break the weak hyperbolicity condition. This is beautifully 
answered  in  \cite {Bo-Me}. This is a ``wall-crossing'' problem as the authors explain,
 and they show that crossing a wall means performing a precise surgery, which they describe. From the viewpoint of singularities, what we do 
when crossing a wall is to put two complete intersection singularities in a 1-parameter family of singularities which are all 
complete intersections, except that they bifurcate when crossing the wall.

\section*{}

%%%%%%%%%%%%%
%%%%%%%%%%%
%%%%%%%%%%%% 

\end{document}